\documentclass[letter,11pt]{amsart}

\usepackage{times}
\usepackage[pdftex]{graphicx}
\usepackage{amssymb,amsfonts,amsmath,amsthm}
\usepackage{float}
\usepackage{xcolor}
\usepackage{enumitem,comment}
\newtheorem{theorem}{Theorem}
\newtheorem{proposition}{Proposition}
\newtheorem{lemma}{Lemma}
\newtheorem{ass}{Assumption}
\newtheorem{definition}{Definition}
\newtheorem{remark}{Remark}
\newtheorem{cor}{Corollary}

\newcommand\EE {\mathbb E}

\newcommand\NN {\mathbb N}
\newcommand\RR {\mathbb R}
\newcommand\PP {\mathbb P}

\newcommand\ZZ {\mathbb Z}

\newcommand\cF {\mathcal F}
\newcommand\cL {\mathcal L}

\def\bone{\mathbf{1}}
\def\pa{\partial}

\newcommand\1 {\mathbf 1}

\newcommand{\ed}{\end{document}}
\newcommand{\be}{\begin{equation}}
\newcommand{\ee}{\end{equation}}
\newcommand{\bq}{\begin{eqnarray}}
\newcommand{\eq}{\end{eqnarray}}

\vspace{4in}

\definecolor{Red}{rgb}{0.9,0,0.0}
\definecolor{Blue}{rgb}{0,0.0,1.0}

\begin{document}
\title[Consistency of MLE in partially observed diffusions]{Consistency of MLE in partially observed diffusion models on a torus} 
\author{Ibrahim Ekren and Sergey Nadtochiy}
\address{Mathematics Department, University of Michigan, Ann Arbor, MI 48109.} 
\email{iekren@umich.edu}
\address{Department of Applied Mathematics, Illinois Institute of Technology, Chicago, IL 60616.}
\email{snadtochiy@iit.edu}
\footnotetext[1]{S.~Nadtochiy is partially supported by the NSF grant DMS-2205751.}
\footnotetext[2]{I.~Ekren is partially supported by the NSF grant DMS-2406240.}

\begin{abstract}
In this paper, we consider a general partially observed diffusion model with periodic coefficients and with non-degenerate diffusion component. The coefficients of such a model depend on an unknown (static and deterministic) parameter which needs to be estimated based on the observed component of the diffusion process. We show that, given enough regularity of the diffusion coefficients, a maximum likelihood estimator of the unknown parameter converges to the true parameter value as the sample size grows to infinity.
\end{abstract}

\maketitle

\section{Introduction and main results}


In this article, we consider the unobserved $\RR^q$-valued signal $X$ following the diffusion process
\begin{align}
& dX_t = b^\theta(X_t) dt + \sigma^\theta(X_t)\,dB_t,\quad t\geq0,\quad X_0\sim \nu,\label{eq.Background.dyn.X}
\end{align}
and the observation process $Y$ given by
\begin{align}
& dY_t = h^\theta(X_t) dt + dW_t,\quad t\geq0,\quad Y_0=0,\label{eq.Background.dyn.Y}
\end{align}
where $B$ and $W$ are standard Brownian motions, taking their respective values in $\RR^{d}$ and $\RR^m$,
$X_0$ is independent of $B$, and $W$ is independent of $X$.
The coefficients of the above diffusion processes depend on the unknown parameter $\theta$ taking values in a metric space $(\Theta,\bar d)$.\footnote{Note that the generality of the parameter space is not a key contribution of the present work: e.g., even if $\Theta$ were a subset of a Euclidean space, the arguments used in the proof of our main result would not simplify.}
Any statistical inference about the value of the unknown parameter $\theta$ is done using only an observed path of $Y$. The structure of \eqref{eq.Background.dyn.Y} is interpreted as follows: the observation acquired at time $t$ is given by $h^\theta(X_t)\,dt$ corrupted by the independent noise $dW_t$. The latter structure is standard in the theory of (continuous-time and continuous-space) stochastic filtering: cf. \cite{BainCrisan} for the discussion of such models and their applications, and \cite{Elliott}, \cite{HMM} for the versions of such models where either the time or the signal space are discrete.

\smallskip

The goal of the present work is to establish consistency of the maximum likelihood estimator (MLE) of the unknown parameter $\theta$, as the ``sample size" $t$ of the observation $Y_{[0,t]}$ goes to infinity. The latter question has been analyzed in the context of discrete-time partially observed (hidden) Markov models, e.g., by \cite{MLELeroux}, \cite{MLEDouc}, \cite{MLEGenonCatalot}, \cite{MLERamon}, \cite{HMM}, which ultimately show that the consistency of MLE for $\theta$ follows from the ergodicity of the signal and from the continuity of the model with respect to $\theta$. The case of continuous time and discrete signal space was considered in \cite{MLEChiganski}, and the consistency of MLE-type estimators in specific partially observed diffusion models was established in \cite{MLEKutoyants}, \cite{MLEKhasminski}. To the best of our knowledge, the present work provides the first proof of consistency of MLE in general (i.e., with general coefficients $(b,\sigma,h)$) partially observed diffusion models.

The main challenge of the analysis of MLE in partially observed Markov models is that the associated likelihood function does not have a sufficiently tractable representation, making it impossible to apply the methods based on explicit computation, described, e.g., in \cite{KutoyantsBook}, \cite{IbrahimovKhasminskiBook}. Instead, we follow the approach outlined, e.g., in \cite{MLEChiganski} and \cite{MLEDouc}:
 (i) establish the exponential stability of the conditional distribution of the signal given an observation, known as the \emph{stochastic filter}, with respect to its initial condition (this stability is essentially known, even in the diffusion setting), (ii) show the uniform robustness of the filter with respect to the unknown parameter (this is a new result, in the diffusion setting), and (iii) prove that any two distinct distributions of the observation, restricted to $[0,t]$, separate exponentially fast as $t\rightarrow\infty$. On a technical level, the proofs for discrete-time hidden Markov models rely on the stationarity/ergodicity of the observation process and on the associated ergodic theorems in order to deduce (iii) from the ergodicity of the signal and from (i). In the continuous-time models, the observation process is typically non-stationary.\footnote{The continuous-time setting can be viewed as a discrete-time model with a stationary observation process and with a signal process taking values in an infinite-dimensional space. However, such an infinote-dimensional representation yields additional technical challenges, e.g., in verifying the ``mixing" property of the signal.} Nevertheless, the special form of the likelihood function in the diffusion setting makes it very convenient to apply martingale methods in order to obtain a somewhat simpler proof of (iii). In particular, we manage to avoid making any assumptions on the ergodicity of the filter and on its invariant distribution, and obtain the results analogous to those established for discrete-time models: i.e., the consistency of MLE follows from the ergodicity of the signal and from the continuity of the model with respect to the parameter.




\smallskip

\begin{remark}
In some applications it may be natural to consider the observation function $h^\theta$ that depends on both $X_t$ and $Y_t$, and to allow $W$ to be correlated with $B$. Although, strictly speaking, we exclude such settings herein, it is worth mentioning that the former is a straightforward extension of the results and methods described in this article (provided $h^\theta$ is periodic in both variables), while the latter requires additional work. It is also easy to see that a simple linear transformation applied to $Y$ extends our results to the case where `$dW_t$' is replaced by `$\Sigma\,dW_t$', with an invertible constant matrix $\Sigma$. The case of $\Sigma$ depending on $Y_t$ can be handled easily as well (provided $\Sigma$ is periodic), while the case of $\Sigma$ depending on $X_t$ and the case of degenerate $\Sigma$ are fundamentally different from the setting considered herein.
\end{remark}

\medskip

Before introducing the filter, the likelihood function and MLE, we need to make several assumptions on the input $(b,\sigma,h)$. 

\begin{ass}\label{periodic}
We assume that the functions $b,\sigma, h$ are $1$-periodic, in the sense that they are invariant with respect to the translations of their variable by $(n_1,\ldots,n_q)^\top$, for any $n_1,\ldots,n_q\in\ZZ$.
\end{ass}

The above assumption allows us to interpret the signal $X$ as a process on a torus, which is needed in order to establish the ergodicity of $X$ and to prove the exponential stability of the filter. Indeed, to date, there exist no results on the exponential stability of a stochastic filter in diffusion models with unbounded domains, which would not impose other restrictive assumptions on the model coefficients, such as the linearity of $h$, the requirement that $\sigma(x)$ does not depend on $x$, and additional assumptions on the initial value of the filter. We choose to assume the periodicity of the coefficients instead (see, e.g., \cite{ChiganskyStabSurvey}, \cite{ChiganskyRamonStab}, \cite{Stannat}, \cite{RamonThesis}, \cite{RamonStab} and the references therein for more on filter stability).

For any vector $e:=(e_1,\ldots,e_q)$ with strictly positive entries, the results established herein extend trivially to a setting where the coefficients are $e$-periodic (i.e., invariant with respect to the translations of their variable by $(n_1e_1,\ldots,n_qe_q)^\top$, for any $n_1,\ldots,n_q\in\ZZ$).


\smallskip

Next, we make the following technical assumption.

\begin{ass}\label{holder}
We define $a^\theta:=\sigma^\theta(\sigma^\theta)^\top$ and assume that the functions $(a,b,h):\Theta\times \RR^q\mapsto \RR^{q\times q}\times\RR^q\times\RR^m$ satisfy
\begin{align*}
\inf_{|\xi|=1,\,\theta\in\Theta,\,x\in\RR^q} \xi^\top\, a^\theta(x)\,\xi>0\,\,\mbox{ and }\,\,\sup_{\theta\in\Theta}\left(\|a^\theta\|_{C^1}+\|b^\theta\|_{C^1}+\|h^\theta\|_{C^{2}}\right)<\infty.
\end{align*}
\end{ass}

Under Assumption \ref{holder}, thanks to \cite[Theorem 7.2.1]{stroock1997multidimensional}, for any $\nu\in\mathcal{P}(\RR^q)$ and any $\theta\in\Theta$, there exists a weak solution to \eqref{eq.Background.dyn.X} that is unique in law and satisfies the strong Markov property. Herein and in the remainder of the paper, we denote by $\mathcal{P}(\Omega)$ the space of probability measures on $\Omega$, equipped with the Borel sigma-algebra if $\Omega$ is a topological space and unless stated otherwise. We denote by $\PP^{\theta,\nu}$ the law of $(X,Y)$, viewed as a random element in $C(\RR_+\rightarrow\RR^q\times\RR^m)$ equipped with the cylindrical sigma-algebra.

\begin{remark}
The regularity of $(b^\theta,\sigma^\theta)$, stated in Assumption \ref{holder}, is needed for the following reasons: to ensure the existence of a weak solution to \eqref{eq.Background.dyn.X}, as well as its uniqueness in law, Markov property and continuity (in law) with respect to the initial distribution, and to prove Lemmas \ref{ass:mixing}, \ref{le:v.reg} and \ref{le:expStab.X}. In principle, for a specific diffusion model at hand, the latter results may be established under a relaxed set of assumptions: e.g., one may be able to deduce Lemma \ref{ass:mixing} for discontinuous $(b^\theta,\sigma^\theta)$ using the methods of \cite{KrylovNonlin}, while the well-posedness of \eqref{eq.Background.dyn.X} and of the SDE in the proof of Lemma \ref{le:expStab.X} does not require any regularity of $b^\theta$. The ellipticity of $a^\theta$ stated in Assumption \ref{holder} is crucial for the approach chosen herein: in particular, along with Assumption \ref{periodic}, it ensures the ergodicity of $X$.

The $C^2$ property of $h^\theta$ is used to establish the representation \eqref{eq:rep1} in Lemma \ref{le:rep.via.U}, which is important for our analysis. An alternative representation can be established using the results of \cite{PardouxTrick}, without any regularity of $h^\theta$. The latter representation is, unfortunately, less convenient to work with, in particular, making it much more challenging to prove Lemma \ref{le:v.reg} -- this is why we chose to use \eqref{eq:rep1} instead and, hence, to require the $C^2$ property of $h^\theta$ stated in Assumption \ref{holder}.
\end{remark}

\smallskip

Assumptions \ref{periodic} and \ref{holder} imply the existence of a sufficiently regular invariant density for the signal $X$ viewed as a process on the associated torus: i.e., for the process 
\begin{align*}
X^e_t:=(X^1_t\,\mathrm{mod}\,1,\ldots, X^q_t\,\mathrm{mod}\,1)^\top.
\end{align*} 

\begin{lemma}\label{ass:regular.stationary.density}
For any $\theta\in\Theta$, there exists a bounded away from zero and continuously differentiable $\psi^{\theta}_0:\,{[0,1)^q}\rightarrow\RR_+$, with $\int_{[0,1)^q} \psi^\theta_0(x)\,dx=1$, such that the distribution of $(X^e_{t+\cdot})_{\RR_+}$ under $\PP^{\theta,\nu^\theta_0}$ does not depend on $t\geq0$, where $\nu^\theta_0(dx):=\psi^\theta_0(x)\,dx$. 
\end{lemma}
The proof of the above lemma is standard, and it is given in Subsection \ref{subse:regular.stationary.density}.

\smallskip

In order to have any hope for a consistent estimation of $\theta$, we need to ensure that the coefficients of the model are continuous with respect to this parameter.

\begin{ass}\label{ass:unifcont}
We assume that $\Theta$ is compact and that the following uniform continuity property holds:
\begin{align*}
\lim_{\delta \downarrow 0}\sup_{\bar d(\theta,\theta')\leq \delta} \left(\|a^\theta-a^{\theta'}\|_{C^1} + \|b^\theta-b^{\theta'}\|_{C} + \|h^\theta-h^{\theta'}\|_{C^{2}}\right)=0.
\end{align*}
\end{ass}

\medskip

Next, we recall that $(X,Y)$ is the canonical element on $C(\RR_+\rightarrow\RR^q\times\RR^m)$ and consider the natural filtrations $\mathbb{F}^X$, $\mathbb{F}^Y$ and $\mathbb{F}:=\mathbb{F}^{X,Y}=(\mathcal{F}^X_t\bigvee\mathcal{F}^Y_t)_{{t\geq0}}$. We also recall that $\PP^{\theta,\nu}$ is a measure on $C(\RR_+\rightarrow\RR^q\times\RR^m)$, given by the unique law of the solution to \eqref{eq.Background.dyn.X}--\eqref{eq.Background.dyn.Y}, and denote by $\PP^{\theta,\nu}_Y$ the $Y$-marginal of $\PP^{\theta,\nu}$, defined on $(C(\RR_+\rightarrow\RR^m),\mathcal{F}^Y_\infty:=\sigma(\bigcup_{t\geq0}\mathcal{F}^Y_t))$. We denote the restriction of $\PP^{\theta,\nu}_Y$ to $\mathcal{F}^Y_t$ by $\PP^{\theta,\nu}_{Y,t}$. Recall that, in a model where the observation is given by $Y_{[0,t]}$, the time horizon $t$ represents the ``size of a sample", and $t\rightarrow\infty$ corresponds to the large-sample regime. For a fixed $t<\infty$, all $\{\PP^{\theta,\nu}_{Y,t}\}_\theta$ are equivalent -- to each other and to the Wiener measure $\mathbb{W}_{Y,t}$ -- but not for $t=\infty$. In fact, these measures are pairwise singular, for $t=\infty$ and with $\nu=\nu^\theta_0$, provided that $\theta$ determines $\PP^{\theta,\nu^\theta_0}_Y$ uniquely. The latter property is known as the \emph{identifiability} of $\theta$ (cf. \cite{MLELeroux}, \cite{IDsurvey}). 

To cover a more general scenario, in which the identifiability may fail, for each parameter value $\theta\in\Theta$, we introduce the equivalence class
\begin{align}
\overline{\Theta}(\theta):=\{\theta'\in\Theta:\,\PP^{\theta,\nu^\theta_0}_Y=\PP^{\theta',\nu^{\theta'}_0}_Y\}.\label{eq.theta.equivClass}
\end{align}
Note that all elements of $\overline{\Theta}(\theta)$ are statistically indistinguishable, as they consist of the parameter values that generate the same law of the observation: indeed, one would not be able to differentiate any two sample paths produced by the models with different parameter values in $\overline{\Theta}(\theta)$ from two sample paths generated by the $\theta$-model. Hence, any statistical inference for $\theta$ can only be done up to its equivalence class.

\begin{remark}
The use of $\nu^\theta_0$ as the initial condition, in the above definition of $\overline{\Theta}(\theta)$, is based on the implicit assumption that the signal process $X$ had been running for an ``infinite amount of time" before the observation started, so that it had reached its stationary distribution by the time $t=0$.
In view of Lemma \ref{le:expStab.X}, one can easily modify \eqref{eq.theta.equivClass} by allowing for non-stationary initial distributions of the signal, with the statement of Theorem \ref{thm:main} modified accordingly.
\end{remark}

\smallskip

To define the likelihood function and its maximizer, we introduce the (stochastic) filter 
\begin{align*}
\pi^{\theta,\nu}_t(dx):=\PP^{\theta,\nu}(X_t\in dx\,\vert\,\mathcal{F}^Y_{t})
\end{align*} 
and view it as a $\mathcal{P}(\RR^q)$-valued random element.
We recall the classical optional projection result (cf. Theorem 7.12 in \cite{LiptzerShiryaev}):
\begin{align}
dY_t = \pi^{\theta,\nu}_t[h^\theta]\,dt + d\tilde W^{\theta,\nu}_t,\quad \PP^{\theta,\nu}\text{-a.s.}\label{eq.Y.projected}
\end{align}
where $\pi^{\theta,\nu}_t[h^\theta]:=\int_{\RR^q}h^\theta(x)\,\pi^{\theta,\nu}_t(dx)=\EE^{\theta,\nu}\left(h^\theta(X_t)\,\vert\,\mathcal{F}^Y_t\right)$ and $\tilde W^{\theta,\nu}$ is a standard Brownian motion with respect to $(\mathbb{F}^Y,\PP^{\theta,\nu})$.
Then, using Girsanov's theorem (or further results in \cite{LiptzerShiryaev}), we deduce that the Radon-Nikodym derivative of $\PP^{\theta,\nu}_{Y,t}$ with respect to $\mathbb{W}_{Y,t}$ is given by
\begin{align}
& L^{\theta,\nu}_{t}:=\exp\left(-\frac{1}{2}\int_0^{t} |\pi^{\theta,\nu}_s[h^\theta]|^2\,ds + \int_0^t \pi^{\theta,\nu}_s[h^\theta]^\top dY_s\right),
\label{eq.L.def}
\end{align}
which is the likelihood of $Y_{[0,t]}$ with respect to the (reference) Wiener measure.

\begin{definition}
A $\Theta$-valued stochastic process $(\hat\theta^\nu_t)_{t\geq0}$, adapted to $\mathbb{F}^Y$, is a \emph{maximum likelihood estimator} (MLE) associated with the initial distribution $\nu$, if it satisfies the following, $\mathbb{W}_Y$-a.s., for any $t\geq0$:
\begin{align*}
& L^{\hat\theta^\nu_t,\nu}_{t} \geq L^{\theta',\nu}_{t}\quad\forall\,\theta'\in\Theta.
\end{align*}
\end{definition}

\begin{remark}
Since $\Theta$ is assumed to be compact and in view of the continuity of $\theta\mapsto L^{\theta,\nu}_t$, which is ensured, e.g., by Lemma \ref{le:v.reg}, the existence of MLE follows from the measurable selection theorem (cf. Theorem 18.13 in \cite{Aliprantis}), provided the metric space $\Theta$ is separable.
\end{remark}

\medskip

The following theorem is the main result of this paper: it states that MLE is consistent, in the sense that, under $\PP^{\theta,\nu^\theta_0}$, it converges in probability to the equivalence class of the true parameter $\overline{\Theta}(\theta)$.

\begin{theorem}\label{thm:main}
Under Assumptions \ref{periodic}, \ref{holder}, \ref{ass:unifcont}, any MLE $\hat\theta^\nu$ is consistent, in the sense that
\begin{align*}
\lim_{t\rightarrow\infty}\PP^{\theta,\nu^\theta_0}\left(\bar d(\hat\theta^\nu_t,\overline{\Theta}(\theta))\geq\varepsilon\right)=0,
\end{align*}
for any $\varepsilon>0$, $\nu\in\mathcal{P}(\RR^q)$ and $\theta\in\Theta$, where $\bar d(\theta',\overline{\Theta}(\theta)):=\inf_{\theta''\in\overline{\Theta}(\theta)}\bar d(\theta',\theta'')$.
\end{theorem}

If the mapping $\theta\mapsto\PP^{\theta,\nu^\theta_0}_Y$ is injective, then, the set $\overline{\Theta}(\theta)$ is a singleton and the above result yields a more standard form of consistency.

\begin{cor}
Under Assumptions \ref{periodic}, \ref{holder}, \ref{ass:unifcont}, and assuming the identifiability (i.e., that the mapping $\theta\mapsto\PP^{\theta,\nu^\theta_0}_Y$ is injective), the following holds for any MLE $\hat\theta^\nu$:
\begin{align*}
\lim_{t\rightarrow\infty}\PP^{\theta,\nu^\theta_0}\left(\bar d(\hat\theta^\nu_t,\theta)\geq\varepsilon\right)=0,
\end{align*}
for any $\varepsilon>0$, $\nu\in\mathcal{P}(\RR^q)$ and $\theta\in\Theta$.
\end{cor}

\smallskip

\begin{remark}
The recent paper \cite{MLESergey} provides a sufficient condition for the identifiability of $\theta$, which ultimately reduces to the verification that the mapping from $\theta$ to the stationary distribution of $h^\theta(X_\cdot)$ is injective. Note also that $h^\theta$ can be estimated with certainty from a finite sample using the Zakai equation, which means that the identifiability can be be ensured by checking that the mapping $\theta\mapsto h^\theta$ is injective. More generally, one can ensure identifiability by showing that the mapping from $\theta$ to the law of $(Y_{t_1},\ldots,Y_{t_n})$ under $\PP^{\theta,\nu^\theta_0}$ is injective, for some $0\leq t_1<\cdot<t_n$, by finding an appropriate functional that separates the distinct values of $\theta$, as shown in the proof of Lemma \ref{le:singularity} (see also \cite{IDsurvey}, \cite{IDalgebra}, \cite{IDNonParam}, \cite{MLELeroux} for the sufficient conditions for identifiability in hidden Markov models in discrete time).
\end{remark}



\medskip

In the course of proving Theorem \ref{thm:main}, we establish a new result on the uniform robustness of a stochastic filter in diffusion models, i.e., its continuity with respect to the parameter $\theta$ uniformly over the sample size (note that \cite{RobContSpaceTime} is not sufficiently strong for our purposes), which is valuable in its own right and is stated in the following proposition.

\begin{proposition}\label{le:uniform.robust.filter}
Under Assumptions \ref{periodic}, \ref{holder}, \ref{ass:unifcont}, for any $\theta\in\Theta$, $\nu\in\mathcal{P}([0,1)^q)$ and $\epsilon>0$,
\begin{align}
&\lim_{\delta\downarrow0}\sup_{t\geq1}\PP^{\theta,\nu^\theta_0}\left(\sup_{ \bar d(\theta',\theta'')\leq\delta} |\pi^{\theta',\nu}_t[h^{\theta'}] - \pi^{\theta'',\nu}_t[h^{\theta''}]|\geq\epsilon\right)=0.
\label{eq.uniform.robust.filter}
\end{align}
\end{proposition}

\medskip

The rest of the paper is devoted to the proof of Theorem \ref{thm:main}. Section \ref{se:properties} shows that the time evolution of the filter can be described by a (random) probability kernel whose density is bounded away from zero (Lemma \ref{le:v.reg}) and uses it to show that the aforementioned time evolution is a contraction w.r.t. the Hilbert (pseudo-)metric (Proposition \ref{le:H.contract}). Subsection \ref{subse:exp.stab} uses this contraction property to prove the exponential stability of the filter -- i.e., its continuity w.r.t. the initial condition (Corollaries \ref{cor:exp.stab}, \ref{cor:TV.bound.exp} and Lemma \ref{cor:H.contracts.on.pi}). Although such stability is essentially known (see \cite{StabContSpaceTimeOcone}, \cite{StabContSpaceTimeDaPrato}), we chose to prove it herein because we need a slightly stronger version of this result and because its proof follows from the auxiliary results that are important for other parts of the paper.
In Subsection \ref{subse:conv.ave.sq.error}, we establish the Markov property of the filter (Lemma \ref{le:Pi.isMarkov}) and use it, along with the exponential stability, to show that the normalized logarithm of the likelihood ratio, computed for any two parameter values $\theta,\theta'$, converges exponentially fast to a constant, denoted $\Lambda(\theta,\theta')$ and also known as the \emph{contrast function}, as the sample size goes to infinity -- first in the sense of its expectation only (Proposition \ref{prop:exp.conv.sq.error}), and then in the sense of the convergence in probability (Proposition \ref{prop:strong.ergodicity}). Proposition \ref{prop:Lambda.pos} uses the martingale theory to show that the aforementioned exponential convergence and the ergodicity of the signal imply $\Lambda(\theta,\theta')\neq0$ for any $\theta'\notin\overline{\Theta}(\theta)$, without having any tractable representation of $\Lambda(\theta,\theta')$ and without the stationarity or ergodicity of the observation process. Subsection \ref{subse:robust} proves Proposition \ref{le:uniform.robust.filter} and combines this result with the observation that $\Lambda(\theta,\theta')\neq0$ for $\theta'\notin\overline{\Theta}(\theta)$ to show that $\Lambda(\theta,\theta')$ is bounded away from zero uniformly over all $\theta'$ that are bounded away from $\overline{\Theta}(\theta)$ (Corollary \ref{cor:Lambda.bdd.away.from.zero}).
The latter corollary is then used to prove Proposition \ref{prop:last.step}, and the subsequent discussion shows how to deduce Theorem \ref{thm:main} from this proposition.

\section{Preliminary results}\label{se:properties}

We begin with the auxiliary notation.
First, we denote the space of $1$-periodic functions on $\RR^q$ by $\mathrm{P}_1$.
For a function $f$ of $(t,x)$, the notation $f\in\mathrm{P}_1$ means that $f(t,\cdot)\in\mathrm{P}_1$ for all $t$ in the domain of the definition of $f$.
For any measure $\mu$ and any admissible function $f$, we denote $\mu[f]:=\int f\,d\mu$.
We let Assumptions  \ref{periodic}, \ref{holder}, \ref{ass:unifcont} hold throughout the rest of the paper, even if they are not cited explicitly.

\smallskip

We also make use of the following norms. For any $f\in C^{0,k}([T_1,T_2]\times \RR^q\rightarrow \RR)$, any $[T_0,T]\subset[T_1,T_2]$ and any $k\in\NN\cup\{0\}$, we define
\begin{align*}
\|f\|_{C^{k}_{T_0,T}}:=\sup_{t\in [T_0,T]}\|f(t,\cdot)\|_{C^k(\RR^q)}
\end{align*}
and write $C_{T_0,T}$ in place of $C^0_{T_0,T}$.

Given any path $y\in C([T_1,T_2]\rightarrow\RR^m)$, we define its oscillation on an interval $[T_0,T]\subset [T_1,T_2]$ as
\begin{align*}
\mathrm{Osc}_{T_0,T}(y):=\sup_{T_0\leq s\leq t\leq T}|{y_s-y_t}|.
\end{align*}
\smallskip
For any non-negative functions $f,g:\,\RR^q\rightarrow\RR_+$, we define the Hilbert (pseudo-)metric:
\begin{align}\label{eq:defH}
    H(f,g):=\left(1-\frac{\inf_{x:g(x)>0}\frac{f(x)}{g(x)}}{\sup_{x:g(x)>0}\frac{f(x)}{g(x)}}\right)/\left(1+\frac{\inf_{x:g(x)>0}\frac{f(x)}{g(x)}}{\sup_{x:g(x)>0}\frac{f(x)}{g(x)}}\right).
\end{align}
We refer to \cite{BirkhoffContraction} for the general properties of this (pseudo-)metric. We also note that the topology induced by this metric is equivalent to the one induced by the more standard Hilbert metric (see \cite{le2004stability} for the use of the latter metric in stochastic filtering).
With a slight abuse of notation, we apply the above Hilbert metric to measures on subsets of Euclidian spaces whenever these measures have densities with respect to the Lebesgue measure (if one of them does not, the value of the metric is set to $\infty$).







We denote the infinitesimal generator associated with \eqref{eq.Background.dyn.X} by 
\begin{align}
    \cL^\theta \phi:=\frac{1}{2}\sum_{i,j=1}^q  a_{i,j}^\theta\,\pa_{x_i x_j}\phi + \sum_{i=1}^q b_i^\theta\,\pa_{x_i}\phi.
\end{align}


For convenience, we define, for $t\geq0$ and $s\geq 0$,
\begin{align}
& L^{\theta,\nu}_{t,s}:=L^{\theta,\nu}_{t+s}/L^{\theta,\nu}_{t}=\exp\left(-\frac{1}{2}\int_t^{t+s} |\pi^{\theta,\nu}_r[h^\theta]|^2\,dr + \int_t^{t+s} (\pi^{\theta,\nu}_r[h^\theta])^\top dY_r\right),
\label{eq.L.def.2}\\
&E^\theta_{t,s}(x):=\exp\left(\frac{1}{2} |h^\theta(x)|^2\,s - h^\theta(x)^\top\,(Y_{t+s}-Y_t)\right),\label{eq.E.def}
\end{align}
and notice that 
\begin{align*}
dE^\theta_{t,s}(x)= E^\theta_{t,s}(x)\,\left(|h^\theta(x)|^2\,ds-h^\theta(x)^\top\,dY_{t+s}\right),\quad \mathbb{W}_Y\text{-a.s.},
\end{align*} 
with the differentials being with respect to $s$.

\medskip

Let us now proceed to the first result of this section, which shows that Assumptions \ref{periodic} and \ref{holder} imply certain regularity of the solutions to the Cauchy problem associated with $\cL^\theta$. The proof of this lemma is given in Subsection \ref{subse:ass:mixing}.

\begin{lemma}\label{ass:mixing}
For any $T>0$, $f,g\in C^{0,1}([0,T]\times\RR^q)\cap\mathrm{P}_1$, $\theta\in\Theta$ and $0\leq\phi\in C(\RR^q)\cap\mathrm{P}_1$, there exists a unique bounded classical solution $u\in C^{1,2}((0,T)\times\RR^q)\cap C([0,T]\times\RR^q)\cap \mathrm{P}_1$ to the partial differential equation (PDE)
\begin{align}
\partial_t u + \cL^\theta\,u + g^\top\,\nabla u + f\,u=0,\quad u(T,\cdot)=\phi.\label{eq.u.PDE}
\end{align}
In addition, if $U^{\theta,T}_{f,g}[\phi]$ denotes the function $u$ defined above, for any given $\theta\in\Theta$, $\phi\in C(\RR^q)\cap\mathrm{P}_1$, $T>0$, and $f$, $g$ as above, then, the mapping $(\phi,f,g)\mapsto U^{\theta,T}_{f,g}[\phi]$ is Borel measurable, with the associated spaces being equipped with the topology of uniform convergence on compacts.

Moreover, for any $T>0$, there exist non-increasing (in each variable) strictly positive measurable functions $\epsilon_0,\epsilon_1:\,(0,1)\times\RR^2_+\rightarrow(0,\infty)$ and an increasing (in each variable) function $\bar C:\,\RR^3_+\rightarrow(0,\infty)$, such that, for any $z_1,z_2\in\RR_+$,
\begin{align*}
&\int_0^{1}\frac{1}{\epsilon_0(t,z_1,z_2)}\,dt<\infty
\end{align*}
and, for any $\theta\in\Theta$, $\phi\in C^1(\RR^q)\cap\mathrm{P}_1$, $T_0\in(0,T)$, $i=1,\ldots,q$, and any $f$, $g$ as above, we have:
\begin{align}
&{\|U^{\theta,T}_{f,g}[\phi]\|_{C_{0,T-T_0}} +} \|U^{\theta,T}_{f,g}[\partial_{x_i}\phi]\|_{C_{0,T-T_0}}\nonumber \\
&\phantom{????????????????????}+ \|\partial_{x_i}U^{\theta,T}_{f,g}[\phi]\|_{C_{0,T-T_0}}
\leq \frac{\|\phi\|_{L^1([0,1]^q)}}{\epsilon_1\left(T_0,\|f\|_{C^{1}_{0,T}},\|g\|_{C^{1}_{0,T}}\right)},\label{eq.u.positive}\\
&\|\partial_{x_i}U^{\theta,T}_{f,g}[\phi](T-T_0,\cdot)\|_{L^1([0,1)^q)}
\leq \frac{\|\phi\|_{L^1([0,1]^q)}}{\epsilon_0\left(T_0,\|f\|_{C^{1}_{0,T}},\|g\|_{C^{1}_{0,T}}\right)},\label{eq:derder}\\
&\|U^{\theta,T}_{f,g}[\partial_{x_i}\phi]\|_{C_{0,T-T_0}} + \|\partial_{x_i}U^{\theta,T}_{f,g}[\phi]\|_{C_{0,T-T_0}}
\leq \frac{\|\phi\|_{C}}{\epsilon_0\left(T_0,\|f\|_{C^{1}_{0,T}},\|g\|_{C^{1}_{0,T}}\right)},\label{eq:dersup}\\
&\|U^{\theta,T}_{f,g}[\phi]\|_{C_{0,T}} \leq \bar C\left(T,\|f\|_{C^{1}_{0,T}},\|g\|_{C^{1}_{0,T}}\right)\,\|\phi\|_{C},\label{eq:supsup}\\
&\phi\geq 0\,\,\,\Rightarrow\,\,\,\epsilon_1\left(T_0,\|f\|_{C^{1}_{0,T}},\|g\|_{C^{1}_{0,T}}\right)\,\|\phi\|_{L^1([0,1]^q)}
\leq \inf_{t\in[0,T-T_0],x\in\RR^q}U^{\theta,T}_{f,g}[\phi](t,x).\label{eq:lowerassum}
\end{align}
\end{lemma}

\smallskip

The following lemma shows that the filter has a Lebesgue density and provides a convenient PDE representation for the latter. 

\begin{lemma}\label{le:rep.via.U}
For any $(\theta,\nu)\in \Theta\times \mathcal{P}(\RR^q)$, there exists a progressively-measurable $L^2(\RR^q)$-valued process $(p^{\theta,\nu}_{t})_{t>0}$,
such that $\pi^{\theta,\nu}_{t}(dz)=p^{\theta,\nu}_{t}(z)\,dz$, $\PP^{\theta,\nu}$-a.s., for any $t>0$.
In addition, for any $t\geq0$ and $T>0$, the following representation holds for all $\psi \in  C(\RR^q)\cap \mathrm{P}_1$: 
\begin{align}
\int_{\RR^q} E^\theta_{t,T}(x)\,p^{\theta,\nu}_{t+T}(x)\,\psi(x)\,dx=\frac{\int_{\RR^q} U^{\theta,T}_{f^\theta_t,g^\theta_t}[\psi](0,x)\,\pi^{\theta,\nu}_{t}(dx)}{\int_{\RR^q} U^{\theta,T}_{f^\theta_t,g^\theta_t}[1/E^\theta_{t,T}](0,x)\,\pi^{\theta,\nu}_{t}(dx)},\label{eq:rep1}
\end{align}
where $U^{\theta,T}_{f^\theta_t,g^\theta_t}$ is defined in Lemma \ref{ass:mixing} and
\begin{align*}
g^\theta_t(s,x)&:=a^\theta(x)\,\nabla \log E^\theta_{t,s}(x), \\
f^\theta_t(s,x)&:=\cL^\theta \log E^\theta_{t,s}(x)+\frac{1}{2} \nabla^\top \log E^\theta_{t,s}(x) \,a^\theta(x) \,\nabla\log E^\theta_{t,s}(x) .
\end{align*}
\end{lemma}
\begin{proof}
Let us fix arbitrary $\nu\in\mathcal{P}(\RR^q)$ and $\theta\in\Theta$.
Applying Lemma 3.18 and Theorem 3.24 in \cite{BainCrisan}, 
we conclude that, for any $s\geq0$, there exists a progressively measurable measure-valued process $(\rho^{\theta,\nu}_s)$, such that
\begin{align*}
&\pi^{\theta,\nu}_s(dx) = \rho^{\theta,\nu}_s(dx)/L^{\theta,\nu}_s,
\end{align*}
where $L^{\theta,\nu}_s$ is given by \eqref{eq.L.def} and $\rho^{\theta,\nu}$ satisfies the following Zakai equation, for any bounded $\phi\in C(\RR^q)$:
\begin{align*}
    d\rho^{\theta,\nu}_s[\phi]&=\rho^{\theta,\nu}_s[\cL^\theta\phi]\,ds + \rho^{\theta,\nu}_s[h\,\phi]^\top\,dY_s,
    \quad s\geq0,\quad \rho^{\theta,\nu}_0=\nu,\quad \PP^{\theta,\nu}\text{-a.s.}.
\end{align*}

\smallskip

Next, we fix $t\geq0$ and define $v_s(dx):=E^\theta_{t,s}(x)\,\rho^{\theta,\nu}_{t+s}(dx)$, where we recall \eqref{eq.E.def}.
Then, for any bounded $\phi\in C(\RR^q)$,
\begin{align*}
& dv_s[\phi]=\rho^{\theta,\nu}_{t+s}[\cL^\theta(\phi E^\theta_{t,s})]\,ds = v_s[(E^\theta_{t,s})^{-1}\cL^\theta(\phi E^\theta_{t,s})]\,ds\\
&=v_s\left[(E^\theta_{t,s})^{-1} \left(\frac{1}{2}\sum_{i,j=1}^q a_{ij}^\theta \pa_{x_ix_j}(\phi E^\theta_{t,s}) + \sum_{i=1}^q b_i^\theta\pa_{x_i}(\phi E^\theta_{t,s}) \right) \right]\,ds
=v_s[\cL^\theta \phi]\,ds \\
&+ v_s\left[\frac{1}{2}\sum_{i,j=1}^q \left(2a_{ij}^\theta\,\pa_{x_i}\phi\,\pa_{x_j}\log E^\theta_{t,s} + a_{ij}^\theta\,\phi\,\pa_{x_ix_j}E^\theta_{t,s}/E^\theta_{t,s}\right)  + \sum_{i=1}^q b_i^\theta\,\phi\,\pa_{x_i}\log E^\theta_{t,s}\right]\,ds\\
&= v_s[\cL^\theta \phi]\,ds + v_s[(g^\theta_t)^\top\nabla \phi]\,ds + v_s[f^\theta_t\phi]\,ds.
\end{align*}

\smallskip

Next, we choose arbitrary $T>0$ and $0\leq\psi\in C(\RR^q)\cap\mathrm{P}_1$, and consider the unique bounded 1-periodic classical solution $u=U^{\theta,T}_{f^\theta_t,g^\theta_t}[\psi]$ to
\begin{align*}
\pa_s u + \cL^\theta u + (g^\theta_t)^\top\nabla u + f^\theta_t\,u = 0,\quad u(T,\cdot)=\psi,
\end{align*}
which is well defined thanks to Lemma \ref{ass:mixing}.
It is easy to see that $dv_s[u(s,\cdot)]=0$. Thus,
\begin{align}
&L^{\theta,\nu}_{t}\pi^{\theta,\nu}_t\left[U^{\theta,T}_{f^\theta_t,g^\theta_t}[\psi](0,\cdot)\right]=L^{\theta,\nu}_{t}\pi^{\theta,\nu}_t[u(0,\cdot)]=\rho^{\theta,\nu}_t[u(0,\cdot)]=v_0[u(0,\cdot)]\nonumber\\
&=v_T[u(T,\cdot)]=\rho^{\theta,\nu}_{t+T}[E^\theta_{t,T}\,\psi]
= L^{\theta,\nu}_{t+T}\,\pi^{\theta,\nu}_{t+T}[E^\theta_{t,T}\,\psi],\label{eq.Lem3.pf.eq1}
\end{align}
from which we deduce
\begin{align}
&L^{\theta,\nu}_{t,T}\,\pi^{\theta,\nu}_{t+T}[\psi]
= \pi^{\theta,\nu}_t\left[U^{\theta,T}_{f^\theta_t,g^\theta_t}[\psi/E^\theta_{t,T}](0,\cdot)\right].\label{eq.Lem3.pf.eq2}
\end{align}
The above display and the property \eqref{eq.u.positive} imply the existence of a density $p^{\theta,\nu}_t$ of $\pi^{\theta,\nu}_t$, which inherits the progressive-measurability property of $\rho^{\theta,\nu}$.

\smallskip

Applying \eqref{eq.Lem3.pf.eq2} with $\psi\equiv1$, we obtain 
\begin{align}\label{eq:repL}
L^{\theta,\nu}_{t,T}=\int_{\RR^q} U^{\theta,T}_{f^\theta_t,g^\theta_t}[1/E^\theta_{t,t+T}](0,x)\,p^{\theta,\nu}_{t}(x)\,dx.
\end{align}
Combining the above with \eqref{eq.Lem3.pf.eq1}, we obtain \eqref{eq:rep1}.
\end{proof}

\medskip

Next, for any measure $\mu\in \mathcal{P}(\RR^q)$, we define its projection $\hat \mu \in \mathcal{P}([0,1)^q)$ via 
\begin{align}
& \hat\mu(A) := \sum_{e\in\mathbb{Z}^q}\mu(A+e),\label{eq.meas.proj.def}
\end{align}
for any Borel set $A\subset [0,1)^q$. In particular, we introduce $\hat\pi^{\theta,\nu}_t$ as the projection of the filter $\pi^{\theta,\nu}_t$ on $[0,1)^q$ so that, for any $\nu\in\mathcal{P}(\RR^q)$ and $x\in [0,1)^q$, 
\begin{align*}
\hat\pi^{\theta,\nu}_t(dx) =\hat{p}^{\theta,\nu}_t(x)\,dx,
\quad \hat{p}^{\theta,\nu}_t(x):=\sum_{e\in\mathbb{Z}^q} p^{\theta,\nu}_t(x+e).
\end{align*}
Additionally, any $\nu\in\mathcal{P}([0,1)^q)$ is trivially extended to a measure in $\mathcal{P}(\RR^q)$ via setting $\nu(\RR^q{\setminus}[0,1)^q)=0$, so that $\hat\pi^{\theta,\nu}_t$ is well defined. 
It is clear that, for any bounded measurable $\phi\in \mathrm{P}_1$ and $\mu\in \mathcal{P}(\RR^q)$,
\begin{align}\label{eq:truncation}
\hat\pi^{\theta,\hat\mu}_t[\phi] = \pi^{\theta,\mu}_t[\phi].
\end{align}

\medskip

Next, we obtain a kernel representation for the time-propagation of $\hat p^{\theta,\nu}$.
 
\begin{lemma}\label{le:v.reg}
For any $(t,T,\theta)\in \RR_+\times (0,\infty)\times \Theta$, there exists a Borel measurable mapping $C([0,T]\rightarrow\RR^m)\times[0,1)^q\times [0,1)^q\ni((Y_{t+\cdot}-Y_t)_{[0,T]},x,z)\mapsto \hat{K}^{\theta}_{t,T}(x,z)$, such that, for any $t\geq0$, $T>0$, $\theta\in\Theta$ and $\nu\in\mathcal{P}([0,1)^q)$, the following holds for Lebesgue-a.e. $z\in[0,1)^q$ and $\mathbb{W}_Y$-a.e. $Y$:
\begin{align}\label{eq:rep2}
\hat p^{\theta,\nu}_{t+T}(z)=&\frac{\int_{[0,1)^q} \hat{K}^{\theta}_{t,T}(x,z)\hat \pi^{\theta,\nu}_{t}(dx)}{\iint_{[0,1)^q} \hat{K}^{\theta}_{t,T}(x,y)\hat\pi^{\theta,\nu}_{t}(dx)dy},
\end{align}
where we suppress the dependence of $ \hat{K}^{\theta}_{t,T}$ on $Y$, to ease the notation.

Moreover, there exists a measurable function $\epsilon:(0,\infty)\times\RR_+\rightarrow(0,1)$, which is nonincreasing in each variable and is such that, for any $(t,T,\theta,\theta',x,z)\in[0,\infty)\times (0,\infty)\times \Theta^2\times[0,1)^q\times[0,1)^q$, the following hold for $\mathbb{W}_Y$-a.e. $Y$, $\tilde Y$:
\begin{align}
&\epsilon\left(T,\mathrm{Osc}_{t,t+T}(Y)\right)\leq \hat K^\theta_{t,T}(Y_{[t,t+T]}-Y_t;x,z)\leq \frac{1}{\epsilon\left(T,\mathrm{Osc}_{t,t+T}(Y)\right)},\label{eq:boundK}\\
&|\hat K^\theta_{t,T}(Y_{[t,t+T]}-Y_t;x,z)-\hat K^{\theta'}_{t,T}(\tilde Y_{[t,t+T]}-\tilde Y_t;x,z)|\nonumber\\
&\phantom{?????????}\leq \frac{\|a^\theta-a^{\theta'}\|_{C^1}+\|b^\theta-b^{\theta'}\|_{C}+\|h^{\theta}-h^{\theta'}\|_{C^{2}}+\|Y-\tilde Y\|_{C(t,t+T)}}{\epsilon\left(T,\mathrm{Osc}_{t,t+T}(Y)\vee\mathrm{Osc}_{t,t+T}(\tilde Y)\right)}\label{eq:contdv}.
\end{align}
\end{lemma}

The proof of the above lemma is given in Subsection \ref{subse:le:v.reg}.

\medskip

Thanks to Lemma \ref{le:v.reg}, we have
\begin{align}
& \hat{p}^{\theta,\nu}_{t+T}=A^\theta_{t,T}[\hat{\pi}^{\theta,\nu}_{t}],\label{eq.hatp.via.kernel}
\end{align}
where
\begin{align*}
& \hat{A}^{\theta}_{t,T}[\mu](z):=\int_{[0,1)^q} \hat{K}^{\theta}_{t,T}(x,z)\,\mu(dx),
\quad A^\theta_{t,T}[\mu](z):=\frac{\hat{A}^{\theta}_{t,T}[\mu](z)}{\iint_{[0,1)^q} \hat{K}^{\theta}_{t,T}(x,y)\,\mu(dx)dy}
\end{align*}
are $\sigma((Y_{t+r}-Y_{t})_{r\in[0,T]})$-measurable and satisfy \eqref{eq:boundK}--\eqref{eq:contdv} with the associated $\epsilon$.
The above representation allows us to view $\hat p^{\theta,\nu}$ and $\hat\pi^{\theta,\nu}$ as non-anticipative functions of the paths $Y$, and hence as adapted stochastic processes on $(C(\RR_+\rightarrow\RR^m),\mathbb{F}^Y)$.
More importantly, using the above results, we show that $A^{\theta}_{t,T}$ is a contraction under the Hilbert metric.

\begin{proposition}\label{le:H.contract}
Let $\epsilon$ be as in Lemma \ref{le:v.reg}.
Then, the nonincreasing mapping $\gamma:\,x\in (0,\infty)\mapsto\ln \frac{1+\epsilon^2(1,x)}{1-\epsilon^2(1,x)}\in(0,\infty)$ satisfies the following, for any $\nu,\nu'\in\mathcal{P}([0,1)^q)$, $\theta\in\Theta$ and $t>0$:
\begin{align*}
H(\hat{p}^{\theta,\nu}_{t+1},\hat{p}^{\theta,\nu'}_{t+1}) \leq H(\hat{p}^{\theta,\nu}_{t},\hat{p}^{\theta,\nu'}_{t})\,e^{-\gamma(\mathrm{Osc}_{t,t+1}(Y))}
\quad \mathbb{W}_Y\text{-a.s.}
\end{align*}
\end{proposition}
\begin{proof}
Notice that
\begin{align*}
&H(\hat{p}^{\theta,\nu}_{t+1},\hat{p}^{\theta,\nu'}_{t+1}) 
= H(A^\theta_{t,1}[\hat{\pi}^{\theta,\nu}_{t}],A^\theta_{t,1}[\hat{\pi}^{\theta,\nu'}_{t}])\\
&
= H\left(\frac{\hat{A}^\theta_{t,1}[\hat{\pi}^{\theta,\nu}_{t}]}{\iint_{[0,1)^q} \hat{K}^{\theta}_{t,1}(x,y)\,\hat{\pi}^{\theta,\nu}_{t}(dx)dy},\frac{\hat{A}^\theta_{t,1}[\hat{\pi}^{\theta,\nu'}_{t}]}{\iint_{[0,1)^q} \hat{K}^{\theta}_{t,1}(x,y)\,\hat{\pi}^{\theta,\nu'}_{t}dy}\right)= H(\hat{A}^\theta_{t,1}[\hat{\pi}^{\theta,\nu}_{t}],\hat{A}^\theta_{t,1}[\hat{\pi}^{\theta,\nu'}_{t}]).
\end{align*}
Then, it follows from \cite[Propositions 7 and 11]{BirkhoffContraction} and the estimate \eqref{eq:boundK} that
\begin{align*}
&H(\hat{A}^\theta_{t,1}[\hat{\pi}^{\theta,\nu}_{t}],\hat{A}^\theta_{t,1}[\hat{\pi}^{\theta,\nu'}_{t}])\leq H(\hat{\pi}^{\theta,\nu}_{t},\hat{\pi}^{\theta,\nu'}_{t})\,e^{-\gamma(\mathrm{Osc}_{t,t+1}(Y))}
=H(\hat{p}^{\theta,\nu}_{t},\hat{p}^{\theta,\nu'}_{t})\,e^{-\gamma(\mathrm{Osc}_{t,t+1}(Y))}.
\end{align*}
\end{proof}

\section{Proof of consistency}
\subsection{Exponential stability of the filter}
\label{subse:exp.stab}

In this subsection, we derive several useful results following from the contraction property established in Proposition \ref{le:H.contract}.
We start with the following corollary of Proposition \ref{le:H.contract}.

\begin{cor}\label{cor:exp.stab}
For any $\nu,\nu',\nu''\in\mathcal{P}([0,1)^q)$ and $\theta,\theta'\in\Theta$, we have:
\begin{align}
& H(\hat{\pi}^{\theta',\nu'}_{s+t}, \hat{\pi}^{\theta',\nu''}_{s+t}) \leq H(\hat{\pi}^{\theta',\nu'}_s, \hat{\pi}^{\theta',\nu''}_s) e^{-\sum_{i=0}^{\lfloor t\rfloor-1}\gamma(\mathrm{Osc}_{s+i,s+i+1}(Y))},\quad s\geq0,\,\,t\geq1,\quad \PP^{\theta,\nu}_Y\text{-a.s.}\label{eq.exp.stab.local}
\end{align}
\end{cor}

\smallskip

Next, we introduce the pair of measures
\begin{align*}
& \Pi^{\theta_1,\theta_2,\nu_1,\nu_2}_t(dx):=\left(A^{\theta_1}_{0,t}[\nu_1](x)\,dx, A^{\theta_2}_{0,t}[\nu_2](x)\,dx\right)^\top=(\hat{\pi}^{\theta_1,\nu_1}_t(dx),\hat{\pi}^{\theta_2,\nu_2}_t(dx))^\top
\end{align*}
and deduce the following lemma.

\begin{lemma}\label{cor:H.contracts.on.pi}
For any $\nu,\nu_1,\nu_2,\nu_1',\nu_2'\in\mathcal{P}([0,1)^q)$ and $\theta,\theta_1,\theta_2\in\Theta$, we have, for any $s\geq0$ and $t\geq1$, $\PP^{\theta,\nu}_Y$-a.s.:
\begin{align*}
&\|\Pi^{\theta_1,\theta_2,\nu_1,\nu_2}_{s+t}-\Pi^{\theta_1,\theta_2,\nu_1',\nu_2'}_{s+t}\|_{\mathrm{TV}}\\
&\phantom{???????????????????}\leq 1\wedge\left(\frac{2H(\hat{\pi}^{\theta_1,\nu_1}_s,\hat{\pi}^{\theta_1,\nu_1'}_s)\,\Gamma_{s,t}}{1-H(\hat{\pi}^{\theta_1,\nu_1}_s,\hat{\pi}^{\theta_1,\nu_1'}_s)\,\Gamma_{s,t}} 
+ \frac{2H(\hat{\pi}^{\theta_2,\nu_2}_s,\hat{\pi}^{\theta_2,\nu_2'}_s)\,\Gamma_{s,t}}{1-H(\hat{\pi}^{\theta_2,\nu_2}_s,\hat{\pi}^{\theta_2,\nu_2'}_s)\,\Gamma_{s,t}}\right)\\ 
&\Gamma_{s,t}:=e^{-\sum_{i=0}^{\lfloor t\rfloor-1}\gamma(\mathrm{Osc}_{s+i,s+i+1}(Y))},
\end{align*}
where $\|\cdot\|_{\mathrm{TV}}$ is the total variation norm (when applied to a vector of measures, $\|\cdot\|_{\mathrm{TV}}$ denotes the Euclidean norm of the vector of the total variation norms of the associated components).
\end{lemma}
\begin{proof}
First, we notice that
\begin{align*}
&\frac{H(\hat{\pi}^{\theta',\nu'}_{s+t}, \hat{\pi}^{\theta',\nu''}_{s+t})}{1-H(\hat{\pi}^{\theta',\nu'}_{s+t}, \hat{\pi}^{\theta',\nu''}_{s+t})}
\leq \frac{H(\hat{\pi}^{\theta',\nu'}_{s}, \hat{\pi}^{\theta',\nu''}_{s}) e^{-\sum_{i=0}^{\lfloor t\rfloor-1}\gamma(\mathrm{Osc}_{s+i,s+i+1}(Y))}}{1-H(\hat{\pi}^{\theta',\nu'}_{s}, \hat{\pi}^{\theta',\nu''}_{s})e^{-\sum_{i=0}^{\lfloor t\rfloor-1}\gamma(\mathrm{Osc}_{s+i,s+i+1}(Y))}}.
\end{align*}
Then,
\begin{align*}
&\|\Pi^{\theta_1,\theta_2,\nu_1,\nu_2}_{s+t}-\Pi^{\theta_1,\theta_2,\nu_1',\nu_2'}_{s+t}\|_{\mathrm{TV}}\leq \|\hat{\pi}^{\theta_1,\nu_1}_{s+t}-\hat{\pi}^{\theta_1,\nu_1'}_{s+t}\|_{\mathrm{TV}}+\|\hat{\pi}^{\theta_2,\nu_2}_{t+s}-\hat{\pi}^{\theta_2,\nu_2'}_{t+s}\|_{\mathrm{TV}}\\
&\leq \int \hat{p}^{\theta_1,\nu_1'}_{s+t}(x)\left|\frac{\hat{p}^{\theta_1,\nu_1}_{s+t}(x)}{\hat{p}^{\theta_1,\nu_1'}_{s+t}(x)}-1\right|dx+\int \hat{p}^{\theta_2,\nu_2'}_{t+s}(x)\left|\frac{\hat{p}^{\theta_2,\nu_2}_{t+s}(x)}{\hat{p}^{\theta_2,\nu_2'}_{t+s}(x)}-1\right|dx\\
&\leq \sup\frac{\hat{p}^{\theta_1,\nu_1}_{s+t}}{\hat{p}^{\theta_1,\nu_1'}_{s+t}} - \inf\frac{\hat{p}^{\theta_1,\nu_1}_{s+t}}{\hat{p}^{\theta_1,\nu_1'}_{s+t}}
+\sup\frac{\hat{p}^{\theta_2,\nu_2}_{s+t}}{\hat{p}^{\theta_2,\nu_2'}_{s+t}} - \inf\frac{\hat{p}^{\theta_2,\nu_2}_{s+t}}{\hat{p}^{\theta_2,\nu_2'}_{s+t}}\\
&\leq \frac{2H(\hat{\pi}^{\theta_1,\nu_1}_{s+t},\hat{\pi}^{\theta_1,\nu_1'}_{s+t})}{1-H(\hat{\pi}^{\theta_1,\nu_1}_{s+t},\hat{\pi}^{\theta_1,\nu_1'}_{s+t})}
+ \frac{2H(\hat{\pi}^{\theta_2,\nu_2}_{s+t},\hat{\pi}^{\theta_2,\nu_2'}_{s+t})}{1-H(\hat{\pi}^{\theta_2,\nu_2}_{s+t},\hat{\pi}^{\theta_2,\nu_2'}_{s+t})}\\
&\leq 2 \left(\frac{H(\hat{\pi}^{\theta_1,\nu_1}_{s},\hat{\pi}^{\theta_1,\nu_1'}_{s})}{1-H(\hat{\pi}^{\theta_1,\nu_1}_{s},\hat{\pi}^{\theta_1,\nu_1'}_{s})\,\Gamma_{s,t}}
+ \frac{H(\hat{\pi}^{\theta_2,\nu_2}_{s},\hat{\pi}^{\theta_2,\nu_2'}_{s})}{1-H(\hat{\pi}^{\theta_2,\nu_2}_{s},\hat{\pi}^{\theta_2,\nu_2'}_{s})\,\Gamma_{s,t}}\right) \,\Gamma_{s,t}.
\end{align*}
\end{proof}

\medskip

In what follows, we sometimes need to analyze the expected value of the quantity estimated in the above lemma. This is achieved via the following corollary.

\begin{cor}\label{cor:TV.bound.exp}
There exists a constant $\bar\gamma>0$ such that, for any $\nu,\nu_1,\nu_2,\nu_1',\nu_2'\in\mathcal{P}([0,1)^q)$, $\theta,\theta_1,\theta_2\in\Theta$, $s\geq0$ and $t\geq1$, we have:
\begin{align*}
&\EE^{\theta,\nu}\|\Pi^{\theta_1,\theta_2,\nu_1,\nu_2}_{s+t}-\Pi^{\theta_1,\theta_2,\nu_1',\nu_2'}_{s+t}\|_{\mathrm{TV}}\\
&\phantom{???????????????????}\leq 2 e^{-(\lfloor t\rfloor-1)\bar\gamma}\, \EE^{\theta,\nu}\left(\frac{H(\hat{\pi}^{\theta_1,\nu_1}_s,\hat{\pi}^{\theta_1,\nu_1'}_s)}{1-H(\hat{\pi}^{\theta_1,\nu_1}_s,\hat{\pi}^{\theta_1,\nu_1'}_s)} 
+ \frac{H(\hat{\pi}^{\theta_2,\nu_2}_s,\hat{\pi}^{\theta_2,\nu_2'}_s)}{1-H(\hat{\pi}^{\theta_2,\nu_2}_s,\hat{\pi}^{\theta_2,\nu_2'}_s)}\right).
\end{align*}
\end{cor}
\begin{proof}
Recalling \eqref{eq.Y.projected} and using the boundedness of $h$, we deduce the existence of a constant $C>0$ s.t.
\begin{align*}
&\mathrm{Osc}_{s+i,s+i+1}(Y) \leq C + \mathrm{Osc}_{s+i,s+i+1}(\tilde W^{\theta,\nu})\quad\PP^{\theta,\nu}\text{-a.s.,}
\end{align*}
for any $\theta$, $\nu\in\mathcal{P}([0,1)^q)$, $i=0,1,\ldots$ and  $s>0$. Using the monotonicity of $\gamma$, as well as the independence and stationarity of the increments of Brownian motion, we obtain:
\begin{align*}
& \EE^{\theta,\nu}\left[ e^{-\sum_{i=0}^{\lfloor t\rfloor-1}\gamma(\mathrm{Osc}_{s+i,s+i+1}(Y))}\,\vert\,\mathcal{F}^Y_{s}\right]
\leq \EE^{\theta,\nu}\left[ e^{-\sum_{i=0}^{\lfloor t\rfloor-1}\gamma(C + \mathrm{Osc}_{s+i,s+i+1}(\tilde W^{\theta,\nu}))}\,\vert\,\mathcal{F}^Y_{s}\right]\\
& = \prod_{i=0}^{\lfloor t\rfloor-1}  \EE^{\theta,\nu} e^{-\gamma(C + \mathrm{Osc}_{s+i,s+i+1}(\tilde W^{\theta,\nu}))}
=: e^{-(\lfloor t\rfloor-1)\bar\gamma}.
\end{align*}
Then, the statement of the corollary follows from Lemma \ref{cor:H.contracts.on.pi} and from the inequality $H/(1-H\,\Gamma)\leq H/(1-H)$.
\end{proof}

\subsection{Convergence of the normalized log-likelihood ratio}
\label{subse:conv.ave.sq.error}

The main goal of this subsection is to show that the exponential stability of the filter and the identifiability property imply that the ``normalized finite-variation component of the log-likelihood ratio" (this terminology is justified by \eqref{eq.ratio.terminology}),
\begin{align*}
\frac{1}{t}\int_0^t \left|\hat{\pi}^{\theta,\nu}_s[h^\theta]-\hat{\pi}^{\theta',\nu}_s[h^{\theta'}]\right|^2ds,
\end{align*} 
converges to a strictly positive limit as $t\rightarrow\infty$, for any $\theta'\notin\overline{\Theta}(\theta)$.
First, we deduce the Markov property of the process $\Pi$.

\begin{lemma}\label{le:Pi.isMarkov}
For any $t,s\geq0$, any bounded $f\in C(\RR^q)\times C(\RR^q)$ and $g\in C(\RR^2)$, any $\nu,\nu'\in\mathcal{P}([0,1)^q)$ and any $\theta,\theta'\in\Theta$, we have:
\begin{align*}
& \EE^{\theta,\nu}\left[g\left(\Pi^{\theta,\theta',\nu,\nu'}_{t+s}[f]\right)\,\vert\,\mathcal{F}^Y_s\right] = \left.\EE^{\theta,\nu_1}g\left(\Pi^{\theta,\theta',\nu_1,\nu_2}_{t}[f]\right)\right|_{\nu_1=\hat{\pi}^{\theta,\nu}_s,\,\nu_2=\hat{\pi}^{\theta',\nu'}_s},\quad \PP^{\theta,\nu}\text{-a.s.}
\end{align*}
\end{lemma}
\begin{proof}
The representation \eqref{eq.hatp.via.kernel} implies
\begin{align*}
& \Pi^{\theta,\theta',\nu,\nu'}_{t+s}[f] = \left( \int A^{\theta}_{s,t}[\hat{\pi}^{\theta,\nu}_s](x)\,f^1(x)\,dx, \int A^{\theta'}_{s,t}[\hat{\pi}^{\theta',\nu'}_s](x)\,f^2(x)\,dx \right)^\top,
\end{align*}
where $\hat{\pi}^{\theta,\nu}_s$, $\hat{\pi}^{\theta',\nu'}_s$ are $\mathcal{F}^Y_s$-measurable, and $A^{\theta}_{s,t}$, $A^{\theta'}_{s,t}$ are given by measurable functions of $(Y_{\cdot} - Y^{\theta,\nu}_{s})_{[s,s+t]}$, the same across all $s$.
Using the above and treating $\Pi^{\theta,\theta',\nu,\nu'}$ as a function of the paths of $Y$, we obtain:
\begin{align*}
& \EE^{\theta,\nu}\left[g\left(\Pi^{\theta,\theta',\nu,\nu'}_{t+s}[f]\right)\,\vert\,\mathcal{F}^Y_s\right]\\
&= \EE^{\mathbb{W}_Y}\left[g\left(\Pi^{\theta,\theta',\nu,\nu'}_{t+s}[f]\right)\,\exp\left(-\frac{1}{2}\int_s^{t+s}|\hat{\pi}^{\theta,\nu}_r[h^\theta]|^2\,dr
+ \int_s^{t+s} \hat{\pi}^{\theta,\nu}_r[h^\theta]^\top\,dY_r\right)\,\vert\,\mathcal{F}^Y_s\right]\\
&= \EE^{\mathbb{W}_Y}\left[g\left(\Pi^{\theta,\theta',\nu_1,\nu_2}_{t}[f]\right)\,\exp\left(-\frac{1}{2}\int_0^{t}|\hat{\pi}^{\theta,\nu_1}_r[h^\theta]|^2\,dr\right.\right.\\
&\phantom{?????}\left.\left.+ \int_0^{t} \hat{\pi}^{\theta,\nu_1}_r[h^\theta]^\top\,dY_r\right)\right]_{\nu_1:=\hat{\pi}^{\theta,\nu}_s,\,\nu_2:=\hat{\pi}^{\theta',\nu'}_s}
=\left.\EE^{\theta,\nu_1}g\left(\Pi^{\theta,\theta',\nu_1,\nu_2}_{t}[f]\right)\right|_{\nu_1=\hat{\pi}^{\theta,\nu}_s,\,\nu_2=\hat{\pi}^{\theta',\nu'}_s}.
\end{align*}
\end{proof}

\medskip

Next, we prove the exponential convergence of the expected normalized finite-variation component of the log-likelihood ratio.
The first step in this proof is given by the following lemma.

\begin{lemma}\label{le:exp.stab.improved}
For any bounded $f\in C(\RR^q)\times C(\RR^q)$ and any locally Lipschitz $g:\,\RR^2\rightarrow\RR$, there exists a constant $\tilde C>0$, such that, for any $s,t\geq0$ and $\theta,\theta'\in\Theta$,
\begin{align*}
&\left|\EE^{\theta,\nu^\theta_0} g\left(\Pi^{\theta,\theta',\nu^\theta_0,\nu^{\theta'}_0}_{t}[f]\right) - \EE^{\theta,\nu^\theta_0} g\left(\Pi^{\theta,\theta',\nu^\theta_0,\nu^{\theta'}_0}_{t+s}[f]\right)\right|\\
&\leq \tilde C\,e^{-\bar\gamma\,t}\,
\sup_{r\in[0,1]}\EE^{\theta,\nu^\theta_0}\left[\frac{H(\nu^\theta_0,\hat{\pi}^{\theta,\nu^\theta_0}_{r})}{1-H(\nu^\theta_0,\hat{\pi}^{\theta,\nu^\theta_0}_{r})} + \frac{H(\nu^{\theta'}_0,\hat{\pi}^{\theta',\nu^{\theta'}_0}_{r})}{1-H(\nu^{\theta'}_0,\hat{\pi}^{\theta',\nu^{\theta'}_0}_{r})}\right].
\end{align*}
\end{lemma}
\begin{proof}
First, we observe that, for any $\phi\in L^\infty([0,1)^q)$,
\begin{align*}
& \EE^{\theta,\nu^\theta_0} \hat{\pi}^{\theta,\nu^\theta_0}_t[\phi] = \EE^{\theta,\nu^\theta_0} \phi(X_t) = \EE^{\theta,\nu^\theta_0} \phi(X_0),\mbox{ and}\\
& \EE^{\theta,\nu^\theta_0}\left[ \left.\EE^{\theta,\nu_1}g\left(\Pi^{\theta,\theta',\nu^\theta_0,\nu^{\theta'}_0}_{t}[f]\right)\right|_{\nu_1=\hat{\pi}^{\theta,\nu^\theta_0}_s} \right] = \EE^{\theta,\nu^\theta_0}g\left(\Pi^{\theta,\theta',\nu^\theta_0,\nu^{\theta'}_0}_{t}[f]\right)
\end{align*}
holds due to {the linearity of the operator $\nu_1\mapsto \EE^{\theta,\nu_1}$} and the definition of $\nu^\theta_0$.
Using the above, along with Corollary \ref{cor:TV.bound.exp} and Lemma \ref{le:Pi.isMarkov}, we obtain, for all $s\geq 0$ and $t\geq 1$:
\begin{align}
& \left|\EE^{\theta,\nu^\theta_0} g\left(\Pi^{\theta,\theta',\nu^\theta_0,\nu^{\theta'}_0}_{t}[f]\right) - \EE^{\theta,\nu^\theta_0} g\left(\Pi^{\theta,\theta',\nu^\theta_0,\nu^{\theta'}_0}_{t+s}[f]\right)\right|\notag\\
&=\left|\EE^{\theta,\nu^\theta_0}\left( \left.\EE^{\theta,\nu_1} g\left(\Pi^{\theta,\theta',\nu^\theta_0,\nu^{\theta'}_0}_{t}[f]\right)\right|_{\nu_1=\hat{\pi}^{\theta,\nu^\theta_0}_s} - \EE^{\theta,\nu^\theta_0}\left[g\left(\Pi^{\theta,\theta',\nu,\nu^{\theta'}_0}_{t+s}[f]\right)\,\vert\,\mathcal{F}^Y_s\right] \right)\right|\notag\\
&= \left|
\EE^{\theta,\nu^\theta_0}\left(
\left.\EE^{\theta,\nu_1} g\left(\Pi^{\theta,\theta',\nu^\theta_0,\nu^{\theta'}_0}_{t}[f]\right)\right|_{\nu_1=\hat{\pi}^{\theta,\nu^\theta_0}_s} 
- \left.\EE^{\theta,\nu_1} g\left(\Pi^{\theta,\theta',\nu_1,\nu_2}_{t}[f]\right)\right|_{\nu_1=\hat{\pi}^{\theta,\nu^\theta_0}_s,\nu_2=\hat{\pi}^{\theta',\nu'}_s} \right)
\right|\notag\\ 
&\leq \EE^{\theta,\nu^\theta_0}\left(\EE^{\theta,\nu_1}\left. \left[\left| g\left(\Pi^{\theta,\theta',\nu^\theta_0,\nu^{\theta'}_0}_t[f]\right) - g\left(\Pi^{\theta,\theta',\nu_1,\nu_2}_{t}[f]\right)\right|\right]\right|_{\nu_1=\hat{\pi}^{\theta,\nu^\theta_0}_s,\nu_2=\hat{\pi}^{\theta',\nu^{\theta'}_0}_s}\right)\notag\\
&\leq 2C_{f,g}\,e^{-\bar\gamma\,(\lfloor t\rfloor-1)}\, \EE^{\theta,\nu^\theta_0}\left[\frac{H(\nu^\theta_0,\hat{\pi}^{\theta,\nu^\theta_0}_s)}{1-H(\nu^\theta_0,\hat{\pi}^{\theta,\nu^\theta_0}_s)} + \frac{H(\nu^{\theta'}_0,\hat{\pi}^{\theta',\nu^{\theta'}_0}_s)}{1-H(\nu^{\theta'}_0,\hat{\pi}^{\theta',\nu^{\theta'}_0}_s)}\right].\label{eq:boundts}
\end{align}
Next, we improve the above result by removing the dependence of the right hand side on $s$. The telescopic sum gives
\begin{align*}
&\left|\EE^{\theta,\nu^\theta_0} g\left(\Pi^{\theta,\theta',\nu^\theta_0,\nu^{\theta'}_0}_{t}[f]\right) - \EE^{\theta,\nu^\theta_0} g\left(\Pi^{\theta,\theta',\nu^\theta_0,\nu^{\theta'}_0}_{t+s}[f]\right)\right|\\
&\leq \sum_{i=1}^{\lfloor s\rfloor} \left|\EE^{\theta,\nu^\theta_0} g\left(\Pi^{\theta,\theta',\nu^\theta_0,\nu^{\theta'}_0}_{t+i-1}[f]\right) - \EE^{\theta,\nu^\theta_0} g\left(\Pi^{\theta,\theta',\nu^\theta_0,\nu^{\theta'}_0}_{t+i}[f]\right)\right|\\
&+ \left|\EE^{\theta,\nu^\theta_0} g\left(\Pi^{\theta,\theta',\nu^\theta_0,\nu^{\theta'}_0}_{t+\lfloor s\rfloor}[f]\right) - \EE^{\theta,\nu^\theta_0} g\left(\Pi^{\theta,\theta',\nu^\theta_0,\nu^{\theta'}_0}_{t+\lfloor s\rfloor+(s-\lfloor s\rfloor)}[f]\right)\right|.
\end{align*}
Thus, \eqref{eq:boundts} leads to
\begin{align*}
&\left|\EE^{\theta,\nu^\theta_0} g\left(\Pi^{\theta,\theta',\nu^\theta_0,\nu^{\theta'}_0}_{t}[f]\right) - \EE^{\theta,\nu^\theta_0} g\left(\Pi^{\theta,\theta',\nu^\theta_0,\nu^{\theta'}_0}_{t+s}[f]\right)\right|\\
&\leq \sum_{i=1}^{\lfloor s\rfloor} C_1\,e^{-\bar\gamma\,(\lfloor t\rfloor+i-2)}\, \EE^{\theta,\nu^\theta_0}\left[\frac{H(\nu^\theta_0,\hat{\pi}^{\theta,\nu^\theta_0}_1)}{1-H(\nu^\theta_0,\hat{\pi}^{\theta,\nu^\theta_0}_1)} + \frac{H(\nu^{\theta'}_0,\hat{\pi}^{\theta',\nu^{\theta'}_0}_1)}{1-H(\nu^{\theta'}_0,\hat{\pi}^{\theta',\nu^{\theta'}_0}_1)}\right]\\
&+ C_1\,e^{-\bar\gamma\,(t+\lfloor s\rfloor-1)}\, \EE^{\theta,\nu^\theta_0}\left[\frac{H(\nu^\theta_0,\hat{\pi}^{\theta,\nu^\theta_0}_{s-\lfloor s\rfloor})}{1-H(\nu^\theta_0,\hat{\pi}^{\theta,\nu^\theta_0}_{s-\lfloor s\rfloor})} + \frac{H(\nu^{\theta'}_0,\hat{\pi}^{\theta',\nu^{\theta'}_0}_{s-\lfloor s\rfloor})}{1-H(\nu^{\theta'}_0,\hat{\pi}^{\theta',\nu^{\theta'}_0}_{s-\lfloor s\rfloor})}\right]\\
&\leq \tilde C\,e^{-\bar\gamma\,t}\,
\sup_{r\in[0,1]}\EE^{\theta,\nu^\theta_0}\left[\frac{H(\nu^\theta_0,\hat{\pi}^{\theta,\nu^\theta_0}_{r})}{1-H(\nu^\theta_0,\hat{\pi}^{\theta,\nu^\theta_0}_{r})} + \frac{H(\nu^{\theta'}_0,\hat{\pi}^{\theta',\nu^{\theta'}_0}_{r})}{1-H(\nu^{\theta'}_0,\hat{\pi}^{\theta',\nu^{\theta'}_0}_{r})}\right].
\end{align*}
\end{proof}

\smallskip

The next lemma shows that the multiplicative factor in the estimate provided by Lemma \ref{le:exp.stab.improved} is finite.
Its proof is given in Subsection \ref{subse:expH.est}.

\begin{lemma}\label{le:bound.supH}
For any $\theta,\theta'\in\Theta$, we have
\begin{align}\label{eq:unifbound}
\sup_{r\in[0,1]}\EE^{\theta,\nu^\theta_0}\frac{H(\nu^{\theta'}_0,\hat{\pi}^{\theta',\nu^{\theta'}_0}_{r})}{1-H(\nu^{\theta'}_0,\hat{\pi}^{\theta',\nu^{\theta'}_0}_{r})}  <\infty.
\end{align}
\end{lemma}

\medskip

An application of Lemmas \ref{le:exp.stab.improved} and \ref{le:bound.supH} to 
\begin{align*}
g(\mu_1[f^1],\mu_2[f^2])= (\mu_1[h^\theta_i]-\mu_2[h^{\theta'}_i])^2,\quad i=1,\ldots,q,
\end{align*}
yields the desired exponential convergence of the expected normalized finite-variation component of the log-likelihood ratio, as $t\rightarrow\infty$. This is summarized in the following proposition.

\begin{proposition}\label{prop:exp.conv.sq.error}
There exists a function $\Lambda:\,\Theta\times\Theta\rightarrow\RR$ and a constant $\hat C>0$, such that, for any $\theta,\theta'\in\Theta$,
\begin{align*}
& \left|\EE^{\theta,\nu^\theta_0}\left|\hat{\pi}^{\theta,\nu^\theta_0}_t[h^{\theta}] - \hat{\pi}^{\theta',\nu^{\theta'}_0}_t[h^{\theta'}]\right|^2 - \Lambda(\theta,\theta')\right|
\leq \hat C\,e^{-\overline{\gamma}\,t},\quad t\geq0.
\end{align*}
\end{proposition}

\medskip

Next, we show that the exponential convergence stated in Proposition \ref{prop:exp.conv.sq.error} implies the strict positivity of $\Lambda(\theta,\theta')$ for $\theta'\notin\overline{\Theta}(\theta)$.

\begin{proposition}\label{prop:Lambda.pos}
$\Lambda(\theta,\theta')>0$ for all $\theta\in\Theta$ and $\theta'\notin\overline{\Theta}(\theta)$.
\end{proposition}
\begin{proof}
We fix $\theta\in\Theta$, $\theta'\notin\overline{\Theta}(\theta)$ and use Lemma \ref{le:singularity} to deduce that the measures $\PP^{\theta,\nu^\theta_0}_Y$ and $\PP^{\theta',\nu^{\theta'}_0}_Y$ are mutually singular. Then, we recall Proposition \ref{prop:exp.conv.sq.error} to obtain
\begin{align*}
& \lim_{t\rightarrow\infty}\EE^{\theta,\nu^\theta_0}\left|\hat\pi^{\theta,\nu^\theta_0}_t[h^{\theta}] - \hat\pi^{\theta',\nu^{\theta'}_0}_t[h^{\theta'}]\right|^2 = \Lambda(\theta,\theta'),
\end{align*}
with the convergence rate being exponential.
We argue by contradiction and assume that $\Lambda(\theta,\theta')=0$.
The latter assumption and the exponential convergence rate in Proposition \ref{prop:exp.conv.sq.error} imply
\begin{align*}
\EE^{\theta,\nu^\theta_0}\int_0^\infty \left|\hat\pi^{\theta,\nu^\theta_0}_t[h^{\theta}] - \hat\pi^{\theta',\nu^{\theta'}_0}_t[h^{\theta'}]\right|^2\, dt < \infty,
\end{align*}
which in turn yields
\begin{align}
&\int_0^\infty \left|\hat\pi^{\theta,\nu^\theta_0}_t[h^{\theta}] - \hat\pi^{\theta',\nu^{\theta'}_0}_t[h^{\theta'}]\right|^2\, dt < \infty\quad \PP^{\theta,\nu^\theta_0}\text{-a.s.}
\label{eq.Lambda.pos.eq1}
\end{align}
Next, we recall \eqref{eq.L.def} and consider the likelihood ratio:
\begin{align}
& R_t:=\frac{L^{\theta',\nu^{\theta'}_0}_t}{L^{\theta,\nu^\theta_0}_t}=\exp\left(-\frac{1}{2}\int_0^t |\hat\pi^{\theta',\nu^{\theta'}_0}_s[h^{\theta'}]-\hat\pi^{\theta,\nu^\theta_0}_s[h^{\theta}]|^2\,ds \right.\label{eq.ratio.terminology}\\
&\phantom{???????????????????????????????}\left.+ \int_0^t (\hat\pi^{\theta',\nu^{\theta'}_0}_s[h^{\theta'}]-\hat\pi^{\theta,\nu^\theta_0}_s[h^{\theta}])^\top d\tilde W^{\theta,\nu^\theta_0}_s\right).\nonumber
\end{align}
Equation \eqref{eq.Lambda.pos.eq1} implies that the quadratic variation of the continuous martingale
\begin{align*}
&\int_0^t (\hat\pi^{\theta',\nu^{\theta'}_0}_s[h^{\theta'}]-\hat\pi^{\theta,\nu^\theta_0}_s[h^{\theta}])^\top d\tilde W^{\theta,\nu^\theta_0}_s
\end{align*}
converges $\PP^{\theta,\nu^\theta_0}$-a.s. to a finite limit, as $t\rightarrow\infty$. Then, the above martingale enjoys the same property, and, hence, the likelihood ratio $R_t$ converges $\PP^{\theta,\nu^\theta_0}$-a.s. to a strictly positive limit, denoted $\xi$, as $t\rightarrow\infty$.
We note that $\xi$ is measurable with respect to $\mathcal{F}^Y_\infty=\sigma(\bigcup_{t\geq0}\mathcal{F}^Y_t)$.

Next, we consider any $\epsilon>0$ such that $\PP^{\theta,\nu^\theta_0}_Y(\xi>\epsilon)=\PP^{\theta,\nu^\theta_0}(\xi>\epsilon)>0$
and notice that $R_t=d\PP^{\theta',\nu^{\theta'}_0}_{Y,t}/d\PP^{\theta,\nu^\theta_0}_{Y,t}$. Then, we consider an arbitrary event $A\in\bigcup_{t>0}\mathcal{F}^Y_t$ and apply Fatou's lemma to obtain:
\begin{align*}
& \epsilon\,\PP^{\theta,\nu^\theta_0}_Y(A\cap\{\xi>\epsilon\})\leq\EE^{\theta,\nu^\theta_0}\left[ \xi\,\bone_{A\cap\{\xi>\epsilon\}} \right] \leq \EE^{\theta,\nu^\theta_0}\left[ \lim_{t\rightarrow\infty}R_t\,\bone_{A} \right]\\
&\phantom{????????????????????????????????????????????}\leq \lim_{t\rightarrow\infty}\EE^{\theta,\nu^\theta_0}\left[ R_t\,\bone_{A} \right]
= \PP^{\theta',\nu^{\theta'}_0}_Y(A).
\end{align*}
Applying the monotone class theorem, we conclude that the inequality $\epsilon\,\PP^{\theta,\nu^\theta_0}_Y(A\cap\{\xi>\epsilon\})\leq \PP^{\theta',\nu^{\theta'}_0}_Y(A)$ holds for all $A\in\mathcal{F}^Y_\infty$.
Recall that the event $\{\xi>\epsilon\}\in\mathcal{F}^Y_\infty$ has a strictly positive probability under $\PP^{\theta,\nu^\theta_0}_Y$, and we have shown that the restriction of $\PP^{\theta,\nu^\theta_0}_Y$ to $\{\xi>\epsilon\}$ is absolutely continuous w.r.t. the restriction of $\PP^{\theta',\nu^{\theta'}_0}_Y$ to $\{\xi>\epsilon\}$.
This contradicts the mutual singularity of $\PP^{\theta,\nu^\theta_0}_Y$ and $\PP^{\theta',\nu^{\theta'}_0}_Y$ and completes the proof.
\end{proof}

\medskip

Next, we use the above results and the uniform stability of the signal (Lemma \ref{le:expStab.X}) to show that the variance of the normalized finite-variation component of the log-likelihood ratio vanishes for large times.

\begin{lemma}\label{le:vanishing.variance}
For any bounded $f\in C(\RR^q)\times C(\RR^q)$, any locally Lipschitz $g:\,\RR^2\rightarrow\RR$ and any $\theta,\theta'\in\Theta$, we have:
\begin{align*}
\lim_{T\rightarrow\infty}\text{Var}^{\theta,\nu^\theta_0}\left(\frac{1}{T}\int_0^T g(\Pi^{\theta,\theta',\nu^\theta_0,\nu^{\theta'}_0}_{t}[f]) dt\right) = 0.
\end{align*}
\end{lemma}
\begin{proof}
Denoting $e_t:=\EE^{\theta,\nu^\theta_0}g(\Pi^{\theta,\theta',\nu^\theta_0,\nu^{\theta'}_0}_{t}[f])$, we obtain
\begin{align*}
&\text{Var}^{\theta,\nu^{\theta}_0}\left(\int_0^T g(\Pi^{\theta,\theta',\nu^{\theta}_0,\nu^{\theta'}_0}_{t}[f]) dt\right)\\
&= \EE^{\theta,\nu^{\theta}_0}\int_0^T\int_0^T (g(\Pi^{\theta,\theta',\nu^{\theta}_0,\nu^{\theta'}_0}_{t}[f])-e_t) (g(\Pi^{\theta,\theta',\nu^{\theta}_0,\nu^{\theta'}_0}_{s}[f])-e_s) \,dt\,ds\\
&= 2\EE^{\theta,\nu^{\theta}_0}\int_{0\leq s\leq t\leq T} \left(\EE^{\theta,\nu^{\theta}_0}\left[g(\Pi^{\theta,\theta',\nu^{\theta}_0,\nu^{\theta'}_0}_{t}[f])\,\vert\,\mathcal{F}^Y_s\right]
- e_{t-s} + (e_{t-s}-e_t)\right)\\
&\phantom{???????????????????????????????????????????????}\times(g(\Pi^{\theta,\theta',\nu^{\theta}_0,\nu^{\theta'}_0}_{s}[f])-e_s) \,dt\,ds.
\end{align*}
Since
\begin{align*}
& |e_{t-s}-e_t|\leq C_1\,e^{-\bar\gamma\,(t-s)},\quad |g(\Pi^{\theta,\theta',\nu^{\theta}_0,\nu^{\theta'}_0}_{s}[f])-e_s|\leq C_2,
\end{align*}
it suffices to estimate the expected value of
\begin{align}
&\left|\EE^{\theta,\nu^{\theta}_0}\left[g(\Pi^{\theta,\theta',\nu^{\theta}_0,\nu^{\theta'}_0}_{t}[f])\,\vert\,\mathcal{F}^Y_s\right]
- e_{t-s}\right|\nonumber\\
&{\leq} \left|\EE^{\theta,\nu_1} g(\Pi^{\theta,\theta',\nu_1,\nu_2}_{t-s}[f])
- \EE^{\theta,\nu_1} g(\Pi^{\theta,\theta',\nu^{\theta}_0,\nu^{\theta'}_0}_{t-s}[f])\right|_{\nu_1:=\hat\pi^{\theta,\nu^{\theta}_0}_s,\,\nu_2:=\hat\pi^{\theta',\nu^{\theta'}_0}_s}\label{eq.vanishVar.pf.eq2}\\
& + \left|\EE^{\theta,\nu_1} g(\Pi^{\theta,\theta',\nu^{\theta}_0,\nu^{\theta'}_0}_{t-s}[f])
- \EE^{\theta,\nu^{\theta}_0} g(\Pi^{\theta,\theta',\nu^{\theta}_0,\nu^{\theta'}_0}_{t-s}[f])\right|_{\nu_1:=\hat\pi^{\theta,\nu^{\theta}_0}_s}.\nonumber
\end{align}
Let us estimate the first summand in the right hand side of the above, using Lemma \ref{cor:H.contracts.on.pi}:
\begin{align*}
&\left|\EE^{\theta,\nu_1} g(\Pi^{\theta,\theta',\nu_1,\nu_2}_{t-s}[f])
- \EE^{\theta,\nu_1} g(\Pi^{\theta,\theta',\nu^{\theta}_0,\nu^{\theta'}_0}_{t-s}[f])\right|_{\nu_1:=\hat{\pi}^{\theta,\nu^{\theta}_0}_s,\,\nu_2:=\hat{\pi}^{\theta',\nu^{\theta'}_0}_s}\\
&\leq C_3 \left.\EE^{\theta,\nu_1}\left[ 1\wedge\left(\frac{H(\nu_1,\nu^{\theta}_0)\,\Gamma_{0,t-s}}{1-H(\nu_1,\nu^{\theta}_0)\,\Gamma_{0,t-s}} 
+ \frac{H(\nu_2,\nu^{\theta'}_0)\,\Gamma_{0,t-s}}{1-H(\nu_2,\nu^{\theta'}_0)\,\Gamma_{0,t-s}}\right) \right]\right|_{\nu_1:=\hat{\pi}^{\theta,\nu^{\theta}_0}_s,\,\nu_2:=\hat{\pi}^{\theta',\nu^{\theta'}_0}_s}.
\end{align*}
Recalling that $\mathrm{Osc}_{s+i,s+i+1}(Y)\leq C_4 + \mathrm{Osc}_{s+i,s+i+1}(\tilde{W}^{\theta,\nu_1})$, and that the latter are i.i.d. random variables across all $i=0,1,\ldots$, under $\PP^{\theta,\nu_1}$, we use the inequality $H\Gamma(1-H\Gamma)\leq(1/\Gamma-1)^{-1}$ to conclude that, for all $\nu_1,\nu_2\in\mathcal{P}([0,1)^q)$,
\begin{align*}
&\left|\EE^{\theta,\nu_1} g(\Pi^{\theta,\theta',\nu_1,\nu_2}_{t-s}[f])
- \EE^{\theta,\nu_1} g(\Pi^{\theta,\theta',\nu^{\theta}_0,\nu^{\theta'}_0}_{t-s}[f])\right|
 \leq C_5\,\tilde\EE\left[ 1\wedge \left(e^{\sum_{i=0}^{\lfloor t-s\rfloor-1}\xi_i} - 1\right)^{-1} \right],
\end{align*}
where $\{\xi_i\}$ are strictly positive i.i.d. random variables on a new probability space $(\tilde\Omega,\tilde{\mathcal{F}},\tilde\PP)$.

\smallskip

Thus, we have
\begin{align*}
&\frac{1}{T^2}\EE^{\theta,\nu^{\theta}_0}\int_{0\leq s\leq t\leq T} \left|\EE^{\theta,\nu_1} g(\Pi^{\theta,\theta',\nu_1,\nu_2}_{t-s}[f])
- \EE^{\theta,\nu_1} g(\Pi^{\theta,\theta',\nu^{\theta}_0,\nu^{\theta'}_0}_{t-s}[f])\right|_{\nu_1:=\hat{\pi}^{\theta,\nu^{\theta}_0}_s,\,\nu_2:=\hat{\pi}^{\theta',\nu^{\theta'}_0}_s}
\,dt\,ds\\
&\leq \frac{C_ 4}{T} + C_5\,\tilde \EE \frac{1}{T^2}\int_0^T \int_{0}^{t} 1\wedge \left(e^{\sum_{i=0}^{\lfloor t-s\rfloor-1}\xi_i} - 1\right)^{-1}ds\,dt.
\end{align*}
To estimate the second term in the right hand side of the above, we first notice that the expression inside the expectation is bounded by $1/2$. Hence, due to the dominated convergence theorem, it suffices to show that the latter expression converges to zero a.s. as $T\rightarrow\infty$.
To this end, we notice that
\begin{align*}
&\int_{0}^t 1\wedge \left(e^{\sum_{i=0}^{\lfloor t-s\rfloor-1}\xi_i} - 1\right)^{-1}ds
\leq 2 + \sum_{n=0}^{\lfloor t\rfloor} \left(e^{\sum_{i=0}^{n}\xi_i} - 1\right)^{-1}\\
&\leq 2 + \left(e^{\xi_0} - 1\right)^{-1}+ \sum_{n=1}^{\lfloor t\rfloor}\left(e^{\sum_{i=0}^{n}\xi_i} - e^{\sum_{i=1}^{n}\xi_i}\right)^{-1}\\
&\leq 2 + \left(e^{\xi_0} - 1\right)^{-1}\left(1+ \sum_{n=1}^{\lfloor t\rfloor} e^{-\sum_{i=1}^{n}\xi_i}\right)
\leq 2 + \left(e^{\xi_0} - 1\right)^{-1}\left(1+ \sum_{n=1}^{\infty} e^{-\sum_{i=1}^{n}\xi_i}\right),
\end{align*}
where the above infinite series is finite a.s., as it can be verified by a direct computation that it has a finite expectation.
All in all, we obtain
\begin{align*}
& \frac{1}{T^2}\int_0^T \int_{0}^t 1\wedge \left(e^{\sum_{i=0}^{\lfloor t-s\rfloor-1}\xi_i} - 1\right)^{-1}ds\,dt
 \leq \frac{1}{T}\left( 2 + \left(e^{\xi_0} - 1\right)^{-1}\left(1+ \sum_{n=1}^{\infty} e^{-\sum_{i=1}^{n}\xi_i}\right) \right),
\end{align*}
and the right hand side of the above converges to zero a.s. as $T\rightarrow\infty$, which completes the estimation of the first summand in the right hand side of \eqref{eq.vanishVar.pf.eq2}.

\smallskip

To estimate the second term in the right hand side of \eqref{eq.vanishVar.pf.eq2}, we use Lemma \ref{le:expStab.X} in Subsection \ref{subse:expStab.X} to obtain the existence of $(\tilde X,\bar X,\tilde B,\tilde Y,\bar Y)$ which are defined on a new probability space $(\tilde\Omega,\tilde{\mathcal{F}},\tilde\PP)$, with
\begin{align*}
&\tilde Y_t = \int_0^t h^\theta(\tilde X_s)\,ds + \tilde B_t,\quad \bar Y_t = \int_0^t h^\theta(\bar X_s)\,ds + \tilde B_t,\\
&\phantom{????????????????????????????????????????} \lim_{T\rightarrow\infty}\tilde\EE \int_0^T |h^\theta(\tilde X_s)-h^\theta(\bar X_s)|\,ds<\infty,
\end{align*}
with $\tilde B$ being a standard $\tilde\PP$-Brownian motion independent of $(\tilde X,\bar X)$, with $\tilde\PP\circ\tilde X^{-1}=\PP^{\theta,\nu_1}_X$, $\tilde\PP\circ\bar X^{-1}=\PP^{\theta,\nu^{\theta}_0}_X$, and with the above convergence being uniform over all $\nu_1\in\mathcal{P}([0,1)^q)$.
Thus, 
\begin{align*}
& \left|\EE^{\theta,\nu_1} g(\Pi^{\theta,\theta',\nu^{\theta}_0,\nu^{\theta'}_0}_{t-s}[f])
- \EE^{\theta,\nu^{\theta}_0} g(\Pi^{\theta,\theta',\nu^{\theta}_0,\nu^{\theta'}_0}_{t-s}[f])\right|\\
&\leq \tilde\EE\left| g(\langle A^\theta_{0,t-s}(\tilde Y_{[0,t-s]})[\nu^{\theta}_0],f^1\rangle,\langle A^{\theta'}_{0,t-s}(\tilde Y_{[0,t-s]})[\nu^{\theta'}_0],f^2\rangle)\right.\\
&\left.\quad\quad\quad- g(\langle A^\theta_{0,t-s}(\bar Y_{[0,t-s]})[\nu^{\theta}_0],f^1\rangle,\langle A^{\theta'}_{0,t-s}(\bar Y_{[0,t-s]})[\nu^{\theta'}_0],f^2\rangle)\right|,
\end{align*}
where $\langle\cdot,\cdot\rangle$ denotes the scalar product in $L^2(\RR^q)$ and we explicitly indicated the dependence of the operator $A$ on $\tilde Y$ and $\bar Y$. We continue, using the triangle inequality for $H$, as well as Proposition \ref{le:H.contract} and Lemma \ref{le:v.reg}:
\begin{align*}
& H(A^\theta_{0,t-s}(\tilde Y_{[0,t-s]})[\nu^\theta_0],A^\theta_{0,t-s}(\bar Y_{[0,t-s]})[\nu^{\theta}_0])\\
&\phantom{????????????}\leq H\left(A^\theta_{t-s-1,1}(\tilde Y_{[t-s-1,t-s]}-\tilde Y_{t-s-1})[A^\theta_{0,t-s-1}(\tilde Y_{[0,t-s-1]})[\nu^{\theta}_0]],\right.\\
& \phantom{?????????????????????????}\left.A^\theta_{t-s-1,1}(\tilde Y_{[t-s-1,t-s]}-\tilde Y_{t-s-1})[A^\theta_{0,t-s-1}(\bar Y_{[0,t-s-1]})[\nu^{\theta}_0]]\right)\\
&\phantom{????????????}+ H(A^\theta_{t-s-1,1}(\tilde Y_{[t-s-1,t-s]}-\tilde Y_{t-s-1})[A^\theta_{0,t-s-1}(\bar Y_{[0,t-s-1]})[\nu^{\theta}_0]],\\
&A^\theta_{t-s-1,1}(\bar Y_{[t-s-1,t-s]}-\bar Y_{t-s-1})[A^\theta_{0,t-s-1}(\bar Y_{[0,t-s-1]})[\nu^{\theta}_0]])\\
&\phantom{????????????}\leq e^{-\tilde\gamma(\mathrm{Osc}_{t-s-1,t-s}(\tilde B))}\,H(A^\theta_{0,t-s-1}(\tilde Y_{[0,t-s-1]})[\nu^{\theta}_0],A^\theta_{0,t-s-1}(\bar Y_{[0,t-s-1]})[\nu^{\theta}_0])\\
&\phantom{????????????}+ 1\wedge\left(\frac{1}{\tilde\epsilon(\mathrm{Osc}_{t-s-1,t-s}(\tilde B))}\,\int_{t-s-1}^{t-s}|h^{\theta}(\tilde X_r)-h^{\theta}(\bar X_r)|\,dr\right),
\end{align*}
where $\tilde\epsilon$ and $\tilde\gamma$ are strictly positive and non-increasing in each variable functions (the same for all $\nu_1\in\mathcal{P}([0,1)^q)$), and where we identified (with a slight abuse of notation) certain probability density functions with their associated measures.
Using the notation $\xi_i:=\tilde\gamma(\mathrm{Osc}_{t-s-i,t-s-i+1}(\tilde B))$ and $\eta_i:=\tilde\epsilon(\mathrm{Osc}_{t-s-i,t-s-i+1}(\tilde B))$, we iterate the above to obtain
\begin{align*}
& H(A^\theta_{0,t-s}(\tilde Y_{[0,t-s]})[\nu^{\theta}_0],A^\theta_{0,t-s}(\bar Y_{[0,t-s]})[\nu^{\theta}_0])
\leq e^{-\sum_{i=1}^{\lfloor t-s\rfloor}\xi_i}\\
&+\sum_{n=1}^{\lfloor t-s\rfloor}e^{-\sum_{i=1}^{n-1}\xi_i}\,\left[ 1\wedge\left(\frac{1}{\eta_n}\,\int_{t-s-n}^{t-s-n+1}|h^{\theta}(\tilde X_r)-h^{\theta}(\bar X_r)|\,dr\right)\right].
\end{align*}
Next, we recall that $\{\eta_n\}$ are i.i.d. and independent of $(\tilde X,\bar X)$, to deduce, for any $\varepsilon,\delta>0$:
\begin{align*}
& \tilde\PP \left(\frac{1}{\eta_n}\,\int_{t-s-n}^{t-s-n+1}|h^{\theta}(\tilde X_r)-h^{\theta}(\bar X_r)|\,dr\geq\varepsilon\right)\\
&= \int_0^\infty \tilde\PP \left(\int_{t-s-n}^{t-s-n+1}|h^{\theta}(\tilde X_r)-h^{\theta}(\bar X_r)|\,dr\geq\varepsilon\,x\right)\,\tilde\PP(\eta_1\in dx)\\
&\leq \tilde\PP(\eta_1<\delta) + \frac{1}{\varepsilon\,\delta}\,\tilde\PP(\eta_1\geq\delta)\,\tilde\EE \int_{t-s-n}^{t-s-n+1}|h^{\theta}(\tilde X_r)-h^{\theta}(\bar X_r)|\,dr,\\
\end{align*}
where we used Chebyshev's inequality.
Lemma \ref{le:expStab.X} implies that, for any $\epsilon>0$, there exist $t_0,\delta>0$ such that the supremum of the right hand side of the above over all $\nu_1\in\mathcal{P}([0,1)^q)$ does not exceed $\epsilon$ whenever $t-s-n\geq t_0$.
The latter yields, for any $n>0$,
\begin{align*}
&J_{n,t-s}:=\tilde\EE \left[1\wedge\left(\frac{1}{\eta_n}\,\int_{t-s-n}^{t-s-n+1}|h^{\theta}(\tilde X_r)-h^{\theta}(\bar X_r)|\,dr\right)\right] \rightarrow 0\text{ as }t-s\rightarrow\infty,
\end{align*}
uniformly over all $\nu_1\in\mathcal{P}([0,1)^q)$.
Then, using the Cauchy-Schwartz inequality, we obtain:
\begin{align*}
& \tilde\EE \,H(A^\theta_{0,t-s}(\tilde Y_{[0,t-s]})[\nu^{\theta}_0],A^\theta_{0,t-s}(\bar Y_{[0,t-s]})[\nu^{\theta}_0])\\
&\phantom{????????????????????????????????}\leq \left(\tilde\EE e^{-\xi_1}\right)^{\lfloor t-s\rfloor}
+\sum_{n=1}^{\infty} \left(\tilde\EE e^{-2\xi_1}\right)^{(n-1)/2}\,\sqrt{J_{n,t-s}}.
\end{align*}
The dominated convergence theorem implies that the right hand side of the above vanishes as $t-s\rightarrow\infty$, uniformly over all $\nu_1\in\mathcal{P}([0,1)^q)$.

\smallskip

Next, the inequality $\|\mu_1-\mu_2\|_{\mathrm{TV}}\leq 2H(\mu_1,\mu_2)/(1-H(\mu_1,\mu_2))$, as well as the uniform boundedness of $\|\cdot\|_{\mathrm{TV}}$ and of $H$, imply
\begin{align*}
& \lim_{r\rightarrow\infty}\tilde\EE \,\|A^\theta_{0,r}(\tilde Y_{[0,r]})[\nu^{\theta}_0] - A^\theta_{0,r}(\bar Y_{[0,r]})[\nu^{\theta}_0]\|_{\mathrm{TV}}
=0,
\end{align*}
uniformly over all $\nu_1\in\mathcal{P}([0,1)^q)$.
Similarly, we deduce that the above holds with $\theta$ replaced by $\theta'$.
Therefore,
\begin{align*}
& \lim_{r\rightarrow\infty}\sup_{\nu_1\in\mathcal{P}([0,1)^q)}\left|\EE^{\theta,\nu_1} g(\Pi^{\theta,\theta',\nu^{\theta}_0,\nu^{\theta'}_0}_{r}[f])
- \EE^{\theta,\nu^{\theta}_0} g(\Pi^{\theta,\theta',\nu^{\theta}_0,\nu^{\theta'}_0}_{r}[f])\right|=0.
\end{align*}
Thus,
\begin{align*}
&\frac{1}{T^2}\EE^{\theta,\nu^{\theta}_0}\int_{0\leq s\leq t\leq T} \left|\EE^{\theta,\nu_1} g(\Pi^{\theta,\theta',\nu^{\theta}_0,\nu^{\theta'}_0}_{t-s}[f])
- \EE^{\theta,\nu^{\theta}_0} g(\Pi^{\theta,\theta',\nu^{\theta}_0,\nu^{\theta'}_0}_{t-s}[f])\right|_{\nu_1:=\hat{\pi}^{\theta,\nu^{\theta}_0}_s}
\,dt\,ds\\
&\leq \EE^{\theta,\nu^{\theta}_0} \frac{1}{T}\int_{0}^{T} \sup_{\nu_1\in\mathcal{P}([0,1)^q)}\left|\EE^{\theta,\nu_1} g(\Pi^{\theta,\theta',\nu^{\theta}_0,\nu^{\theta'}_0}_{r}[f])
- \EE^{\theta,\nu^{\theta}_0} g(\Pi^{\theta,\theta',\nu^{\theta}_0,\nu^{\theta}_0}_{r}[f])\right|\,dr.
\end{align*}
It remains to notice that the integrand in the right hand side of the above display is absolutely bounded by a constant and that it converges to zero as $r\rightarrow\infty$, $\PP^{\theta,\nu^{\theta}_0}$-a.s.. Applying the dominated convergence theorem, we complete the proof.

\end{proof}

\smallskip

The following proposition shows that the normalized finite-variation component of the log-likelihood ratio converges in probability, for any initial value of the filter.

\begin{proposition}\label{prop:strong.ergodicity}
For any $\nu,\nu'\in\mathcal{P}([0,1)^q)$, $\theta,\theta'\in\Theta$, and any $\varepsilon>0$,
\begin{align*}
& \lim_{T\rightarrow\infty}\PP^{\theta,\nu^\theta_0}\left(\left|\frac{1}{T} \int_0^T\left|\hat{\pi}^{\theta,\nu}_t[h^{\theta}] - \hat{\pi}^{\theta',\nu'}_t[h^{\theta'}]\right|^2\,dt - \Lambda(\theta,\theta')\right|\geq\varepsilon\right) = 0.
\end{align*}
\end{proposition}
\begin{proof}
For $\nu=\nu^\theta_0$ and $\nu'=\nu^{\theta'}_0$, the statement of the proposition follows directly from Proposition \ref{prop:exp.conv.sq.error} and from Lemma \ref{le:vanishing.variance} applied to 
\begin{align*}
g(\mu_1[f^1],\mu_2[f^2])= (\mu_1[h^\theta_i]-\mu_2[h^{\theta'}_i])^2,\quad i=1,\ldots,q.
\end{align*}

\smallskip

Next, we choose any $\nu,\nu'\in\mathcal{P}([0,1)^q)$ and apply Lemma \ref{cor:H.contracts.on.pi} to obtain
\begin{align}
& \left|\left|\hat{\pi}^{\theta,\nu}_t[h^{\theta}] - \hat{\pi}^{\theta',\nu'}_t[h^{\theta'}]\right|^2
-\left|\hat{\pi}^{\theta,\nu^\theta_0}_t[h^{\theta}] - \hat{\pi}^{\theta',\nu^{\theta'}_0}_t[h^{\theta'}]\right|^2\right|\label{eq.filter.stronglyErgodic.pf.eq1}\\
&\leq C_1\,\Gamma_{1,t-1}\,\left(\frac{H(\hat{\pi}^{\theta,\nu}_1,\hat{\pi}^{\theta,\nu^\theta_0}_1)}{1-H(\hat{\pi}^{\theta,\nu}_1,\hat{\pi}^{\theta,\nu^\theta_0}_1)} 
+ \frac{H(\hat{\pi}^{\theta',\nu'}_1,\hat{\pi}^{\theta',\nu^{\theta'}_0}_1)}{1-H(\hat{\pi}^{\theta',\nu'}_1,\hat{\pi}^{\theta',\nu^{\theta'}_0}_1)}\right),\nonumber
\end{align}
$\PP^{\theta,\nu^\theta_0}$-a.s..
Notice also that
\begin{align*}
&\frac{1}{T} \int_1^T \Gamma_{1,t-1}\,dt \leq \frac{1}{T} \sum_{n=0}^{\lfloor T\rfloor}e^{-\sum_{i=0}^{n}\gamma(C_2+\mathrm{Osc}_{i+1,i+2}(\tilde{W}^{\theta,\nu^\theta_0}))},
\end{align*}
and that the latter sum converges to a finite limit a.s., as $T\rightarrow\infty$ (since the expectation of this limit is finite, as can be verified directly).
This observation implies that the right hand side of \eqref{eq.filter.stronglyErgodic.pf.eq1} converges to zero $\PP^{\theta,\nu^\theta_0}$-a.s. and completes the proof of the proposition.
\end{proof}

\subsection{Uniform robustness of the filter}
\label{subse:robust}
In this subsection, we show that the normalized finite-variation component of the log-likelihood ratio, $\frac{1}{t}\int_0^t \left|\hat{\pi}^{\theta,\nu}_s[h^\theta]-\hat{\pi}^{\theta',\nu}_s[h^{\theta'}]\right|^2ds$, converges uniformly over $\theta'$, as $t\rightarrow\infty$, and conclude the proof of Theorem \ref{thm:main}. 

\smallskip

We begin with the proof of Proposition \ref{le:uniform.robust.filter}.

{\it Step 1.} In the first step, we estimate the modulus of continuity of the filter w.r.t. $\theta$, as measured by the Hilbert projective metric, uniformly in time.
We fix $\nu$, $\theta$ and denote by $\omega_1$ the modulus of uniform continuity in Assumption \ref{ass:unifcont}:
\begin{align*}
\omega_1(\delta):=\sup_{\bar d(\theta,\theta')\leq \delta}\left(\|a^\theta-a^{\theta'}\|_{C^1}+\|b^\theta-b^{\theta'}\|_{C}+\|h^\theta-h^{\theta'}\|_{C^{2}}\right).
\end{align*}
With any $t\geq 1$ we associate the unique $t_0\in[1,2)$ and $n\in\mathbb{N}\cup\{0\}$ s.t. $t=t_0+n$. Then, applying the triangle inequality to the Hilbert metric $H$, we obtain:
\begin{align*}
& H(\hat\pi^{\theta',\nu}_t,\hat\pi^{\theta'',\nu}_t)
= H(A^{\theta'}_{t-1,1}[\hat\pi^{\theta',\nu}_{t-1}],A^{\theta''}_{t-1,1}[\hat\pi^{\theta'',\nu}_{t-1}])\\
&\leq H(A^{\theta'}_{t-1,1}[\hat\pi^{\theta',\nu}_{t-1}],A^{\theta'}_{t-1,1}[\hat\pi^{\theta'',\nu}_{t-1}])
+ H(A^{\theta'}_{t-1,1}[\hat\pi^{\theta'',\nu}_{t-1}],A^{\theta''}_{t-1,1}[\hat\pi^{\theta'',\nu}_{t-1}])\\
& \leq H(\hat\pi^{\theta',\nu}_{t-1},\hat\pi^{\theta'',\nu}_{t-1})\,e^{-\gamma(\mathrm{Osc}_{t-1,t}(Y))}
+ H(A^{\theta'}_{t-1,1}[\hat\pi^{\theta'',\nu}_{t-1}],A^{\theta''}_{t-1,1}[\hat\pi^{\theta'',\nu}_{t-1}]),
\end{align*}
where the last inequality is due to \eqref{eq.exp.stab.local}.
Iterating the above, we obtain
\begin{align*}
 H(\hat\pi^{\theta',\nu}_t,\hat\pi^{\theta'',\nu}_t)
\leq &H(A^{\theta'}_{t-1,1}[\hat\pi^{\theta'',\nu}_{t-1}],A^{\theta''}_{t-1,1}[\hat\pi^{\theta'',\nu}_{t-1}])\\
+\sum_{i=1}^{n-1}& H(A^{\theta'}_{t_0+i-1,1}[\hat\pi^{\theta'',\nu}_{t_0+i-1}],A^{\theta''}_{t_0+i-1,1}[\hat\pi^{\theta'',\nu}_{t_0+i-1}])\,e^{-\sum_{j=i}^{n-1} \gamma(\mathrm{Osc}_{t_0+j,t_0+j+1}(Y))}\\
+H&(\hat\pi^{\theta',\nu}_{t_0},\hat\pi^{\theta'',\nu}_{t_0})\,e^{-\sum_{j=0}^{n} \gamma(\mathrm{Osc}_{t_0+j,t_0+j+1}(Y))}.
\end{align*}
For $t, \delta\geq 0$, we define the $\sigma((Y_{r}-Y_{t})_{r\in[t,t+1]})$-measurable 
random variable
\begin{align*}
\eta^\delta_t:=\sup_{\nu'\in\mathcal{P}([0,1)^q),\, \bar d(\theta',\theta'')\leq\delta} H(A^{\theta'}_{t,1}[\nu'],A^{\theta''}_{t,1}[\nu']),
\end{align*}
and, for $s\in[1,2)$, the $\cF^Y_s$-measurable random variable
\begin{align*}
\tilde \eta_s^\delta:=\sup_{\nu'\in\mathcal{P}([0,1)^q),\, \bar d(\theta',\theta'')\leq\delta} H(A^{\theta'}_{0,s}[\nu'],A^{\theta''}_{0,s}[\nu']),
\end{align*}
so that
\begin{align*}
 &\sup_{ \bar d(\theta',\theta'')\leq \delta}H(\hat\pi^{\theta',\nu}_t,\hat\pi^{\theta'',\nu}_t)\\
&\leq \tilde \eta^\delta_{t_0} e^{-\sum_{j=0}^{n} \gamma(\mathrm{Osc}_{t_0+j,t_0+j+1}(Y))}+ \eta^{\delta}_{t-1}+\sum_{i=1}^{n-1} \eta^{\delta}_{t_0+i-1}e^{-\sum_{j=i}^{n-1}\gamma(\mathrm{Osc}_{t_0+j,t_0+j+1}(Y))}. 
\end{align*}
Due to the boundedness of $h$, there exists a constant $C$ such that
\begin{align*} 
\mathrm{Osc}_{s,s+1}(Y)\leq C+\mathrm{Osc}_{s,s+1}(W),
\end{align*}
where we denote $W:=\tilde W^{\theta,\nu^\theta_0}$.
Using the monotonicity of $\gamma$, we obtain
\begin{align}
 &\sup_{ \bar d(\theta',\theta'')\leq \delta}H(\hat\pi^{\theta',\nu}_t,\hat\pi^{\theta'',\nu}_t)\label{eq:est1}\\
&\leq \tilde \eta^\delta_{t_0} e^{-\sum_{j=0}^{n} \gamma(C+\mathrm{Osc}_{t_0+j,t_0+j+1}(W))}+ \eta^{\delta}_{t-1}+\sum_{i=1}^{n-1} \eta^{\delta}_{t_0+i-1}e^{-\sum_{j=i}^{n-1}\gamma(C+\mathrm{Osc}_{t_0+j,t_0+j+1}(W))}.\nonumber
\end{align}
Note that $\tilde \eta^\delta_{t_0}$ is independent of $\sum_{j=0}^{n} \gamma(C+\mathrm{Osc}_{t_0+j,t_0+j+1}(W))$, that $\eta^{\delta}_{t_0+i-1}$ is independent of $\sum_{j=i}^{n-1}\gamma(C+Osc_{t_0+j,t_0+j+1}(W))$, and that the family $\{\gamma(C+\mathrm{Osc}_{n,n+1}(W))\}_{n\geq0}$ is i.i.d.. Defining
\begin{align*}
\bar \gamma:=-\ln \EE^{\theta,\nu^\theta_0}\left[e^{-\gamma(C+\mathrm{Osc}_{t_0,t_0+1}(W))}\right]=-\ln \EE^{\theta,\nu^\theta_0}\left[e^{-\gamma(C+\mathrm{Osc}_{0,1}(W))}\right]>0
\end{align*} 
and taking the expectation in \eqref{eq:est1}, we obtain
\begin{align*}
 &\EE^{\theta,\nu^\theta_0}\left[\sup_{ \bar d(\theta',\theta'')\leq \delta}H(\hat\pi^{\theta',\nu}_t,\hat\pi^{\theta'',\nu}_t)\right]\\
&\leq e^{-\bar\gamma(n+1)}\sup_{s\in[1,2)}\EE^{\theta,\nu^\theta_0}\left[\tilde \eta^\delta_{s}\right] + \EE^{\theta,\nu^\theta_0}\left[\eta^{\delta}_{t-1}\right]+\sum_{i=1}^{n-1} \EE^{\theta,\nu^\theta_0}\left[\eta^{\delta}_{t_0+i-1}\right]e^{-\bar\gamma(n-i)},
\end{align*}
so that 
\begin{align}
 &\sup_{t\geq 1}\EE^{\theta,\nu^\theta_0}\left[\sup_{ \bar d(\theta',\theta'')\leq \delta}H(\hat\pi^{\theta',\nu}_t,\hat\pi^{\theta'',\nu}_t)\right]\leq \sup_{s\in[1,2)}\EE^{\theta,\nu^\theta_0}\left[\tilde \eta^\delta_{s}\right]+\frac{\sup_{t\geq 1}\EE^{\theta,\nu^\theta_0}\left[\eta^{\delta}_{t}\right]}{1-e^{-\bar\gamma}}.
 \label{eq.robust.mainProp.Step1.result}
\end{align}

\smallskip

{\it Step 2.} In this step, we show the convergence $\eta^{\delta}_{t}\to0$, as $\delta\to 0$, and complete the proof of \eqref{eq.uniform.robust.filter}.
We choose a constant $C>0$ as above, so that $\mathrm{Osc}_{t,t+1}(Y) \leq  C + \mathrm{Osc}_{t,t+1}(W)$, and recall $\epsilon(T,x)$ defined in Lemma \ref{le:v.reg}. In order to estimate $H(A^{\theta'}_{t,1}[\nu'],A^{\theta''}_{t,1}[\nu'])$, we denote 
\begin{align*}
R:=\frac{\int_{[0,1)^q}\hat A^{\theta''}_{t,1}[\nu'](y)dy}{\int_{[0,1)^q}\hat A^{\theta'}_{t,1}[\nu'](y)dy}\geq 0
\end{align*}
and obtain
\begin{align*}
&\frac{1}{R}\left|R-\frac{A^{\theta'}_{t,1}[\nu'](z)}{A^{\theta''}_{t,1}[\nu'](z)}\right|= \left|1-\frac{\hat A^{\theta'}_{t,1}[\nu'](z)}{\hat A^{\theta''}_{t,1}[\nu'](z)}\right|\\
&\leq  \frac{  \left|\hat A^{\theta'}_{t,1}[\nu'](z)-\hat A^{\theta''}_{t,1}[\nu'](z)\right|}{ \hat A^{\theta''}_{t,1}[\nu'](z)}\leq \frac{  \left| \hat A^{\theta'}_{t,1}[\nu'](z)-\hat A^{\theta''}_{t,1}[\nu'](z)\right|}{\epsilon(1, C + \mathrm{Osc}_{t,t+1}(W))},
\end{align*}
where we used the lower bound in \eqref{eq:boundK}.
By the definition of $\hat A^{\theta'}_{t,1}[\nu']$ and due to \eqref{eq:contdv}, we have 
\begin{align*}
|\hat A^{\theta'}_{t,1}[\nu'](z)-\hat A^{\theta''}_{t,1}[\nu'](z)|&\leq\Big|\int \hat{K}^{\theta'}_{t,1}(x,z)-\hat{K}^{\theta''}_{t,1}(x,z)\nu'(dx)\Big|\\
&\leq \sup_{x\in [0,1)^q}|\hat{K}^{\theta'}_{t,1}(x,z)-\hat{K}^{\theta''}_{t,1}(x,z)|
\leq \frac{\omega_1 (\bar d(\theta',\theta''))}{\epsilon(1, C + \mathrm{Osc}_{t,t+1}(W))}.
\end{align*}
Thus, 
    \begin{align*}
        &\left|1-\frac{A^{\theta'}_{t,1}[\nu'](z)}{R\,A^{\theta''}_{t,1}[\nu'](z)}\right|\leq \frac{\omega_1 (\bar d(\theta',\theta''))}{\epsilon^2(1, C + \mathrm{Osc}_{t,t+1}(W))}.
    \end{align*}
Next, we use the definition of $H$, to obtain 
\begin{align*}
&0\leq H(A^{\theta'}_{t,1}[\nu'],A^{\theta''}_{t,1}[\nu'])=H(A^{\theta'}_{t,1}[\nu'],R\,A^{\theta''}_{t,1}[\nu'])\leq 1-\frac{\inf_z \frac{A^{\theta'}_{t,1}[\nu'](z)}{R\,A^{\theta''}_{t,1}[\nu'](z)}}{\sup_z \frac{A^{\theta'}_{t,1}[\nu'](z)}{R\,A^{\theta''}_{t,1}[\nu'](z)}}\\
&=1-\frac{1-(1-\inf_z \frac{A^{\theta'}_{t,1}[\nu'](z)}{R\,A^{\theta''}_{t,1}[\nu'](z)})}{1-(1-\sup_z \frac{A^{\theta'}_{t,1}[\nu'](z)}{R\,A^{\theta''}_{t,1}[\nu'](z)})}
\leq 1-\frac{1- \frac{\omega_1 (\bar d(\theta',\theta''))}{\epsilon^2(1, C + \mathrm{Osc}_{t,t+1}(W))}}{1+ \frac{\omega_1 (\bar d(\theta',\theta''))}{\epsilon^2(1, C + \mathrm{Osc}_{t,t+1}(W))}}\\
&\phantom{????????????????????????????????????}
\leq { \frac{2\omega_1 (\bar d(\theta',\theta''))}{\epsilon^2(1, C + \mathrm{Osc}_{t,t+1}(W))}},
\end{align*}
which yields 
\begin{align*}
0\leq \eta^\delta_t\leq 1\wedge  \frac{2\omega_1 (\bar d(\theta',\theta''))}{\epsilon^2(1, C + \mathrm{Osc}_{t,t+1}(W))} .  
\end{align*}

Next, we fix an arbitrary $N>0$ and obtain: 
\begin{align*}
\EE^{\theta,\nu^\theta_0}\left[\eta^{\delta}_{t}\right]&\leq \EE^{\theta,\nu^\theta_0}\left[\eta^{\delta}_{t}\1_{\{\mathrm{Osc}_{t,t+1}(W)<N\}}\right]+\PP^{\theta,\nu^\theta_0}\left[\mathrm{Osc}_{t,t+1}(W)>N\right]\\
    &\leq\frac{2\omega_1(\delta)}{\epsilon^2(1,C+N)}+\PP^{\theta,\nu^\theta_0}\left[\mathrm{Osc}_{0,1}(W)>N\right].
\end{align*}
Using the above, we easily deduce that
\begin{align*}
\sup_{t\geq 0}\EE^{\theta,\nu^\theta_0}\left[\eta^{\delta}_{t}\right]\leq \tilde \omega(\delta),
\end{align*}
for some non-decreasing $\tilde \omega:\RR_+\rightarrow \RR_+$ with $\tilde \omega(0)=0$.
Similarly we obtain
\begin{align*} 
\sup_{s\in[1,2)}\EE^{\theta,\nu^\theta_0}\left[\tilde \eta^\delta_{s}\right]\leq\tilde \omega(\delta),
\end{align*}
which, in view of \eqref{eq.robust.mainProp.Step1.result}, yields
\begin{align*}
\sup_{t\geq 1}\EE^{\theta,\nu^\theta_0}\left[\sup_{ \bar d(\theta',\theta'')\leq \delta}H(\hat\pi^{\theta',\nu}_t,\hat\pi^{\theta'',\nu}_t)\right]\to 0\mbox{ as }\delta\downarrow 0.
\end{align*}
To obtain \eqref{eq.uniform.robust.filter} from the above display, we note that there exists a constant $C_1>0$ such that
\begin{align*}
& \lim_{\delta\downarrow0} \sup_{t\geq1} \EE^{\theta,\nu^\theta_0}\sup_{ \bar d(\theta',\theta'')\leq\delta} |\pi^{\theta',\nu}_t[h^{\theta'}] - \pi^{\theta'',\nu}_t[h^{\theta''}]|
\leq \lim_{\delta\downarrow0}\sup_{ \bar d(\theta',\theta'')\leq\delta}\|h^{\theta'}-h^{\theta''}\|\\
&+ C_1\,\lim_{\delta\downarrow0}\sup_{t\geq1} \EE^{\theta,\nu^\theta_0}\sup_{ \bar d(\theta',\theta'')\leq\delta}H(\hat\pi^{\theta',\nu}_t,\hat\pi^{\theta'',\nu}_t)= 0,
\end{align*}
which implies \eqref{eq.uniform.robust.filter} and completes the proof of Proposition \ref{le:uniform.robust.filter}.

\medskip

The following corollary is an immediate consequence of Proposition \ref{le:uniform.robust.filter} and of the boundedness of $h$.

\begin{cor}\label{cor:uniform.cont.Z}
For any $\theta\in\Theta$ and $\nu\in\mathcal{P}([0,1)^q)$,
\begin{align*}
&\lim_{\delta\downarrow0}\sup_{T\geq0}\EE^{\theta,\nu^\theta_0}\sup_{ \bar d(\theta',\theta'')\leq\delta} \left|\frac{1}{T}\int_0^T \left[\left|\hat{\pi}^{\theta,\nu}_t[h^{\theta}] - \hat{\pi}^{\theta',\nu}_t[h^{\theta'}]\right|^2 - \left|\hat{\pi}^{\theta,\nu}_t[h^{\theta}] - \hat{\pi}^{\theta'',\nu}_t[h^{\theta''}]\right|^2\right]\,dt\right|=0.
\end{align*}
\end{cor}

Using the above corollary and Propositions \ref{le:uniform.robust.filter}, \ref{prop:strong.ergodicity}, \ref{prop:Lambda.pos}, we prove the following result.

\begin{cor}\label{cor:Lambda.bdd.away.from.zero}
For any $\theta\in\Theta$ and $\epsilon>0$, we have: $\inf_{\theta':\, \bar d(\theta',\overline{\Theta}(\theta))\geq\epsilon}\Lambda(\theta,\theta')>0$.
\end{cor}
\begin{proof}
We argue by contradiction and assume that there exist $\theta\in\Theta$ and $\epsilon>0$ such that $\inf_{\theta':\, \bar d(\theta',\overline{\Theta}(\theta))\geq\epsilon}\Lambda(\theta,\theta')=0$. Then, due to the compactness of $\Theta$ (which yields the compactness of $\{\theta':\, \bar d(\theta',\overline{\Theta}(\theta))\geq\epsilon\}$), there exist $\theta_0\in\Theta$, with $\bar d(\theta_0,\overline{\Theta}(\theta))\geq\epsilon$, and a sequence $\theta_k\rightarrow\theta_0$ such that $\Lambda(\theta,\theta_k)\rightarrow0$.
Then,
\begin{align*}
&|\Lambda(\theta,\theta_0)-\Lambda(\theta,\theta_k)|\leq \EE^{\theta,\nu^\theta_0}\left|\Lambda(\theta,\theta_0) - \frac{1}{T}\int_0^T\left|\hat\pi^{\theta,\nu^\theta_0}_t[h^{\theta}] - \hat\pi^{\theta_0,\nu^{\theta}_0}_t[h^{\theta_0}]\right|^2\,dt\right|\\
&+ \EE^{\theta,\nu^\theta_0}\left|\Lambda(\theta,\theta_k) - \frac{1}{T}\int_0^T\left|\hat\pi^{\theta,\nu^\theta_0}_t[h^{\theta}]-\hat\pi^{\theta_k,\nu^{\theta}_0}_t[h^{\theta_k}]\right|^2\,dt\right|\\
& + \EE^{\theta,\nu^\theta_0}\left| \frac{1}{T}\int_0^T\left[\left|\hat\pi^{\theta,\nu^\theta_0}_t[h^{\theta}] - \hat\pi^{\theta_0,\nu^{\theta}_0}_t[h^{\theta_0}]\right|^2
- \left|\hat\pi^{\theta,\nu^\theta_0}_t[h^{\theta}] - \hat\pi^{\theta_k,\nu^{\theta}_0}_t[h^{\theta_k}]\right|^2\right]\,dt\right|.
\end{align*}
For any $\varepsilon>0$, we can choose large enough $k$, such that the last term in the right hand side of the above is bounded by $\varepsilon/3$, due to Corollary \ref{cor:uniform.cont.Z}. For any such $k$, we can choose a large enough $T>0$, such that the first two terms in the right hand side of the above are bounded by $\varepsilon/3$ each, due to Proposition \ref{prop:strong.ergodicity}. Thus, for any $\varepsilon>0$, there exists $k>0$ such that $|\Lambda(\theta,\theta_0)-\Lambda(\theta,\theta_k)|\leq \varepsilon$. It only remains to apply Proposition \ref{prop:Lambda.pos} to deduce that $\Lambda(\theta,\theta_0)>0$ and obtain the desired contradiction.
\end{proof}

\medskip

Next, we establish another auxiliary result.

\begin{proposition}\label{le:uniform.robust.filter.2}
For any $\theta\in\Theta$, $\nu\in\mathcal{P}([0,1)^q)$ and $\epsilon>0$,
\begin{align}
&\lim_{T\rightarrow\infty} \PP^{\theta,\nu^\theta_0}\left(\sup_{\theta'\in\Theta} \frac{1}{T}\left|\int_0^T\hat{\pi}^{\theta',\nu}_t[h^{\theta'}]^\top d\tilde W^{\theta,\nu^\theta_0}_t\right|
\geq\epsilon\right)=0.\label{eq.mtg.timeAve.small}
\end{align}
\end{proposition}
\begin{proof}
First, we notice that
\begin{align}
& \log L^{\theta',\nu}_T = -\frac{1}{2}\int_0^T|\hat\pi^{\theta',\nu}_t[h^{\theta'}]|^2 dt + \int_0^T \hat\pi^{\theta',\nu}_t[h^{\theta'}]^\top \left(\hat\pi^{\theta,\nu_0^\theta}_t[h^{\theta}]\,dt + d W_t\right)\nonumber\\
&= \int_0^T\left(\hat\pi^{\theta',\nu}_t[h^{\theta'}]^\top \hat\pi^{\theta,\nu_0^\theta}_t[h^{\theta}] - \frac{1}{2}|\hat\pi^{\theta',\nu}_t[h^{\theta'}]|^2\right) dt + \int_0^T \hat\pi^{\theta',\nu}_t[h^{\theta'}]^\top dW_t,\nonumber\\
& \frac{1}{T}\int_0^T(\hat\pi^{\theta',\nu}_t[h^{\theta'}]-\hat\pi^{\theta'',\nu}_t[h^{\theta''}])^\top dW_t
= \frac{1}{T}(\log L^{\theta',\nu}_T - \log L^{\theta'',\nu}_T) \label{eq.rev.step3.eq1}\\
&+ \frac{1}{T}\int_0^T\left((\hat\pi^{\theta'',\nu}_t[h^{\theta''}]-\hat\pi^{\theta',\nu}_t[h^{\theta'}])^\top \hat\pi^{\theta,\nu_0^\theta}_t[h^{\theta}] + \frac{1}{2}(|\hat\pi^{\theta',\nu}_t[h^{\theta'}]|^2-|\hat\pi^{\theta'',\nu}_t[h^{\theta''}]|^2)\right) dt.\nonumber
\end{align}

\smallskip

Let us estimate the first term in the right hand side of \eqref{eq.rev.step3.eq1}. To this end, we recall \eqref{eq:repL}, to obtain:
\begin{align*}
& \frac{1}{T}|\log L^{\theta',\nu}_T - \log L^{\theta'',\nu}_T|
\leq \frac{1}{T}|\log L^{\theta',\nu}_{\lfloor T\rfloor-1}/L^{\theta'',\nu}_{\lfloor T\rfloor-1}|\\
&\phantom{?????????????????????????????????}
+ \frac{1}{T} |\log L^{\theta',\nu}_{\lfloor T\rfloor-1,T-\lfloor T\rfloor+1}/L^{\theta'',\nu}_{\lfloor T\rfloor-1,T-\lfloor T\rfloor+1}|,\\
&\frac{1}{T}\EE^{\theta,\nu_0^\theta} \sup_{\bar d(\theta',\theta'')\leq\delta} |\log L^{\theta',\nu}_{\lfloor T\rfloor-1}/L^{\theta'',\nu}_{\lfloor T\rfloor-1}|
\leq \frac{1}{T} \sum_{i=1}^{\lfloor T\rfloor-1} \EE^{\theta,\nu_0^\theta} \sup_{\bar d(\theta',\theta'')\leq\delta} |\log L^{\theta',\nu}_{i-1,1}/L^{\theta'',\nu}_{i-1,1}|\\
&= \frac{1}{T} \sum_{i=1}^{\lfloor T\rfloor-1} \EE^{\theta,\nu_0^\theta} \sup_{\bar d(\theta',\theta'')\leq\delta} \left|\log \frac{\int_{[0,1)^q} U^{\theta',1}_{f^{\theta'}_{i-1},g^{\theta'}_{i-1}}[1/E^{\theta'}_{i-1,i}](0,x)\pi^{\theta',\nu}_{i-1}(dx)}{\int_{[0,1)^q} U^{\theta'',1}_{f^{\theta''}_{i-1},g^{\theta''}_{i-1}}[1/E^{\theta''}_{i-1,i}](0,x)\pi^{\theta'',\nu}_{i-1}(dx)}\right|\\
&\leq \frac{1}{T} \sum_{i=1}^{\lfloor T\rfloor-1} \EE^{\theta,\nu_0^\theta}\sup_{\bar d(\theta',\theta'')\leq\delta,\,x\in[0,1)^q}\left|\log \frac{U^{\theta',1}_{f^{\theta'}_{i-1},g^{\theta'}_{i-1}}[1/E^{\theta'}_{i-1,i}](0,x)}{U^{\theta'',1}_{f^{\theta''}_{i-1},g^{\theta''}_{i-1}}[1/E^{\theta''}_{i-1,i}](0,x)}\right|\\
&= \frac{\lfloor T\rfloor-1}{T} \EE^{\theta,\nu_0^\theta} \sup_{\bar d(\theta',\theta'')\leq\delta,\,x\in[0,1)^q}\left|\log \frac{U^{\theta',1}_{f^{\theta'}_{0},g^{\theta'}_{0}}[1/E^{\theta'}_{0,1}](0,x)}{U^{\theta'',1}_{f^{\theta''}_{0},g^{\theta''}_{0}}[1/E^{\theta''}_{0,1}](0,x)}\right|,
\end{align*}
where the last equality is due to the fact that the increments of $Y$ are stationary under $\PP^{\theta,\nu_0^\theta}$.

Next, we recall \eqref{eq.rev.appendix.eq1} to obtain
\begin{align*}
& U^{\theta',1}_{f^{\theta'}_{0},g^{\theta'}_{0}}[1/E^{\theta'}_{0,1}](0,x) = \int_{[0,1)^q} \hat K^{\theta'}_{0,1}(x,z)\,dz.
\end{align*}
Using the above and Lemma \ref{le:v.reg}, we deduce that
\begin{align*}
& \sup_{\bar d(\theta',\theta'')\leq\delta,\,x\in[0,1)^q}\left|\log \frac{U^{\theta',1}_{f^{\theta'}_{0},g^{\theta'}_{0}}[1/E^{\theta'}_{0,1}](0,x)}{U^{\theta'',1}_{f^{\theta''}_{0},g^{\theta''}_{0}}[1/E^{\theta''}_{0,1}](0,x)}\right| \rightarrow0,
\end{align*}
as $\delta\downarrow0$, $\PP^{\theta,\nu_0^\theta}$-a.s..

Next, we apply the comparison principle for \eqref{eq.u.PDE} to obtain:
\begin{align*}
& \inf_{\theta'\in\Theta,\, x\in [0,1)^q} U^{{\theta'},1}_{f^{\theta'}_{0},g^{\theta'}_{0}}[1/E^{\theta'}_{0,1}](0,x) \\
&\geq e^{-\sup_{s\in[0,1],\,x\in [0,1)^q}|f^{\theta'}_{0}(s,x)|}\inf_{x\in [0,1)^q}1/E^{\theta'}_{0,1}(x)\geq e^{-C_2(1+\sup_{s\in[0,1]}|Y_s|)},\\
& \sup_{\theta'\in\Theta,\, x\in [0,1)^q} U^{{\theta'},1}_{f^{\theta'}_{0},g^{\theta'}_{0}}[1/E^{\theta'}_{0,1}](0,x) \leq e^{C_2(1+\sup_{s\in[0,1]}|Y_s|)}.
\end{align*}
Collecting the above and applying the dominated convergence theorem, we deduce
\begin{align*}
&\lim_{\delta\downarrow0}\sup_{T\geq1}\frac{1}{T}\EE^{\theta,\nu_0^\theta} \sup_{\bar d(\theta',\theta'')\leq\delta} |\log L^{\theta',\nu}_{\lfloor T\rfloor-1}/L^{\theta'',\nu}_{\lfloor T\rfloor-1}| = 0.
\end{align*}
Similarly, we obtain
\begin{align*}
&\lim_{\delta\downarrow0}\sup_{T\geq1}\frac{1}{T}\EE^{\theta,\nu_0^\theta} \sup_{\bar d(\theta',\theta'')\leq\delta} |\log L^{\theta',\nu}_{\lfloor T\rfloor-1,T-\lfloor T\rfloor+1}/L^{\theta'',\nu}_{\lfloor T\rfloor-1,T-\lfloor T\rfloor+1}| = 0,
\end{align*}
which yields
\begin{align*}
&\lim_{\delta\downarrow0}\sup_{T\geq1}\frac{1}{T}\EE^{\theta,\nu_0^\theta} \sup_{\bar d(\theta',\theta'')\leq\delta} |\log L^{\theta',\nu}_{T} - \log L^{\theta'',\nu}_{T}| = 0.
\end{align*}

\smallskip

Next, we notice that \eqref{eq.uniform.robust.filter} implies
\begin{align*}
&\lim_{\delta\downarrow0}\sup_{T\geq1}\frac{1}{T}\EE^{\theta,\nu_0^\theta} \sup_{\bar d(\theta',\theta'')\leq\delta}
\left| \int_0^T\left((\hat\pi^{\theta'',\nu}_t[h^{\theta''}]-\hat\pi^{\theta',\nu}_t[h^{\theta'}])^\top \hat\pi^{\theta,\nu_0^\theta}_t[h^{\theta}] \right.\right.\\
&\phantom{?????????????????????}\left.\left.
+ \frac{1}{2}(|\hat\pi^{\theta',\nu}_t[h^{\theta'}]|^2-|\hat\pi^{\theta'',\nu}_t[h^{\theta''}]|^2)\right) dt \right| = 0.
\end{align*}
Then, \eqref{eq.rev.step3.eq1} yields
\begin{align*}
&\lim_{\delta\downarrow0}\sup_{T\geq1}\EE^{\theta,\nu_0^\theta} \sup_{\bar d(\theta',\theta'')\leq\delta}
\left| \frac{1}{T}\int_0^T(\hat\pi^{\theta',\nu}_t[h^{\theta'}]-\hat\pi^{\theta'',\nu}_t[h^{\theta'}])^\top dW_t \right| = 0.
\end{align*}
The above states that, for any $\epsilon>0$, there exists $\delta>0$ such that
\begin{align*}
&\sup_{T\geq1}\EE^{\theta,\nu_0^\theta} \sup_{\bar d(\theta',\theta'')\leq\delta}
\left| \frac{1}{T}\int_0^T(\hat\pi^{\theta',\nu}_t[h^{\theta'}]-\hat\pi^{\theta'',\nu}_t[h^{\theta''}])^\top dW_t \right| < \epsilon.
\end{align*}
Finally, we use the compactness of $\Theta$ to deduce the existence of a $\delta$-net $\{\theta_i\}_{i=1}^n$. Then, to conclude the proof, we notice that there exists $T_0\geq1$ such that, for all $T\geq T_0$,
\begin{align*}
&\sum_{i=1}^n\EE^{\theta,\nu_0^\theta}
\left| \frac{1}{T}\int_0^T \hat\pi^{\theta_i,\nu}_t[h^{\theta_i}]^\top dW_t \right| < \epsilon,
\end{align*}
which yields
\begin{align*}
&\EE^{\theta,\nu_0^\theta} \sup_{\theta'\in\Theta}
\left| \frac{1}{T}\int_0^T \hat\pi^{\theta',\nu}_t[h^{\theta'}]^\top dW_t \right| 
\leq \sum_{i=1}^n\EE^{\theta,\nu_0^\theta}
\left| \frac{1}{T}\int_0^T \hat\pi^{\theta_i,\nu}_t[h^{\theta_i}]^\top dW_t \right|\\
&+ \EE^{\theta,\nu_0^\theta} \sup_{\bar d(\theta',\theta'')\leq\delta}
\left| \frac{1}{T}\int_0^T(\hat\pi^{\theta',\nu}_t[h^{\theta'}]-\hat\pi^{\theta'',\nu}_t[h^{\theta'}])^\top dW_t \right| < 2\epsilon.
\end{align*}
\end{proof}

\medskip

The following proposition is deduced from Propositions \ref{le:uniform.robust.filter}, \ref{prop:strong.ergodicity}, \ref{le:uniform.robust.filter.2}, and from Corollaries \ref{cor:uniform.cont.Z}, \ref{cor:Lambda.bdd.away.from.zero}.

\begin{proposition}\label{prop:last.step}
For any $\theta\in\Theta$, $\nu\in\mathcal{P}([0,1)^q)$ and $\epsilon>0$, there exists $\delta>0$ s.t.
\begin{align*}
&\lim_{T\rightarrow\infty} \PP^{\theta,\nu^\theta_0} \left( \inf_{\theta':\, \bar d(\theta',\overline{\Theta}(\theta))\geq\epsilon}\left[\frac{1}{2T}\int_0^T \left|\hat\pi^{\theta,\nu}_t[h^{\theta}] - \hat\pi^{\theta',\nu}_t[h^{\theta'}]\right|^2\,dt\right.\right.\\
&\phantom{???????????????????}\left.\left. - \frac{1}{T} \int_0^T \left(\hat\pi^{\theta',\nu}_t[h^{\theta'}]-\hat\pi^{\theta,\nu}_t[h^{\theta}]\right)^\top d\tilde W^{\theta,\nu^\theta_0}_t\right]\leq \delta\right) = 0.
\end{align*}
\end{proposition}
\begin{proof}
We fix arbitrary $\theta\in\Theta$, $\nu\in\mathcal{P}([0,1)^q)$, $\epsilon,\varepsilon>0$, denote 
\begin{align*}
&\iota:=\inf_{\theta':\, \bar d(\theta',\overline{\Theta}(\theta))\geq\epsilon}\Lambda(\theta,\theta')>0,
\end{align*}
with the latter inequality being due to Corollary \ref{cor:Lambda.bdd.away.from.zero}, and choose $\delta:={\iota/12}$. 
Next, we introduce $Z^{\theta'}_t:=\frac{1}{t}\int_0^t\left|\hat\pi^{\theta,\nu}_t[h^{\theta}]-\hat\pi^{\theta',\nu}_t[h^{\theta'}]\right|^2\,ds$ and deduce from Corollary \ref{cor:uniform.cont.Z} the existence of $\delta'>0$ s.t.
\begin{align*}
&\sup_{t\geq0}\PP^{\theta,\nu^\theta_0}\left(\sup_{ \bar d(\theta',\theta'')\leq\delta'} |Z^{\theta'}_t - Z^{\theta''}_t|\geq\iota/6\right)<\varepsilon.
\end{align*}
Since $\Theta$ is compact, there exists a finite $\delta'$-net $\{\theta_i\}_{i=1}^m$.
Then,
\begin{align*}
&\PP^{\theta,\nu^\theta_0}\left( \inf_{\theta':\, \bar d(\theta',\overline{\Theta}(\theta))\geq\epsilon}\left[\frac{1}{2T}\int_0^T \left|\hat\pi^{\theta,\nu}_t[h^{\theta}] - \hat\pi^{\theta',\nu}_t[h^{\theta'}]\right|^2\,dt \right.\right.\\
&\phantom{?????????????????????????}\left.\left.- \frac{1}{T}\int_0^T \left(\hat\pi^{\theta',\nu}_t[h^{\theta'}]-\hat\pi^{\theta,\nu}_t[h^{\theta}]\right)^\top d\tilde W^{\theta,\nu^\theta_0}_t\right]\leq \delta\right)\\
&\leq \PP^{\theta,\nu^\theta_0}\left(\inf_{\theta':\, \bar d(\theta',\overline{\Theta}(\theta))\geq\epsilon}Z^{\theta'}_T \leq 4\delta\right)\\
&\phantom{?????????????????????????} + \PP^{\theta,\nu^\theta_0}\left(\sup_{\theta'}
\frac{1}{T}\left|\int_0^T(\hat\pi^{\theta',\nu}_t[h^{\theta'}]-\hat\pi^{\theta,\nu}_t[h^{\theta}])^\top\,d\tilde W^{\theta,\nu^\theta_0}_t\right|
\geq\delta\right)\\
&= \PP^{\theta,\nu^\theta_0}\left(\inf_{\theta':\, \bar d(\theta',\overline{\Theta}(\theta))\geq\epsilon}Z^{\theta'}_T \leq \iota/3\right)\\
&\phantom{??????????????????????} + \PP^{\theta,\nu^\theta_0}\left(\sup_{\theta'}
\frac{1}{T}\left|\int_0^T(\hat\pi^{\theta',\nu}_t[h^{\theta'}]-\hat\pi^{\theta,\nu}_t[h^{\theta}])^\top\,d\tilde W^{\theta,\nu^\theta_0}_t\right|
\geq\iota/12\right)\\
&\leq \PP^{\theta,\nu^\theta_0}\left(\inf_{i}Z^{\theta_i}_T \leq 2\iota/3\right)
+ \PP^{\theta,\nu^\theta_0}\left(\sup_{ \bar d(\theta',\theta'')\leq\delta'} |Z^{\theta'}_T - Z^{\theta''}_T|\geq\iota/3\right)\\
&\phantom{??????????????????????} 
+ \PP^{\theta,\nu^\theta_0}\left(\sup_{\theta'}
\frac{1}{T}\left|\int_0^T(\hat\pi^{\theta',\nu}_t[h^{\theta'}]-\hat\pi^{\theta,\nu}_t[h^{\theta}])^\top\,d\tilde W^{\theta,\nu^\theta_0}_t\right|
\geq\iota/12\right)\\
&\leq \sum_{i=1}^m \PP^{\theta,\nu^\theta_0}\left(\left|Z^{\theta_i}_T - \Lambda(\theta,\theta_i)\right|\geq{\iota/3}\right)
+ \PP^{\theta,\nu^\theta_0}\left(\sup_{ \bar d(\theta',\theta'')\leq\delta'} |Z^{\theta'}_T - Z^{\theta''}_T|\geq{\iota/3}\right)\\
&\phantom{?????????????????????????} + \PP^{\theta,\nu^\theta_0}\left(\sup_{\theta'}
\frac{1}{T}\left|\int_0^T\hat\pi^{\theta',\nu}_t[h^{\theta'}]^{\top}\,d\tilde W^{\theta,\nu^\theta_0}_t\right|
\geq{\iota/24}\right).
\end{align*}
Using Propositions \ref{le:uniform.robust.filter}, \ref{prop:strong.ergodicity} and \ref{le:uniform.robust.filter.2}, we conclude that the right hand side of the above is less than $3\varepsilon$ for all large enough $T>0$. This concludes the proof of the proposition.
\end{proof}

\medskip

Finally, we conclude the proof of Theorem \ref{thm:main}.
First, using the periodicity of $(b,\sigma,h)$ and the a.s. continuity of the filter w.r.t. the parameter (which follows, e.g., from Lemma \ref{le:v.reg}), we conclude that $L^{\theta',\nu}_t=L^{\theta',\hat\nu}_t$ for all $\theta'\in\Theta$, $\PP^{\theta,\nu^\theta_0}$-a.s., where $\hat\nu$ is the projection of $\nu$ defined in \eqref{eq.meas.proj.def}. Hence, it suffices to prove Theorem \ref{thm:main} for $\nu\in\mathcal{P}([0,1)^q)$.
To this end, we fix an arbitrary $\nu\in\mathcal{P}([0,1)^q)$ and consider a MLE $\hat\theta^\nu$.
Then, for any $\theta\in\Theta$ and $\varepsilon>0$, we recall \eqref{eq.Y.projected}--\eqref{eq.L.def} to obtain
\begin{align}
& \PP^{\theta,\nu_0^\theta}\left( \bar d(\hat\theta^\nu_T,\overline{\Theta}(\theta))\geq\varepsilon\right)
{= \PP^{\theta,\nu_0^\theta}\left( \sup_{\theta':\, \bar d(\theta',\overline{\Theta}(\theta))\geq\varepsilon}L^{\theta',\nu}_T \geq L^{\theta,\nu}_T \right)}\nonumber\\
&{= \PP^{\theta,\nu_0^\theta}\left( \sup_{\theta':\, \bar d(\theta',\overline{\Theta}(\theta))\geq\varepsilon}\frac{1}{T}\log L^{\theta',\nu}_T \geq \frac{1}{T}\log L^{\theta,\nu}_T \right)}
= \PP^{\theta,\nu_0^\theta}\left( \sup_{\theta':\, \bar d(\theta',\overline{\Theta}(\theta))\geq\varepsilon}\frac{1}{T}\log \frac{L^{\theta',\nu}_T}{L^{\theta,\nu}_T}\geq 0\right)\nonumber\\
&=\PP^{\theta,\nu_0^\theta}\left( \inf_{\theta':\, \bar d(\theta',\overline{\Theta}(\theta))\geq\varepsilon}\left[\frac{1}{2T}\int_0^T \left(|\pi^{\theta',\nu}_t[h^{\theta'}]|^2 - |\pi^{\theta,\nu}_t[h^{\theta}]|^2\right)\,dt \right.\right.\nonumber\\
&\phantom{?????????????????????????????????}\left.\left.- \frac{1}{T}\int_0^T \pi^{\theta,\nu^\theta_0}_t[h^\theta]^\top\,\left(\pi^{\theta',\nu}_t[h^{\theta'}]-\pi^{\theta,\nu}_t[h^{\theta}]\right)\,dt\right.\right.\nonumber\\
&\phantom{?????????????????????????????????}\left.\left.- \frac{1}{T}\int_0^T \left(\pi^{\theta',\nu}_t[h^{\theta'}]-\pi^{\theta,\nu}_t[h^{\theta}]\right)^\top\, d\tilde{W}^{\theta,\nu_0^\theta}_t\right]\leq 0\right)\label{eq.finalEq}\\
&=\PP^{\theta,\nu_0^\theta}\left( \inf_{\theta':\, \bar d(\theta',\overline{\Theta}(\theta))\geq\varepsilon}\left[\frac{1}{2T}\int_0^T \left|\hat\pi^{\theta',\nu}_t[h^{\theta'}] - \hat\pi^{\theta,\nu}_t[h^{\theta}]\right|^2\,dt \right.\right.\nonumber\\
&\phantom{??????????????}\left.\left.+ \frac{1}{T}\int_0^T \left(\hat\pi^{\theta,\nu}_t[h^\theta]-\hat\pi^{\theta,\nu^\theta_0}_t[h^\theta]\right)^\top\,\left(\hat\pi^{\theta',\nu}_t[h^{\theta'}]-\hat\pi^{\theta,\nu}_t[h^{\theta}]\right)\,dt\right.\right.\nonumber\\
&\phantom{?????????????????????????????????}
\left.\left.- \frac{1}{T}\int_0^T \left(\hat\pi^{\theta',\nu}_t[h^{\theta'}]-\hat\pi^{\theta,\nu}_t[h^{\theta}]\right)^\top d\tilde{W}^{\theta,\nu_0^\theta}_t\right]\leq 0\right),\nonumber
\end{align}
where the last equality is deduced using
\begin{align*}
& |\pi^{\theta',\nu}_t[h^{\theta'}]|^2 - |\pi^{\theta,\nu}_t[h^{\theta}]|^2
=\left|\hat\pi^{\theta',\nu}_t[h^{\theta'}] - \hat\pi^{\theta,\nu}_t[h^{\theta}]\right|^2
+ 2\hat\pi^{\theta,\nu}_t[h^\theta]^\top\,\left(\hat\pi^{\theta',\nu}_t[h^{\theta'}]-\hat\pi^{\theta,\nu}_t[h^{\theta}]\right).
\end{align*}

Noticing that
\begin{align*}
&\left| \frac{1}{T}\int_0^T \left(\hat\pi^{\theta,\nu}_t[h^\theta]-\hat\pi^{\theta,\nu^\theta_0}_t[h^\theta]\right)^\top\,\left(\hat\pi^{\theta',\nu}_t[h^{\theta'}]-\hat\pi^{\theta,\nu}_t[h^{\theta}]\right)\,dt\right|\\
&\phantom{???????????????????????????????????????}\leq C\,\frac{1}{T}\int_0^T \left|\hat\pi^{\theta,\nu}_t[h^\theta]-\hat\pi^{\theta,\nu^\theta_0}_t[h^\theta]\right|\,dt,
\end{align*}
we apply Lemma \ref{cor:H.contracts.on.pi} to deduce that the right hand side of the above converges to zero, as $T\rightarrow\infty$, $\PP^{\theta,\nu^\theta_0}$-a.s..
Combining the latter observation with Proposition \ref{prop:last.step}, we conclude that the right hand side of \eqref{eq.finalEq} vanishes, as $T\rightarrow\infty$, which yields the desired consistency of MLE $\hat\theta^\nu$ and completes the proof of Theorem \ref{thm:main}.

\section{Appendix}\label{se:appendix}

\subsection{Proof of Lemma \ref{ass:regular.stationary.density}}
\label{subse:regular.stationary.density}

We fix arbitrary $\theta\in\Theta$ and $\nu\in\mathcal{P}([0,1)^q)$.
The existence of an invariant measure for $X^e$ follows from the standard arguments. Namely, consider the following family of probability measures on $[0,1)^q$:
\begin{align*}
&\mu_T(A) := \frac{1}{T} \int_0^T \PP^{\theta,\nu}(X^e_t\in A)\,dt,\quad\forall\,A\in \mathcal{B}([0,1)^q).
\end{align*}
Since the measures $\mu_T$ are defined on a compact space, with respect to the torus metric, the above family is tight, and the Prokhorov's theorem yields the existence of a sequence $\{\mu_{T_n}\}$ converging weakly to a limit $\nu_0\in \mathcal{P}([0,1)^q)$.
Using the Markov property of $X$ (which follows from Assumption \ref{holder}), as well as the periodicity of $(b,\sigma)$ and the uniqueness in law of the solution to \eqref{eq.Background.dyn.X}, we deduce the Markov property of $X^e$. The latter yields, for any $s\geq0$:
\begin{align*}
&\PP^{\theta,\mu_T}(X^e_s\in dx) = \frac{1}{T} \int_0^T \PP^{\theta,\nu}(X^e_{t+s}\in dx)\,dt
= \frac{1}{T} \int_s^{T+s} \PP^{\theta,\nu}(X^e_{t}\in dx)\,dt.
\end{align*}
It is easy to see that the right hand side of the above converges weakly to $\nu_0$ along $\{T_n\}$.
In addition, since $(b,\sigma)$ are Lipschitz, there exists a unique strong solution to \eqref{eq.Background.dyn.X}, and it is easy to verify that this solution is continuous with respect to the initial condition. Thus, we can pass to the limit in the left hand side of the above display, along $\{T_n\}$, to obtain
\begin{align*}
&\PP^{\theta,\nu_0}(X^e_s\in dx) = \nu_0(dx),\quad\forall\,s>0,
\end{align*}
which means that $\nu_0$ is an invariant measure for $X^e$. Using, once more, the uniqueness in law of the solution to \eqref{eq.Background.dyn.X} and the periodicity of $(b,\sigma)$, we deduce that the distribution of $(X^e_{s+\cdot})_{\RR_+}$ under $\PP^{\theta,\nu_0}$ does not depend on $s\geq0$.

\smallskip

To show that $\nu_0$ has a continuously differentiable density that is bounded away form zero, we consider any $T>0$, $\phi\in C(\RR^q)\cap\mathrm{P}_1$ and the associated $u:=U^{\theta,T}_{0,0}[\phi]$, and apply the Feynman-Kac formula along with the representation \eqref{eq.appendix.Gauss.est.eq3}, to deduce that
\begin{align*}
&\int_{[0,1)^q}\phi(y)\,\nu_0[dy]=\EE^{\theta,\nu_0}\phi(X_T) = \int_{\RR^q}\int_{[0,1)^q}Q(0,x;T,y)\,\nu_0(dx)\,\phi(y)\,dy\\
&= \int_{[0,1)^q}\int_{[0,1)^q} \sum_{e\in\ZZ^q} Q(0,x;T,y+e)\,\nu_0(dx)\,\phi(y)\,dy.
\end{align*}
The above yields:
\begin{align*}
&\nu_0[dy]= \int_{[0,1)^q} \sum_{e\in\ZZ^q} Q(0,x;T,y+e)\,\nu_0(dx)\,dy,
\end{align*}
which means that $\nu_0$ has a density $\psi_0:=\int_{[0,1)^q} \sum_{e\in\ZZ^q} Q(0,x;T,\cdot+e)\,\nu_0(dx)$. The lower bound and the continuous differentiability of this density follow from \eqref{eq.appendix.Gauss.est1}--\eqref{eq.appendix.Gauss.est2}.

\subsection{Proof of Lemma \ref{ass:mixing}}
\label{subse:ass:mixing}

Let us fix $T$ and $f$, $g$ as in the statement of the lemma.
Applying \cite[Theorem 7 in Chapter 9.6]{Friedman1969} and \cite[Theorem 1]{Aronson}, we deduce the existence of a function $Q(t,x;s,y)$, defined for all $x,y\in\RR^q$ and $0\leq t<s\leq T$, such that, for each $(s,y)$, the function $Q(\cdot,\cdot;s,y)$ satisfies the PDE in \eqref{eq.u.PDE} in the domain $(t,x)\in[0,s)\times\RR^q$ and for all $i,j=1,\ldots,q,\,\,0\leq k+m\leq1,\,\,0\leq r\leq 2$, we have
\begin{align}
& \left|\partial^{r+k}_{x_i}\partial^m_{y_j} Q(t,x;s,y) \right| \leq \frac{C_1}{(s-t)^{(q+r+m+k)/2}} e^{-\frac{(x-y)^2}{C_2(s-t)}},\label{eq.appendix.Gauss.est1}\\
&Q(t,x;s,y) \geq \frac{C_3}{(s-t)^{q/2}} e^{-\frac{(x-y)^2}{C_2(s-t)}},\label{eq.appendix.Gauss.est2}\\
& \lim_{t\uparrow s}\int_{\RR^q} Q(t,x;s,y)\,\phi(x)\,dx = \phi(y)\quad\forall\,\phi\in C_b(\RR^q),\nonumber
\end{align}
where $\{C_i,\,1/C_i\}$ are locally bounded strictly positive functions of $(\|f\|_{C^1_{0,T}},\|g\|_{C^1_{0,T}})$.
Then, it is easy to check that the function
\begin{align}
&u(t,x):=\int_{\RR^q} Q(t,x;T,y)\,\phi(y)\,dy\label{eq.appendix.Gauss.est.eq3}
\end{align}
is a bounded classical solution to \eqref{eq.u.PDE}.
To show the uniqueness of such a solution, we derive the PDE for $v(t,x):=u(t,x)\,e^{\lambda(T-t)-\sqrt{1+|x|^2}}$ and apply the maximum principle to it (choosing $\lambda$ of appropriate sign, to handle the maximum and the minimum points of the difference of two candidate solutions). The uniqueness of the solution and the periodicity of $(f,g,\phi)$ yield the periodicity of $u(t,\cdot)$.
In addition, the explicit construction of $Q$ via the parametrix method, described, e.g., in \cite{Friedman1969}, yields the desired measurability of $u$ as a function of $(f,g,\phi)$.

\smallskip

It only remains to show how \eqref{eq.appendix.Gauss.est1}--\eqref{eq.appendix.Gauss.est2} yield \eqref{eq.u.positive}--\eqref{eq:lowerassum}, with $\epsilon_0(t,z_1,z_2)=\sqrt{t}\,\epsilon_2(z_1,z_2)$, where $\epsilon_2:\RR^2_+\rightarrow(0,\infty)$ is a non-increasing (in each variable) strictly positive measurable function. Since the proofs of all the inequalities in \eqref{eq.u.positive}--\eqref{eq:lowerassum} are very similar, we only present the proof of \eqref{eq:derder} and of a part of \eqref{eq:dersup}:
\begin{align*}
&\int_{[0,1)^q} |\partial_{x_i}U^{\theta,T}_{f,g}[\phi](T-T_0,x)|\,dx=\int_{[0,1)^q}\left|\int_{\RR^q} \partial_{x_i}Q(T-T_0,x;T,y)\,\phi(y)\,dy \right|\,dx\\
&=\int_{[0,1)^q}\left| \int_{[0,1)^q} \sum_{e\in\ZZ^q} \partial_{x_i}Q(T-T_0,x;T,y+e)\,\phi(y)\,dy \right|\,dx\\
&\leq\int_{[0,1)^q} \int_{[0,1)^q} \sum_{e\in\ZZ^q}\left| \partial_{x_i}Q(T-T_0,x;T,y+e)\right|\,dx |\phi(y)|\,dy\\
&\leq\int_{[0,1)^q} \int_{[0,1)^q} \sum_{e\in\ZZ^q}\frac{C_4}{T_0^{\frac{q+1}{2}}} e^{-\frac{(x-y+e)^2}{C_2\,T_0}}\,dx |\phi(y)|\,dy\leq \int_{[0,1)^q} \int_{\RR^q} \frac{C_5}{T_0^{\frac{q+1}{2}}} e^{-\frac{(x-y)^2}{C_2\,T_0}}\,dx |\phi(y)|\,dy\\
&\leq \int_{[0,1)^q} \frac{C_6}{\sqrt{T_0}}\,|\phi(y)|\,dy
= \frac{C_6}{\sqrt{T_0}}\,\|\phi\|_{L^1([0,1)^q)}\quad\mbox{and, similarly,}\\
& |U^{\theta,T}_{f,g}[\partial_{x_i}\phi](t,x)|
= \left| \int_{\RR^q} Q(t,x;T,y)\,\partial_{y_i}\phi(y)\,dy \right| = \left| \int_{\RR^q} \partial_{y_i}Q(t,x;T,y)\,\phi(y)\,dy \right|\\
&= \left| \int_{[0,1)^q} \sum_{e\in\ZZ^q} \partial_{y_i}Q(t,x;T,y+e)\,\phi(y)\,dy \right|
\leq \int_{[0,1)^q} \frac{C_4}{(T-t)^{\frac{q+1}{2}}} e^{-\frac{(x-y)^2}{C_2\,(T-t)}}\,|\phi(y)|\,dy\\
&=\int_{[0,1)^q} \frac{C_4}{\sqrt{T-t}} e^{-\frac{z^2}{C_2}}\,|\phi(x-\sqrt{T-t}z)|\,dz
\leq \frac{C_6}{\sqrt{T-t}}\,\|\phi\|_C,
\end{align*}
with $C_5$, $C_6$ depending on $\|f\|_{C^1_{0,T}}$ and $\|g\|_{C^1_{0,T}}$.

\subsection{Proof of Lemma \ref{le:v.reg}}
\label{subse:le:v.reg}

To ease the notation, we denote $O_Y:=\mathrm{Osc}_{t,t+T}(Y)\vee \mathrm{Osc}_{t,t+T}(\tilde Y)$.
First, we fix an arbitrary $t\geq0$ and notice that there exists a function $M:(0,\infty)^2\mapsto [1,\infty)$ that is increasing in both variables and such that, for all $T>0$,
\begin{align}
&\|E^\theta_{t,\cdot}\|_{C^{1}_{0,T}}+\|1/E^\theta_{t,\cdot}\|_{C_{0,T}}+\|f^\theta_t\|_{C^{1}_{0,T}}+\|g^\theta_t\|_{C^{1}_{0,T}}\leq M(T, O_Y),\notag\\
&\label{eq:btheta}\|g^{\theta}_t-g^{\theta'}_t\|_{C_{0,T}}+\|f^{\theta}_t-f^{\theta'}_t\|_{C_{0,T}}\leq M(T,O_Y)\\
&\phantom{???????????????}\times\left(\|(a^{\theta}-a^{\theta'}\|_C+\|(b^{\theta}-b^{\theta'}\|_C+\|(h^{\theta}-h^{\theta'}\|_{C^2}+\|Y-\tilde Y\|_{C(t,t+T)}\right).\notag
\end{align}
Next, we recall the functions $\bar C,\epsilon_0,\epsilon_1$ defined in Lemma \ref{ass:mixing} and, using their monotonicity and the integrability of $1/\epsilon_0$, deduce that the function 
\begin{align*}
\bar L(T,x)&:=\int_{0}^{T/2}\frac{M(T,x)(1+M(T,x))}{ \epsilon_0\left(s,M(T,x),M(T,x)\right)\epsilon_1\left(T/2,M(T,x),M(T,x)\right)}ds\\
      &+T\frac{\bar C\left(T,M(T,x),M(T,x)\right)M(T,x)}{\epsilon_1\left(T/2,M(T,x),M(T,x)\right)}(1+M(T,x))\\
      &+\int_{T/2}^T \frac{M(T,x) (1+M(T,x))}{\epsilon_0(T-s,M(T,x),M(T,x))\epsilon_1(T/2,M(T,x),M(T,x))}ds>0
\end{align*}
is finite for each $T,x>0$ and increasing in each variable.
We also define the non-increasing (in each variable) function 
\begin{align*}
\epsilon(T,x):=\frac{1}{\bar L(T,x)}\wedge\frac{\epsilon_1(T,M(T,x),M(T,x))}{M(T,x)}.
\end{align*}

\smallskip
    
Next, we fix $(t,T,\theta,\theta')\in\RR_+\times (0,\infty)\times \Theta^2$ and notice that, for $\mathbb{W}_Y$-a.e. $Y$, we have 
\begin{align*}
\|1/E^\theta_{t,\cdot}\|_{C_{0,T}}+\|f^\theta_t\|_{C^{1}_{0,T}}+\|g^\theta_t\|_{C^{1}_{0,T}}\leq M(T, O_Y)<\infty.
\end{align*}
Then, thanks to \eqref{eq.u.positive}, the linear mapping
\begin{align*}
C^1(\RR^q)\cap\mathrm{P}_1\ni\psi \mapsto U^{\theta,T}_{f^\theta_t,g^\theta_t}[\psi/E^\theta_{t,T}](0,\cdot)\in C^1(\RR^q)\cap\mathrm{P}_1
\end{align*}
satisfies 
\begin{align*}
|U^{\theta,T}_{f^\theta_t,g^\theta_t}[\psi/E^\theta_{t,T}](0,x)|\leq & \frac{\|\psi/E^\theta_{t,T}\|_{L^1([0,1)^q)}}{\epsilon_1\left(T,M(T, O_Y),M(T, O_Y)\right)}\leq \frac{\|\psi\|_{L^1([0,1)^q)}}{\epsilon(T,O_Y)}.
\end{align*}
Thus, the aforementioned linear mapping can be extended to a continuous linear operator in $L^1([0,1)^q)$, $\mathbb{W}_Y$-a.s..
Thanks to Reisz representation theorem, there exists $(Y,x,z)\mapsto\hat K^\theta_{t,T}(x,z)$, which inherits its measurability from $(\phi,f,g)\mapsto U^{\theta,T}_{f,g}[\phi]$, such that the following holds $\mathbb{W}_Y$-a.s.:
\begin{align}
&\sup_{x,z\in[0,1)^q}|\hat K^\theta_{t,T}(x,z)|\leq \frac{1}{\epsilon\left(T,O_Y\right)},
\quad U^{\theta,T}_{f^\theta_t,g^\theta_t}[\psi/E^\theta_{t,T}](0,x)=\int_{[0,1)^q} \hat K^\theta_{t,T}(x,z) \psi(z)dz,
\label{eq.rev.appendix.eq1}
\end{align}
for all $\psi\in C(\RR^q)\cap\mathrm{P}_1$ and $x\in[0,1)^q$.
In addition, for $0\leq \psi\in C(\RR^q)\cap\mathrm{P}_1$, the property \eqref{eq:lowerassum} implies
 \begin{align*}
 \int_{[0,1)^q} \hat K^\theta_{t,T}(x,z) \psi(x)dx&\geq \epsilon_1(T,M(T, O_Y), M(T, O_Y))\int_{[0,1)^q} \psi(x)/E^\theta_{t,T}(x)\,dx\\
 &\geq \epsilon\left(T,O_Y\right)\int_{[0,1)^q} \psi(x)\,dx,
 \end{align*}
which leads to the lower bound in \eqref{eq:boundK}. Similarly, \eqref{eq.u.positive} yields the upper bound in \eqref{eq:boundK}.

\smallskip

To show \eqref{eq:rep2}, we consider an arbitrary $\psi \in C(\RR^q)\cap \mathrm{P}_1$ and apply \eqref{eq:rep1} along with the periodicity of the associated functions, to obtain 
\begin{align*}
&\int_{[0,1)^q}\psi(x)\hat p^{\theta,\nu}_{t+T}(x)\,dx=\int_{\RR^q}\psi(x)\pi^{\theta,\nu}_{t+T}(dx)=\frac{\int_{\RR^q} U^{\theta,T}_{f^\theta_t,g^\theta_t}[\psi/E^\theta_{t,T}](0,x)\pi^{\theta,\nu}_{t}(dx)}{\int_{\RR^q} U^{\theta,T}_{f^\theta_t,g^\theta_t}[1/E^\theta_{t,T}](0,x)\pi^{\theta,\nu}_{t}(dx)}\\
&=\frac{\int_{[0,1)^q} U^{\theta,T}_{f^\theta_t,g^\theta_t}[\psi/E^\theta_{t,T}](0,x) \hat\pi^{\theta,\nu}_{t}(dx)}{\int_{[0,1)^q} U^{\theta,T}_{f^\theta_t,g^\theta_t}[1/E^\theta_{t,T}](0,x) \hat\pi^{\theta,\nu}_{t}(dx)}
=\frac{\int_{[0,1)^q} \int_{[0,1)^q} \hat K^\theta_{t,T}(x,z) \psi(z)\,dz\,\hat\pi^{\theta,\nu}_{t}(dx)}{\int_{[0,1)^q}\int_{[0,1)^q} \hat K^\theta_{t,T}(x,z)\,dz\,\hat\pi^{\theta,\nu}_{t}(dx)},
\end{align*}
which is the desired \eqref{eq:rep2}.

\smallskip

It remains to prove \eqref{eq:contdv}.
To this end, we fix $(t,T,\theta,\theta')\in\RR_+\times (0,\infty)\times \Theta^2$ and $\psi\in C(\RR^q)\cap P_1$, with $\|\psi\|_{L^1([0,1)^q)}\leq 1$, and define
\begin{align*}
\Delta^{\theta,\theta'}(s,x):=U^{\theta,T}_{f^\theta_t,g^\theta_t}[\psi/E^\theta_{t,T}](s,x)-U^{\theta',T}_{f^{\theta'},g^{\theta'}}[\psi/E^{\theta'}_{t,T}](s,x),\quad (s,x)\in [0,T)\times \RR^q,
\end{align*}
so that 
\begin{align*}
\Delta^{\theta,\theta'}(0,x)=\int_{[0,1)^q} \psi(z) \left(\hat{K}^{\theta}_{t,T}(Y_{[t,t+T]}-Y_t;x,z)- \hat{K}^{\theta'}_{t,T}(\tilde Y_{[t,t+T]}-\tilde Y_t;x,z)\right)dz
\end{align*}
and $\Delta^{\theta,\theta'}$ solves the parabolic equation
\begin{align}
&\pa_s \Delta^{\theta,\theta'} + \cL^\theta\Delta^{\theta,\theta'} + g^\theta_t\nabla^\top\Delta^{\theta,\theta'} + \Delta^{\theta,\theta'} f^\theta_t =-(\cL^\theta-\cL^{\theta'}) U^{\theta',T}_{f^{\theta'}_t,g^{\theta'}_t}[\psi/E^{\theta'}_{t,T}]\label{eq:difftheta}\\
&-(g^{\theta}_t-g^{\theta'}_t)\nabla^\top U^{\theta',T}_{f^{\theta'},g^{\theta'}}[\psi/E^{\theta'}_{t,T}] - (f^\theta_t-f^{\theta'}_t) U^{\theta',T}_{f^{\theta'}_t,g^{\theta'}_t}[\psi/E^{\theta'}_{t,T}],\notag\\
&\Delta^{\theta,\theta'}(T,x)=\psi/E^{\theta}_{t,T} - \psi/E^{\theta'}_{t,T}.\notag
\end{align}
It is a standard exercise to verify that the bounded classical solution $\Delta^{\theta,\theta'}$ to the above parabolic problem (which is unique due to Lemma \ref{ass:mixing}) is given by $U^{\theta,T}_{f^\theta_t,g^\theta_t}[\psi/E^\theta_{t,T}-\psi/E^{\theta'}_{t,T}]$ plus the $s$-integral of $U^{\theta,s}_{f^\theta_t,g^\theta_t}[\cdot]$ applied to the right hand side of the PDE in \eqref{eq:difftheta} evaluated at time $s$.
In particular, we obtain:
\begin{align*}
\Delta^{\theta,\theta'}(0,x)& = U^{\theta,T}_{f^\theta_t,g^\theta_t}[\psi/E^\theta_{t,T}-\psi/E^{\theta'}_{t,T}](0,x)\\
&+ \int_0^T U^{\theta,s}_{f^{\theta}_t,g^{\theta}_t}\left[(\cL^\theta-\cL^{\theta'}) U^{\theta',T}_{f^{\theta'}_t,g^{\theta'}_t}[\psi/E^{\theta'}_{t,T}](s,\cdot)\right](0,x)ds\\
&+\int_0^T U^{\theta,s}_{f^\theta_t,g^\theta_t}\left[(g^{\theta}_t(s,\cdot)-g^{\theta'}_t(s,\cdot))\,\nabla^\top U^{\theta',T}_{f^{\theta'}_t,g^{\theta'}_t}[\psi/E^{\theta'}_{t,T}](s,\cdot)\right](0,x)\,ds\\
&+\int_0^T U^{\theta,s}_{f^{\theta}_t,g^{\theta}_t}\left[(f^{\theta}_t(s,\cdot)-f^{\theta'}_t(s,\cdot))\,U^{\theta',T}_{f^{\theta'}_t,g^{\theta'}_t}[\psi/E^{\theta'}_{t,T}](s,\cdot)\right](0,x)\,ds.
\end{align*}
The product rule leads to  
\begin{align*}
&(\cL^\theta-\cL^{\theta'}) U^{\theta',T}_{f^{\theta'}_t,g^{\theta'}_t}[\psi/E^{\theta'}_{t,T}]
=\frac{1}{2}\sum_{i,j=1}^q \pa_{x_i} \left[(a_{i,j}^{\theta} - a_{i,j}^{\theta'})\pa_{x_j} U^{\theta',T}_{f^{\theta'}_t,g^{\theta'}_t}[\psi/E^{\theta'}_{t,T}]\right]\\
&-\frac{1}{2}\sum_{i,j=1}^q (\pa_{x_i}a_{i,j}^{\theta}-\pa_{x_i}a_{i,j}^{\theta'})\pa_{x_j} U^{\theta',T}_{f^{\theta'}_t,g^{\theta'}_t}[\psi/E^{\theta'}_{t,T}]
+(b^\theta-b^{\theta'})\nabla^\top U^{\theta',T}_{f^{\theta'}_t,g^{\theta'}_t}[\psi/E^{\theta'}_{t,T}],
\end{align*}
so that 
\begin{align*}
    \Delta^{\theta,\theta'}(0,x)&=U^{\theta,T}_{f^\theta_t,g^\theta_t}[\psi/E^\theta_{t,T}-\psi/E^{\theta'}_{t,T}](0,x)\\
& + \sum_{i,j=1}^q\int_0^T U^{\theta,s}_{f^\theta_t,g^\theta_t}\left[ \pa_{x_i} \left(\frac{1}{2}(a_{i,j}^\theta-a_{i,j}^{\theta'})\,\pa_{x_j} U^{\theta',T}_{f^{\theta'}_t,g^{\theta'}_t}[\psi/E^{\theta'}_{t,T}](s,\cdot)\right) \right](0,x)ds\\
& - \sum_{i,j=1}^q\int_0^T {U^{\theta,s}_{f^\theta_t,g^\theta_t}\left[ \frac{1}{2} (\pa_{x_i}a_{i,j}^\theta-\pa_{x_i}a_{i,j}^{\theta'})\,\pa_{x_j} U^{\theta',T}_{f^{\theta'}_t,g^{\theta'}_t}[\psi/E^{\theta'}_{t,T}](s,\cdot) \right](0,x)}ds\\
& + \int_0^T U^{\theta,s}_{f^\theta_t,g^\theta_t}\left[ (b^\theta+g^{\theta}_t(s,\cdot)-b^{\theta'}-g^{\theta'}_t(s,\cdot))\,\nabla^\top U^{\theta',T}_{f^{\theta'}_t,g^{\theta'}_t}[\psi/E^{\theta'}_{t,T}](s,\cdot) \right](0,x)ds\\
    &+\int_0^T U^{\theta,s}_{f^\theta_t,g^\theta_t}\left[ (f^{\theta}_t(s,\cdot)-f^{\theta'}_t(s,\cdot))\,U^{\theta',T}_{f^{\theta'}_t,g^{\theta'}_t}[\psi/E^{\theta'}_{t,T}](s,\cdot) \right](0,x)ds\\
    &=:I^0+\int_0^T\left(\sum_{i,j=1}^q I_{i,j}^1(s)-\sum_{i,j=1}^q I_{i,j}^2(s)+I^3(s)+I^4(s)\right)\,ds.
\end{align*}
The estimate \eqref{eq:boundK} yields
\begin{align*}
|I^0|\leq \frac{\|1/E^\theta_{t,T}-1/E^{\theta'}_{t,T}\|_C \|\psi \|_{L^1([0,1)^q)}}{\epsilon(T,O_Y)}\leq  \frac{\|1/E^\theta_{t,T}-1/E^{\theta'}_{t,T}\|_C }{\epsilon(T,O_Y)}.
\end{align*}

\smallskip

Next, we estimate the remaining terms separately for $s\in (0,T/2)$ and $s\in (T/2,T)$. 
First, we fix an arbitrary $s\in (0,T/2)$ and apply \eqref{eq:dersup} with $T_0=T-s,\,T/2$, to obtain
\begin{align*}
   |I^{1}_{i,j}(s)|+|I^{2}_{i,j}(s)|\leq 2\frac{\|a_{i,j}^\theta-a_{i,j}^{\theta'}\|_{C^1}\,\|\nabla U^{\theta',T}_{f^{\theta'}_t,g^{\theta'}_t}[\psi/E^{\theta'}_{t,T}](s,\cdot)\|_C}{\epsilon_0\left(s,M(T,O_Y),M(T,O_Y)\right)}.
\end{align*}
In addition, using \eqref{eq:supsup} and the monotonicity of $\bar C$, we obtain
\begin{align*}
|I^{3}(s)|&
\leq \bar C\left(T,M(T,O_Y),M(T,O_Y)\right) \|\nabla U^{\theta',T}_{f^{\theta'}_t,g^{\theta'}_t}[\psi/E^{\theta'}_{t,T}](s,\cdot)\|_C\\
&\phantom{????????????????????????????????}\times\left(\|b^\theta-b^{\theta'}\|_C+\|g^{\theta}_t(s,\cdot)-g^{\theta'}_t(s,\cdot)\|_C\right),\\
|I^{4}(s)|&\leq \bar C\left(T,M(T,O_Y),M(T,O_Y)\right) \| U^{\theta',T}_{f^{\theta'}_t,g^{\theta'}_t}[\psi/E^{\theta'}_{t,T}](s,\cdot)\|_C\,
\|f^{\theta}_t(s,\cdot)-f^{\theta'}_t(s,\cdot)\|_C.
\end{align*}
Next, we note that $T-s\geq T/2$, $\|1/E^\theta_{t,\cdot}\|_{C_{0,T}}\leq M(T,O_Y)$ and $\|\psi\|_{L^1([0,1)^q)}\leq 1$.
Hence, we can use \eqref{eq.u.positive} to estimate $\nabla U^{\theta',T}_{f^{\theta'}_t,g^{\theta'}_t}[\psi/E^{\theta'}_{t,T}](s,\cdot)$ and $U^{\theta',T}_{f^{\theta'}_t,g^{\theta'}_t}[\psi/E^{\theta'}_{t,T}](s,\cdot)$, thus obtaining, for $s\in (0,T/2)$:
\begin{align*}
|I^{1}_{i,j}(s)|+|I^{2}_{i,j}(s)|&\leq 2\frac{M(T,O_Y)\, \|a_{i,j}^\theta-a_{i,j}^{\theta'}\|_{C^1}}{ \epsilon_0\left(s,M(T,O_Y),M(T,O_Y)\right)\epsilon_1\left(T/2,M(T,O_Y),M(T,O_Y)\right)},\\
|I^{3}(s)|+ |I^{4}(s)|&
\leq \frac{\bar C\left(T,M(T,O_Y),M(T,O_Y)\right)M(T,O_Y)}{\epsilon_1\left(T/2,M(T,O_Y),M(T,O_Y)\right)}\\
    &\times(\|(b^\theta-b^{\theta'}\|_C+\|g^{\theta}_t(s,\cdot)-g^{\theta'}_t(s,\cdot)\|_C+\|f^{\theta}_t(s,\cdot)-f^{\theta'}_t(s,\cdot)\|_C).
 \end{align*}
 
\smallskip

Next, we fix $s\in [T/2,T]$, apply \eqref{eq.u.positive}, \eqref{eq:derder}, and use the inequality $s\geq T/2$ along with the monotonicity of $\epsilon_1$ in $T$, to obtain 
\begin{align*}
   |I^{1}_{i,j}(s)|+|I^{2}_{i,j}(s)|&\leq \frac{\|a_{i,j}^\theta-a_{i,j}^{\theta'}\|_{C^1}\,\|\nabla U^{\theta',T}_{f^{\theta'}_t,g^{\theta'}_t}[\psi/E^{\theta'}_{t,T}](s,\cdot)\|_{L^1([0,1]^q)}}
   { \epsilon_1\left(T/2,M(T,O_Y),M(T,O_Y)\right)}\\
   &\leq \frac{M(T,O_Y)\,\|a_{i,j}^\theta-a_{i,j}^{\theta'}\|_{C^1}}
   { \epsilon_0\left(T-s,M(T,O_Y),M(T,O_Y)\right) \,\epsilon_1\left(T/2,M(T,O_Y),M(T,O_Y)\right)}
\end{align*}
and 
\begin{align*}
    |I^{3}(s)|&
    \leq  \frac{\|\nabla U^{\theta',T}_{f^{\theta'}_t,g^{\theta'}_t}[\psi/E^{\theta'}_{t,T}](s,\cdot)\|_{L^1([0,1)^q)}\,(\|(b^\theta-b^{\theta'}\|_C + \|g^{\theta}_t(s,\cdot)-g^{\theta'}_t(s,\cdot)\|_C)}{\epsilon_1(T/2,M(T,O_Y),M(T,O_Y))},\\
       |I^{4}(s)|&\leq \frac{\| U^{\theta',T}_{f^{\theta'}_t,g^{\theta'}_t}[\psi/E^{\theta'}_{t,T}](s,\cdot)\|_{L^1([0,1)^q)}\,\|f^{\theta}_t(s,\cdot)-f^{\theta'}_t(s,\cdot)\|_C}{\epsilon_1(T/2,M(T,O_Y),M(T,O_Y))}.
\end{align*}
Thus, thanks to \eqref{eq:derder}, we have 
\begin{align*}
&|I^{3}(s)|+   |I^{4}(s)|\\
&\phantom{???????} \leq  \frac{M(T,O_Y)(\|f^{\theta}_t(s,\cdot)-f^{\theta'}_t(s,\cdot)\|_C+\|(b^\theta-b^{\theta'}\|_C+\|g^{\theta}_t(s,\cdot)-g^{\theta'}_t(s,\cdot)\|_C)}{\epsilon_0(T-s,M(T,O_Y),M(T,O_Y))\epsilon_1(T/2,M(T,O_Y),M(T,O_Y))}.
\end{align*}
We now use \eqref{eq:btheta} to estimate 
\begin{align*}
&\|(b^\theta-b^{\theta'}\|_C+\|g^{\theta}_t(s,\cdot)-g^{\theta'}_t(s,\cdot)\|_C+\|f^{\theta}_t(s,\cdot)-f^{\theta'}_t(s,\cdot)\|_C\\
&\leq (1+M(T,O_Y)))\left(\|(a^\theta-a^{\theta'}\|_C+\|(b^\theta-b^{\theta'}\|_C+\|(h^{\theta}-h^{\theta'}\|_{C^2}+\|Y-\tilde Y\|_{C(t,t+T)}\right),
\end{align*}
so that 
\begin{align*}
 &\int_0^{T/2}\sum_{i,j=1}^q |I^{1}_{i,j}(s)|+|I^{2}_{i,j}(s)|+ |I^{3}(s)|+ |I^{4}(s)| ds\\
  &\leq \left(2\int_{0}^{T/2}\frac{M(T,O_Y)}{ \epsilon_0\left(s,M(T,O_Y),M(T,O_Y)\right)\,\epsilon_1\left(T/2,M(T,O_Y),M(T,O_Y)\right)}ds\right.\\
      &\left.\phantom{????????????????????}
      +T\frac{\bar C\left(T,M(T,O_Y),M(T,O_Y)\right)M(T,O_Y)}{\epsilon_1\left(T/2,M(T,O_Y),M(T,O_Y)\right)}\right)\,(1+M(T,O_Y))\\
    &\times \left(\|a^\theta-a^{\theta'}\|_{C^1}+\|(b^\theta-b^{\theta'}\|_C+\|h^\theta-h^{\theta'}\|_{C^2}+\|Y-\tilde Y\|_{C(t,t+T)}\right)\\
 \end{align*}
and 
\begin{align*}
 &\int_{T/2}^T\sum_{i,j=1}^q |I^{1}_{i,j}(s)|+|I^{2}_{i,j}(s)|+ |I^{3}(s)|+ |I^{4}(s)| ds\\
   &\leq \int_{T/2}^T \frac{M(T,O_Y) (1+M(T,O_Y))}{\epsilon_0(T-s,M(T,O_Y),M(T,O_Y))\epsilon_1(T/2,M(T,O_Y),M(T,O_Y))}ds\\
   &\times \left(\|a^\theta-a^{\theta'}\|_{C^1}+\|(b^\theta-b^{\theta'}\|_C+\|h^\theta-h^{\theta'}\|_{C^2}+\|Y-\tilde Y\|_{C(t,t+T)}\right).
   \end{align*}
Given the definition of $\bar L$ and $\epsilon$, for all $\psi\in P_1\cap C(\RR^q)$ with $\|\psi\|_{L^1([0,1)^q)}\leq 1$ we have
\begin{align*}
&\int_{[0,1)^q} \psi(z) \left(\hat{K}^{\theta}_{t,T}(x,z)- \hat{K}^{\theta'}_{t,T}(x,z)\right)dz=\Delta^{\theta,\theta'}(0,x)\\
&\leq \frac{\|a^\theta-a^{\theta'}\|_{C^1}+\|(b^\theta-b^{\theta'}\|_C+\|h^\theta-h^{\theta'}\|_{C^2}+\|Y-\tilde Y\|_{C(t,t+T)}}{\epsilon(T,O_Y)},
\end{align*}
which yields \eqref{eq:contdv}.

\subsection{Proof of Lemma \ref{le:bound.supH}}
\label{subse:expH.est}

Equation \eqref{eq:rep1} and the periodicity of $U^{\theta',T}_{f^{\theta'}_0,g^{\theta'}_0}[\psi]$ imply that, for any $\theta'\in\Theta$, $T>0$ and any $\psi\in C(\RR^q)\cap\mathrm{P}_1$,
\begin{align*}
&L^{\theta',\nu^{\theta'}_0}_T \hat{\pi}^{\theta',\nu^{\theta'}_0}_T[E^{\theta'}_{0,T}\psi] = \int_{[0,1)^q} u(0,x)\,\psi^{\theta'}_0(x)\,dx,
\end{align*}
where $u$ is the unique bounded 1-periodic classical solution to
\begin{align}
\pa_t u + \cL^{\theta'} u + (g_0^{\theta'})^\top\nabla u + f_0^{\theta'}\,u = 0,\quad u(T,\cdot)=\psi,\label{eq.appendix.expH.est.eq1}
\end{align}
which is well defined thanks to Lemma \ref{ass:mixing}.




\smallskip

Next, we notice that $\int_{[0,1)^q} \psi^{\theta'}_0(x)\,\cL^{\theta'} \phi(x)\,dx=0$ for any $\phi\in C^2(\RR^q)\cap\mathrm{P}_1$, {and that $|\nabla \psi^{\theta'}_0(x)|\leq C_1\,\psi^{\theta'}_0(x)$ for all $x\in[0,1)^q$ (since $\psi^{\theta'}_0$ is continuously differentiable and bounded away from zero, as stated in Lemma \ref{ass:regular.stationary.density}).
Then,} multiplying \eqref{eq.appendix.expH.est.eq1} by $\psi^{\theta'}_0$ and integrating by parts, we obtain
\begin{align*}
& \left|\pa_t \int_{[0,1)^q} u(t,x)\,\psi^{\theta'}_0(x)\,dx\right| {\leq} \left|\int_{[0,1)^q} \psi^{\theta'}_0(x)\,\cL^{\theta'} u(t,x)\,dx\right| \\
&+ \left|\int_{[0,1)^q} (g_0^{\theta'})^\top(t,x)\nabla u(t,x)\,\psi^{\theta'}_0(x)\,dx\right|
+ \left|\int_{[0,1)^q} f^{\theta'}_0(t,x)\,u(t,x)\,\psi^{\theta'}_0(x)\,dx\right|\\
&\leq C_1\,{\sup_{x\in[0,1)^q}(|g^{\theta'}_0(t,x)|} \int_{[0,1)^q} u(t,x)\,\psi^{\theta'}_0(x)\,dx + \left|\int_{[0,1)^q} u(t,x)\,\psi^{\theta'}_0(x)\,\nabla^\top g^{\theta'}_0(t,x)\,dx\right|\\ 
&\phantom{?????????????????????????????}+ \left|\int_{[0,1)^q} u(t,x)\,\psi^{\theta'}_0(x)\,f^{\theta'}_0(t,x)\,dx\right|\\
&\leq \left(C_1\,{\sup_{x\in[0,1)^q}(|g^{\theta'}_0(t,x)|} + \sup_{x\in[0,1)^q}(|\nabla^{\top} g^{\theta'}_0(t,x)|+|f^{\theta'}_0(t,x)|)\right) \int_{[0,1)^q} u(t,x)\,\psi^{\theta'}_0(x)\,dx\\
&\leq C_2\left(1 + \sup_{t\in[0,T]}|Y_t|\right) \int_{[0,1)^q} u(t,x)\,\psi^{\theta'}_0(x)\,dx,
\end{align*} 
which yields
\begin{align*}
&\int_{[0,1)^q} u(0,x)\,\psi^{\theta'}_0(x)\,dx \geq \int_{[0,1)^q} \psi(x)\,\psi^{\theta'}_0(x)\,dx\, e^{-C_2\,T\,\left(1 + \sup_{t\in[0,T]}|Y_t|\right)}\\
&\phantom{???????????????????????????????}\geq \|\psi\|_{L^1([0,1)^q)}\, e^{-C_3\,\left(1 + T\,\sup_{t\in[0,T]}|Y_t|\right)},\\
&\int_{[0,1)^q} u(0,x)\,\psi^{\theta'}_0(x)\,dx \leq \|\psi\|_{L^1([0,1)^q)}\, e^{C_3\,\left(1 + T\,\sup_{t\in[0,T]}|Y_t|\right)}.
\end{align*}

\smallskip

Thus,
\begin{align*}
&L^{\theta',\nu^{\theta'}_0}_T \hat{\pi}^{\theta',\nu^{\theta'}_0}_T[E^{\theta'}_{0,T}\psi] = \int_{[0,1)^q} u(0,x)\,\psi^{\theta'}_0(x)\,dx\\
&\phantom{???????}\geq \int_{[0,1)^q} \phi(x)\,\psi^{\theta'}_0(x)\,dx\, e^{-C_3\,\left(1 + T\,\sup_{t\in[0,T]}|Y_t|\right)}\\
&\phantom{???????}\geq \|\psi\|_{L^1([0,1)^q)}\, e^{-C_3\,\left(1 + T\,\sup_{t\in[0,T]}|Y_t|\right)}
\geq \|E^{\theta'}_{0,T}\psi\|_{L^1([0,1)^q)}\, e^{-C_4\,\left(1 + T\,\sup_{t\in[0,T]}|Y_t|\right)},\\
&L^{\theta',\nu^{\theta'}_0}_T \hat{\pi}^{\theta',\nu^{\theta'}_0}_T[E^{\theta'}_{0,T}\psi] \leq \|E^{\theta'}_{0,T}\psi\|_{L^1([0,1)^q)}\, e^{C_4\,\left(1 + T\,\sup_{t\in[0,T]}|Y_t|\right)}.
\end{align*}
The above yields
\begin{align*}
&e^{-C_4\,\left(1 + T\,\sup_{t\in[0,T]}|Y_t|\right)}/L^{\theta',\nu^{\theta'}_0}_T\leq
\hat{p}^{\theta',\nu^{\theta'}_0}_T(x)\leq e^{C_4\,\left(1 + T\,\sup_{t\in[0,T]}|Y_t|\right)}/L^{\theta',\nu^{\theta'}_0}_T\quad \text{a.e. }x\in[0,1)^q,
\end{align*}
$\PP^{\theta,\nu^\theta_0}$-a.s., for any $T>0$ and $\theta'\in\Theta$.
Recalling
\begin{align*}
2\frac{H(\nu^{\theta'}_0,\hat{\pi}^{\theta',\nu^{\theta'}_0}_{T})}{1-H(\nu^{\theta'}_0,\hat{\pi}^{\theta',\nu^{\theta'}_0}_{T})}
= \frac{\sup_{x\in[0,1)^q}\frac{\hat{p}^{\theta',\nu^{\theta'}_0}_T(x)}{\psi^{\theta'}_0(x)}}{\inf_{x\in[0,1)^q}\frac{\hat{p}^{\theta',\nu^{\theta'}_0}_T(x)}{\psi^{\theta'}_0(x)}} - 1,
\end{align*}
we conclude the proof of Lemma \ref{le:bound.supH}.

\subsection{Exponential stability of the signal}
\label{subse:expStab.X}

The coupling result stated below is essentially known, and its proof follows a standard methodology. However, we could not find a precise statement of this result that would apply to diffusions on a (flat) torus, hence, we provide its proof herein.

\begin{lemma}\label{le:expStab.X}
For any $\theta\in\Theta$ and $\nu',\nu\in\mathcal{P}([0,1)^q)$, there exists a probability space $(\tilde\Omega,\tilde{\mathcal{F}},\tilde\PP)$ supporting random elements $\tilde X,\,\bar X:\,\tilde\Omega\rightarrow C(\RR_+\rightarrow\RR^q)$, such that
\begin{align*}
& \tilde\PP\circ\tilde X^{-1}=\PP^{\theta,\nu}_X,\quad \tilde\PP\circ\bar X^{-1}=\PP^{\theta,\nu'}_X,
\quad \lim_{T\rightarrow\infty}\tilde\EE \int_0^T |h^\theta(\tilde X_s)-h^\theta(\bar X_s)|\,ds<\infty,
\end{align*}
with the convergence being uniform over all $\nu',\nu\in\mathcal{P}([0,1)^q)$.
\end{lemma}
\begin{proof}
The main idea of the proof dates back to \cite{Rogers}. In short, we construct two solutions $\tilde X$ and $\bar X$ of the SDE \eqref{eq.Background.dyn.X}, started from the initial distributions $\nu'$ and $\nu$ and driven by two Brownian motions that are coupled in such a way that the infinitesimal increments of the diffusion components of $\tilde X$ and $\bar X$ satisfy the following two properties: (i) their projections on the line connecting $\tilde X_t$ and $\bar X_t$ sum up to zero, and (ii) their orthogonal components are equal. Intuitively, these properties ensure that the diffusion component of $\tilde X-\bar X$ ``is one-dimensional" and, hence, will hit zero with probability one. The fact that $\tilde X-\bar X$ evolves on a flat torus, then, ensures that the latter hitting time has a finite expectation. The remainder of the proof is devoted to making this idea rigorous.

\smallskip

Let us fix $\theta\in\Theta$, $\nu',\nu\in\mathcal{P}([0,1)^q)$ and suppress the dependence on $\theta$ throughout the proof, as $\theta$ does not vary.
Next, we notice that the strict ellipticity of $a$, stated in Lemma \ref{ass:mixing}, implies that $d\geq q$ and that the rows of $\sigma$ are linearly independent.
To ease the presentation, we assume $d=q$ and deduce from Lemma \ref{ass:mixing} that $\sigma$ has a bounded inverse.

Then, we consider a weak solution $(\tilde X,\bar X)$ to the following SDE:
\begin{align*}
& d\tilde X_t = b(\tilde X_t)\,dt + \sum_{i=1}^{q-1} \sigma(\tilde X_t)\,v_i(\sigma^{-1}(\tilde X_t)(\tilde X - \bar X_t - l_0))\,\bone_{\{t<\tau\}}\,d\tilde B^i_t\\ &\phantom{???????????????}
+ \frac{\tilde X - \bar X_t - l_0}{\|\sigma^{-1}(\tilde X_t)\,(\tilde X - \bar X_t - l_0)\|}\,\bone_{\{t<\tau\}} \,d\bar B_t
+ \,\bone_{\{t\geq\tau\}}\,\sigma(\tilde X_t)d\hat B_t,\quad \tilde X_0\sim\nu,\\
& d\bar X_t = b(\bar X_t)\,\bone_{\{t<\tau\}}\,dt + \sum_{i=1}^{q-1} \sigma(\bar X_t)\,v_i(\sigma^{-1}(\bar X_t)(\tilde X - \bar X_t - l_0))\,\bone_{\{t<\tau\}} \,d\tilde B^i_t\\ 
&- \frac{\tilde X - \bar X_t - l_0}{\|\sigma^{-1}(\bar X_t)\,(\tilde X - \bar X_t - l_0)\|}\,\bone_{\{t<\tau\}}\,d\bar B_t
+ b(\tilde X_t)\,\bone_{\{t\geq\tau\}}\,dt
+ \,\bone_{\{t\geq\tau\}}\,\sigma(\tilde X_t)d\hat B_t,\quad \bar X_0\sim\nu',
\end{align*}
defined on a probability space $(\tilde\Omega,\tilde{\mathcal{F}},\tilde\PP)$, supporting independent standard Brownian motions $\hat B$, $\tilde B$ and $\bar B$, taking their respective values in $\RR^q$, $\RR^{q-1}$ and $\RR$. In the above, $(\tilde X_0,\bar X_0)$ are mutually independent and independent of $(\tilde B,\bar B)$, the integer vector $l_0\in\ZZ^q$ is defined as a minimizer of
\begin{align*}
& \min_{l\in\ZZ^q}|\tilde X_0-\bar X_0-l|,
\end{align*} 
the random coupling time $\tau$ is defined as
\begin{align*}
\tau:=\inf\{t\geq0:\,\min_{l\in\ZZ^q}|\tilde X_t-\bar X_t-l|=0\},
\end{align*} 
and $\{v_i(y)\}_{i=1}^{q-1}$ are Lipschitz functions of $y\neq0$ that form an orthonormal basis of the orthogonal complement of $y$. 

Under Assumption \ref{holder}, thanks to \cite[Theorem 7.2.1]{stroock1997multidimensional}, there exists a weak solution to the above SDE, defined up to $\tau$, that is unique in law and satisfies the strong Markov property.
It is easy to see that $(\tilde X_\cdot,\bar X_\cdot)$ have a left limit at $\tau$ and that $\tilde X_\tau - \bar X_\tau \in \ZZ^q$.
Then, we extend $\tilde X$ to the time interval $[\tau,\infty)$ via Girsanov's theorem (using the Lipschitz property of $\sigma$) and, recalling $\tilde X_\tau - \bar X_\tau \in \ZZ^q$ and using the 1-periodicity of $(b,\sigma)$, we conclude that $\tilde X_t - \bar X_t = \tilde X_\tau - \bar X_\tau \in \ZZ^q$ and
\begin{align*}
& d\bar X_t = b(\bar X_t)\,dt + \sigma(\bar X_t)d\hat B_t,
\end{align*}
for $t\geq\tau$.
In addition, it is easy to see that
\begin{align*}
\int_0^\cdot \sigma^{-1}(\tilde X_t)\,\left(\sum_{i=1}^{q-1} \sigma(\tilde X_t)\,v_i(\sigma^{-1}(\tilde X_t)(\tilde X - \bar X_t - l_0)) \,d\tilde B^i_t 
+ \frac{\tilde X - \bar X_t - l_0}{\|\sigma^{-1}(\tilde X_t)\,(\tilde X - \bar X_t - l_0)\|} \,d\bar B_t\right)
\end{align*}
is a standard Brownian motion in $\RR^q$, on the time interval $[0,\tau)$, and the analogous conclusion applies to the driving martingale of $\bar X$. Hence, due to the uniqueness in law of the solution to \eqref{eq.Background.dyn.X}, we conclude that
\begin{align*}
& \tilde\PP\circ\tilde X^{-1}=\PP^{\theta,\nu}_X,\quad \tilde\PP\circ\bar X^{-1}=\PP^{\theta,\nu'}_X.
\end{align*}

\smallskip

It only remains to prove the convergence stated in the lemma. To this end, we introduce $Z^0_t:=\tilde X_t - \bar X_t - l_0$, for $t\in[0,\tau\wedge1)$, and observe
\begin{align*}
& dZ^0_t = (b(\tilde X_t) - b(\bar X_t)) dt + \tilde\sigma^0(\tilde X_t,\bar X_t)\,d\tilde B_t\\
&\phantom{??????????????????????}
+ (\|\sigma^{-1}(\tilde X_t)\,Z^0_t/|Z^0_t|\|^{-1} + \|\sigma^{-1}(\bar X_t)\,Z^0_t/|Z^0_t|\|^{-1})\,\frac{Z^0_t}{|Z^0_t|}\,d\bar B_t,
\end{align*}
where $\tilde\sigma^0$ is uniformly absolutely bounded and satisfies $\|\tilde\sigma^0(\tilde x,\bar x)\|\leq C_1 |\tilde x - \bar x-l_0|$ for all $\tilde x\neq\bar x+l_0$, with a constant $C_1$.
Using It\^o's formula, we obtain
\begin{align*}
& d|Z^0_t|^2 = 2(Z^0_t)^\top(b(\tilde X_t) - b(\bar X_t)) dt + 2(Z^0_t)^\top\tilde\sigma^0(\tilde X_t,\bar X_t)\,d\tilde B_t\\
&+ 2|Z^0_t|\,\bar\sigma(\tilde X_t,\bar X_t)\,d\bar B_t + \mathrm{Tr}\left[ \tilde\sigma^0(\tilde\sigma^0)^\top(\tilde X_t,\bar X_t)\right]\,dt
+ \bar\sigma^2(\tilde X_t,\bar X_t)\,dt,
\end{align*}
where $1/C_2\leq \bar\sigma(\tilde x,\bar x)\leq C_2$ for all $\tilde x\neq\bar x$, with a constant $C_2>1$.
Then, we continue:
\begin{align*}
& d|Z^0_t| = \frac{(Z^0_t)^\top}{|Z^0_t|}(b(\tilde X_t) - b(\bar X_t)) dt 
+ \frac{1}{2|Z^0_t|}\mathrm{Tr}\left[ \tilde\sigma^0(\tilde\sigma^0)^\top(\tilde X_t,\bar X_t)\right]\,dt
+ \frac{1}{2|Z^0_t|}\bar\sigma^2(\tilde X_t,\bar X_t)\,dt\\
& + \frac{(Z^0_t)^\top}{|Z^0_t|}\tilde\sigma^0(\tilde X_t,\bar X_t)\,d\tilde B_t
+ \bar\sigma(\tilde X_t,\bar X_t)\,d\bar B_t\\
&- \frac{1}{8|Z^0_t|^3}\left[4(Z^0_t)^\top\tilde\sigma^0(\tilde\sigma^0)^\top(\tilde X_t,\bar X_t)\,Z^0_t 
+ 4|Z^0_t|^2\,\bar\sigma^2(\tilde X_t,\bar X_t) \right]\,dt\\
& = \frac{(Z^0_t)^\top}{|Z^0_t|}(b(\tilde X_t) - b(\bar X_t))\,dt 
+ \frac{\mathrm{Tr}\left[ \tilde\sigma^0(\tilde\sigma^0)^\top(\tilde X_t,\bar X_t)\right]}{2|Z^0_t|}\,dt
- \frac{(Z^0_t)^\top\tilde\sigma^0(\tilde\sigma^0)^\top(\tilde X_t,\bar X_t)\,Z^0_t}{2|Z^0_t|^3}\,dt\\
& + \frac{(Z^0_t)^\top}{|Z^0_t|}\tilde\sigma^0(\tilde X_t,\bar X_t)\,d\tilde B_t
+ \bar\sigma(\tilde X_t,\bar X_t)\,d\bar B_t.
\end{align*}
Note that there exist constants $C_3>0$ and $C_4>1$ such that
\begin{align*}
& \left| \frac{(Z^0_t)^\top}{|Z^0_t|}(b(\tilde X_t) - b(\bar X_t)) 
+ \frac{\mathrm{Tr}\left[ \tilde\sigma\tilde\sigma^\top(\tilde X_t,\bar X_t)\right]}{2|Z^0_t|}
- \frac{(Z^0_t)^\top\tilde\sigma^0(\tilde\sigma^0)^\top(\tilde X_t,\bar X_t)\,Z^0_t}{2|Z^0_t|^3} \right| \leq C_3,\\
& 1/C_4 \leq |Z^0_t|^{-2} |(Z^0_t)^\top\tilde\sigma^0(\tilde X_t,\bar X_t)|^2
+ \bar\sigma^2(\tilde X_t,\bar X_t) \leq C_4,
\end{align*}
$\tilde\PP$-a.s. for all $t\in[0,\tau\wedge1)$.
Noticing that $|Z^0_0|$ is bounded from above by a universal constant depending only on $q$, we deduce from the above the existence of a left limit $Z^0_{\tau-}$, which we identify with $Z^0_\tau$, and of a constant $\varepsilon>0$, such that
\begin{align*}
\tilde\PP\left(\inf_{t\in[0,\tau\wedge1]}Z^0_t \leq 0\,\vert\,\tilde X_0,\,\bar X_0\right)\geq \varepsilon
\quad\tilde\PP\text{-a.s..}
\end{align*}

Next, we construct $\{Z^i\}_{i=1}^\infty$ recursively, as follows. If $\tau\leq i$, then we set $Z^{i+1}:=Z^{i}$. Otherwise, we define $l_i\in\ZZ^q$ as a minimizer of
\begin{align*}
& \min_{l\in\ZZ^q}|\tilde X_i - \bar X_i - l|,
\end{align*}
and set $Z^{i+1}_t:=\tilde X_t - \bar X_t - l_0$, for $t\in[i,\tau\wedge(i+1))$.
Repeating the above arguments, we deduce
\begin{align*}
\tilde\PP\left(\inf_{t\in[i,\tau\wedge(i+1)]}Z^i_t \leq 0\,\vert\,\tilde X_i,\,\bar X_i,\right)\,\bone_{\{\tau\geq i\}}\geq \varepsilon\,\bone_{\{\tau\geq i\}}
\quad\tilde\PP\text{-a.s.,}
\end{align*}
for all $i\geq0$.

It is clear that the event $\{\inf_{t\in[i,\tau\wedge(i+1)]}Z^i_t \leq 0,\,\tau\geq i\}$ implies $\{\tau\in[i,i+1]\}$.
Then, using the Markov property of $(\tilde X,\bar X)$, we obtain, for any $n\geq1$:
\begin{align*}
\tilde\PP(\tau\geq n) \leq \tilde\EE\left[\bone_{\{\tau\geq n-1\}}\,\tilde\PP\left(\inf_{t\in[n-1,\tau\wedge n]} Z^{n-1}_t > 0\,\vert\,\tilde X_{n-1},\,\bar X_{n-1}\right)\right]
\leq (1-\varepsilon)\,\tilde\PP(\tau\geq n-1).
\end{align*}
Iterating the above, we deduce the exponential decay of $t\mapsto\tilde\PP(\tau\geq t)$ and
\begin{align*}
&\tilde\EE |h(\tilde X_t)-h(\bar X_t)| \leq C_5\,\tilde\PP\left(\min_{l\in\ZZ^q}|\tilde X_t-\bar X_t-l|>0\right)
= C_5\,\tilde\PP\left(\tau>t\right) \leq C_5\, (1-\varepsilon)^{\lfloor t\rfloor},
\end{align*}
which completes the proof of the lemma.
\end{proof}

\subsection{Mutual singularity of the observation laws}
\label{subse:singularity}

The next lemma proves the following remarkable property: for any two probability measures under which the signal process is stationary, if the laws of the observation process under these measures are not identical, they must be mutually singular.

\begin{lemma}\label{le:singularity}
For any $\theta\in\Theta$ and any $\theta'\notin\overline{\Theta}(\theta)$, the measures $\PP^{\theta,\nu^\theta_0}_Y$ and $\PP^{\theta',\nu^{\theta'}_0}_Y$ are mutually singular.
\end{lemma}
\begin{proof}
Since $\PP^{\theta,\nu^\theta_0}_Y\neq\PP^{\theta',\nu^{\theta'}_0}_Y$, there exist $k\geq1$, $0\leq t_1<\cdots<t_k$, and a continuous bounded $F:\,(\RR^{m})^k\rightarrow\RR$, such that
\begin{align}
\EE^{\theta,\nu^\theta_0}F(Y_{t_1},Y_{t_2}-Y_{t_1},\ldots,Y_{t_k}-Y_{t_{k-1}})
\neq \EE^{\theta',\nu^{\theta'}_0}F(Y_{t_1},Y_{t_2}-Y_{t_1},\ldots,Y_{t_k}-Y_{t_{k-1}}).\label{eq.singularity.eq1}
\end{align}
To ease the notation, we assume that $k=2$ and consider
\begin{align*}
\bar\xi_n:=\frac{1}{n}\sum_{i=1}^n \xi_i,
\quad \xi_i = F(Y_{t_1+t_2i}-Y_{t_2i},Y_{t_2+t_2i}-Y_{t_1+t_2i}).
\end{align*}

Next, we notice that the coupling constructed in the proof of Lemma \ref{le:expStab.X} yields a slightly stronger statement: for any $s>0$ and any bounded measurable $G:\,C([0,s]\rightarrow\RR^{2m})\rightarrow\RR$, there exist constants $C_1,C_2>0$ such that
\begin{align*}
&|\EE^{\theta,\nu^\theta_0} G(h(X_\cdot)_{[t,t+s]},(W_\cdot-W_t)_{[0,s]}) - \EE^{\theta,\delta_x} G(h(X_\cdot)_{[t,t+s]},(W_\cdot-W_t)_{[0,s]})| \leq C_1\,e^{-C_2\,t},
\end{align*}
for all $t\geq 0$, $x\in\RR^q$.
Then, noticing that there exists $G$ such that 
\begin{align*}
& \xi_i = G(h(X_\cdot)_{[t_2(i-1),t_2(i-1)+t_2]},(W_\cdot-W_{t_2(i-1)})_{[0,t_2]}),
\end{align*}
we use the Markov property of $X$ to obtain:
\begin{align*}
&\EE^{\theta,\nu^\theta_0} (\bar\xi_n - \EE^{\theta,\nu^\theta_0}\bar\xi_n)^2
= O(1/n) + \frac{2}{n^2}\sum_{1\leq i< j\leq n} \EE^{\theta,\nu^\theta_0} (\xi_i - \EE^{\theta,\nu^\theta_0}\xi_i)(\xi_j - \EE^{\theta,\nu^\theta_0}\xi_j)\\
&=O(1/n) + \frac{2}{n^2}\sum_{1\leq i< j\leq n} \EE^{\theta,\nu^\theta_0}\left[ (\xi_i - \EE^{\theta,\nu^\theta_0}\xi_i)\,\EE^{\theta,\nu^\theta_0}\left( \xi_j - \EE^{\theta,\nu^\theta_0}\xi_j \,\vert\, \mathcal{F}^{X,Y}_{t_2i}\right) \right]\\
&=O(1/n) + \frac{2}{n^2}\sum_{1\leq i< j\leq n} \EE^{\theta,\nu^\theta_0}\left[ (\xi_i - \EE^{\theta,\nu^\theta_0}\xi_i)\,\left(\EE^{\theta,\delta_x}\xi_j - \EE^{\theta,\nu^\theta_0}\xi_j \right)_{x=X_{t_2i}} \right]\\
&\leq O(1/n) + \frac{2C_1}{n^2}\sum_{1\leq i< j\leq n} \EE^{\theta,\nu^\theta_0}\left|\xi_i - \EE^{\theta,\nu^\theta_0}\xi_i\right|\,e^{-C_2\,t_2\,(j-i)} = O(1/n),
\end{align*}
as $n\rightarrow\infty$. Repeating the above arguments under $\PP^{\theta',\nu^{\theta'}_0}$, we conclude that the same convergence of the variance to zero holds under the latter measure. Since $\{\xi_i\}$ are identically distributed, we conclude that $\bar\xi_n$ converges $\PP^{\theta,\nu^{\theta}_0}_Y$-a.s. along a subsequence to 
\begin{align*}
&\EE^{\theta,\nu^{\theta}_0} \xi_1 = \EE^{\theta,\nu^\theta_0}F(Y_{t_1},Y_{t_2}-Y_{t_1},\ldots,Y_{t_k}-Y_{t_{k-1}}).
\end{align*}
as $n\rightarrow\infty$. Analogous conclusion holds under $\PP^{\theta',\nu^{\theta'}_0}_Y$. Recalling \eqref{eq.singularity.eq1}, we obtain two disjoint events in $\mathcal{F}^Y_\infty=\sigma(\bigcup_{t\geq0}\mathcal{F}^Y_t)$, corresponding to the convergence of $\bar\xi_n$ along a common subsequence to two different limits, each having probability one under the associated measure. This yields the desired mutual singularity of $\PP^{\theta,\nu^\theta_0}_Y$ and $\PP^{\theta',\nu^{\theta'}_0}_Y$.  
\end{proof}

\bigskip\bigskip\bigskip

\bibliographystyle{amsalpha}
\bibliography{ref}

\bigskip\bigskip\bigskip

\end{document}